\documentclass[10pt]{article}

\usepackage{url}
\usepackage{mathtools}
\usepackage{amsmath}
\usepackage{amssymb}
\usepackage{mathabx}
\usepackage{amsthm}
\usepackage{empheq}
\usepackage{latexsym}
\usepackage{enumitem}
\usepackage{eurosym}
\usepackage{dsfont}
\usepackage{appendix}
\usepackage{color} 
\usepackage[unicode]{hyperref}
\usepackage{frcursive}
\usepackage[utf8]{inputenc}
\usepackage[T1]{fontenc}
\usepackage{geometry}
\usepackage{multirow}
\usepackage{todonotes}
\usepackage{lmodern}
\usepackage{anyfontsize}
\usepackage{stmaryrd}
\usepackage{cases}
\usepackage{calc}
\usepackage{pdfpages}
\usepackage{graphicx}
\usepackage{float}

\definecolor{red}{rgb}{0.7,0.15,0.15}
\definecolor{green}{rgb}{0,0.5,0}
\definecolor{blue}{rgb}{0,0,0.7}
\hypersetup{colorlinks, linkcolor={red},citecolor={green}, urlcolor={blue}}
			
\makeatletter \@addtoreset{equation}{section}

\usepackage[justification=centering]{caption}
\newtheorem{theorem}{Theorem}[section]
\newtheorem{assumption}[theorem]{Assumption}

\newtheorem{example}[theorem]{Example}

\newtheorem{lemma}[theorem]{Lemma}
\newtheorem{proposition}[theorem]{Proposition}

\newtheorem{definition}[theorem]{Definition}
\newtheorem{remark}[theorem]{Remark}

\DeclareUnicodeCharacter{014D}{\=o}
\def \N{\mathbb{N}}

\def \R{\mathbb{R}}

\def\Ac{{\cal A}}

\def\Cc{{\cal C}}

\def\Nc{{\cal N}}
\def\Oc{{\cal O}}
\def\Pc{{\cal P}}

\newcommand{\abs}[1]{\left|#1\right|}

\urldef{\myself}\url{https://www.kelwatt.fr/prix/tarif-bleu-edf#resume}

\usepackage{tikz}
\usepackage{pgfplots}
\pgfplotsset{compat=1.14}
\usetikzlibrary{positioning}
\usetikzlibrary{decorations.text}
\usetikzlibrary{calc,intersections}
\usepackage{tkz-euclide,subfigure}
\usetkzobj{all}    
\usetikzlibrary{calc}

\title{An adverse selection approach to power pricing\footnote{The authors gratefully acknowledge the support of the ANR project Pacman, ANR-16-CE05-0027.}}
 \author{Cl\'emence {\sc Alasseur} \thanks{EDF R\&D and FIME, Laboratoire de Finance des March\'es de l'\'Energie, clemence.alasseur@edf.fr.}\and Ivar {\sc Ekeland} \thanks{Universit\'e Paris--Dauphine, PSL Research University, CNRS, CEREMADE and Institut de Finance, ekeland@math.ubc.ca.}\and Romuald {\sc \'Elie} \thanks{Universit\'e Paris-Est Marne-la-Vall\'ee, romuald.elie@univ-mlv.fr
.}\and Nicol\'as {\sc Hern\'andez Santib\'a\~{n}ez} \thanks{Department of Mathematics, University of Michigan, nihernan@umich.edu.}\and Dylan {\sc Possama\"{i}} \footnote{Columbia University, IEOR department, dp2917@columbia.edu.}}

\date{\today}

\begin{document}

\maketitle

\begin{abstract}

We study the optimal design of electricity contracts among a population of consumers with different needs. This question is tackled within the framework of Principal--Agent problems in presence of adverse selection. The particular features of electricity induce an unusual structure on the production cost, with no decreasing return to scale. We are nevertheless able to provide an explicit solution for the problem at hand. The optimal contracts are either linear or polynomial with respect to the consumption. Whenever the outside options offered by competitors are not uniform among the different type of consumers, we exhibit situations where the electricity provider should contract with consumers with either low or high appetite for electricity. 
\vspace{5mm}


\vspace{5mm}


\end{abstract}

\section{Introduction}
Electricity is non--storable, except marginally: any quantity which is consumed now must be produced now, and conversely. This means that the installed capacity must be sufficient to supply electricity when demand is maximal. As a consequence, part of this installed capacity will stay idle when demand is lower. This is the overcapacity problem, which is compounded by the fact that electricity consumption is far from stable, and little substitutable. 
Since there is no electricity stored to dampen shocks and smooth discrepancy, adjusting supply to demand is a difficult task. One way to do this is to use prices. Very early on, power companies have hit upon the idea of making electricity more expensive in peak hours, so that consumers who are able to do so would switch their demand to off--peak periods. This falls naturally within the framework of Principal--Agent problems: the Principal (here the power company) offers a tariff that provides incentives, and each Agent (consumer) reacts according to his own needs. 
It does not seem, however, that such an analysis is available at the present time, and the aim of the present paper is to fill this gap.

\medskip
We focus on the problem of finding an appropriate tariff: 
for a given production function, how should a power company price electricity in order to maximise its profit? The company faces a variety of consumers, industrial users and domestic users, some of them are efficient, others less so. Some of them, for instance, live in insulated homes and need less electricity to achieve comfortable temperatures than others. The tariff the producer offers will be time--dependent and consumption--based. It will act in several ways: redirecting part of the consumers to off-peak periods, by pricing properly the peak hours; avoiding overly expensive production costs, by penalising higher consumptions; effectively excluding some of the users, who will find the proposed tariff too expensive, and who will look for better alternatives elsewhere. The empirical effective response of individual agents to hourly pricing contracts is studied in \cite{allcott2011rethinking} or \cite{faruqui2010household} and is heterogeneous among the electricity consumers, as pointed out in \cite{herter2007residential} or \cite{jessoe2014towards}. Besides, more resilient measures of high electricity prices appear to be more efficient in practice for incentivising the reduction of electricity consumption, see \cite{ito2014consumers} or \cite{wolak2011residential}. In any case, the critical importance of the feedback information signals provided to the Agent has been highlighted in the empirical literature \cite{faruqui2013dynamic, gans2013smart, gilbert2014dynamic}. In particular, the development of smart metering systems allows for the implementation of decentralised home automation systems, who partially regulate power consumption, as presented e.g. in \cite{abras2008multi} or \cite{lagorse2010multi}. 
 %
 %
The adoption of dynamic electricity tariffs is also proven to be heterogeneous among the population \cite{kowalska2014turning,qiu2017risk, harding2014goal}. In our model, we encompass this feature by considering that each consumer has different private characteristics that are summed up in a type. The type in the population plays two roles and it leads to an adverse selection problem for the company. First, consumers with high type (most efficient) have a higher utility for a given level of consumption than consumers with low type (least efficient). Second, consumers have outside options to the contract, such as swapping for a competitor or to alternative sources of energy, which provide them reservation utilities. The reservation utility is not necessarily the same for everybody and we model it as a function of the type. The heterogeneity of fallback options also emphasizes the effect of competitors on the energy retailer market, whose impact is discussed in  \cite{joskow2006retail} or \cite{spear2003electricity}.

 \medskip
We frame the problem in a Principal--Agent model, in which the power company proposes a contract, that is, a tariff, and the consumers either turn down the contract and drop out, or accept it, and adjust their consumption accordingly. Electricity pricing has a special feature, which distinguishes it from other Principal--Agent problems. Usually, the profit of the Principal is the sum of the profits she gets from all participating Agents. Here, the cost to the power plant is the cost of producing the aggregate demand, which is not the sum of the costs of producing the individual demands, because of decreasing returns to scale in production. This introduces a mathematical difficulty which 
is not treated frequently in the literature and has deep consequences on the tariff.
There are actually earlier works in a bilevel optimisation setting, where only linear pricing is considered, see for example \cite{hobbs1992nonlinear} where the question is to provide electricity and to sponsor at the same time saving measures, or \cite{afsar2016achieving} where consumers can reduce their consumption at a price of inconvenience. As a related study, we mention as well the Principal--Agent modelling discussed in \cite{gillingham2012split} for reducing energy consumption in a landlord--tenant relationship.

\medskip
Instead of putting an upper bound on production, we consider a production function with steeply increasing marginal cost, as discussed for example in \cite{dierker1976increasing}. In reality, the capacity constraint is not binding: more electricity can always be found, by putting in service less efficient production units, or by resorting to the spot market. It only becomes extremely expensive when the limits are pushed.\footnote{Some power plants have no or very low fuel costs such as renewables (hydro power plants, wind or solar production) or nuclear production. These types of productions are chosen for satisfying base--load consumption. But when, the electricity consumption increases such as in peak hours, other power plants (coal, gas or fuel thermal plants for example) need to be turned on and their cost of production is much more expensive.} The production function can be understood either as reflecting the actual cost of producing electricity from primary energies, in which case the company is a producer, or the financial cost of buying electricity on the open market, in which case the company is a retailer, or a combination of both. 

\medskip
Increasing marginal costs are considered in \cite{spulber1993monopoly}, where in a similar model to ours, the profit--maximising direct revelation mechanism is found, along with different pricing strategies implementing it. In our work, however, we do not aim at finding a menu of tariffs, but instead a single one which is offered to all the clients. For this reason our methodology, which is described below, departs from the usual analysis of multi--dimensional screening, such as \cite{armstrong1999multi,rochet1985taxation,rochet1987necessary,rochet1998ironing}. We are able to solve the problem at hand in the present paper explicitly, in some particular cases (CRRA utility function, power cost). The tariff we find is quite natural. It consists in the sum of three
components:

\begin{itemize}[leftmargin=*]
	\item a fixed component, independent of the consumption, which is a
	subscription to the service;
	
	\item a linear component in consumption, which
	consists simply of pricing the current consumption at the current price. Recall that the price depends on time in order to discourage peak--hour consumption;
	
	\item a non--linear component on the current consumption. This part appears only when the consumption is high enough. Its purpose is to make the high consumers pay
	for making the Principal produce at very high marginal cost. 
\end{itemize}


In Principal--Agent models, one is used to the fact that the optimal contract for the Principal will exclude the low end of the market, that is, the least
efficient Agents. Indeed, such {\it shutdown} contracts are quite commonplace in the adverse selection literature, see for instance the seminal papers of Guesnerie and Laffont \cite{guesnerie1984complete} or Rochet and Chon\'e \cite{rochet1998ironing} in discrete--time (see also the monograph by Laffont and Martimort for more details \cite{laffont2009theory}), or the more recent contributions of Cvitani\'c et al. \cite{cvitanic2013dynamics} and Hern\'andez Santib\'a\~nez et al. \cite{hernandez2017bank} in continuous--time. A remarkable feature of our problem is that, in certain
circumstances, the optimal contract excludes the high end of the market: the
most efficient agents, those who need less electricity to achieve the same
degree of welfare than others, go elsewhere, and only the least efficient
ones remain. 
In that case it becomes worthwhile for the company not to take part
in the competition for the most efficient Agents, because matching their fall back options is too costly, 
and it is preferable to concentrate on the least efficient consumers. 
Exclusion of agents other than the least efficient ones can be found in models with countervailing incentives (see \cite{lewis1989countervailing}). In \cite{jullien2000participation} and \cite{maggi1995countervailing}, it may be optimal to exclude intermediate agents from the contract, which is also the case in our model.

\medskip
We now proceed to describe the main features of our model. The Principal's
cost, as mentioned above, is a convex function of aggregate production. She
offers a contract to the Agents, who may accept or decline it. If they
accept, they commit for a period $T>0$. They decline if the total utility
they derive from the contract is less than the reservation utility corresponding to their type. The Agents' utilities are separable: the utility which an agent of type $x$
derives from consuming a quantity $c$ of electricity at time $t$ and being
charged a (nonlinear) price $p$ is $u( t,x,c) -p(t,c),$ where $u( t,x,\cdot) $ is a concave function of $c$. This
separability assumption is traditional in Principal--Agent problems. In this
case, there are additional justifications, as a large part of the Agents are
industrial users, who consume electricity in order to produce other goods,
so that their utility is simply the profit they derive from this
activity. Note also the time--dependence, which reflects the seasonality of
consumption.

\medskip
In the sequel, we will consider CRRA utilities, of the type $\gamma
^{-1}c^{\gamma }$, with $\gamma <1$, and we will provide explicit solutions
(except in the case $\gamma =0$, or $u( c) =\ln c$). The case $0<\gamma <1$
reflects the "industrial" behaviour, where high consumption is the norm,
subject to decreasing returns to scale. The case $
\gamma <0$ reflects the "household" behaviour, where electricity fulfils
some basic needs, such as lighting or appliances, and $0$ consumption is not
acceptable, while high consumption is not needed\footnote{In this case the utilities are negative but the analysis is not impacted. One can even make them positive by just adding a constant.}. Note, however, that in both cases
there is a "fallback" option, a substitute to electricity when it becomes
too expensive, for instance an alternative energy source, or simply another
provider. 
Despite the particular structure of
the cost function, we are able to solve explicitly the problem at hand. We
observe that the optimal contract rewrites as the combination of a fixed
cost together with two variable costs, proportional to either the
electricity consumption or a power function of it. This tariff structure
happens to be quite simple and quite close to the classical tariff
structures offered by most electricity providers.

\medskip
Whenever the fallback option is the same for every Agent, we observe as
usual in Principal--Agent problems, that the lower end of the market is not
covered: the low types (meaning those households who are less dependent on
electricity, or those industry users who are less efficient) will not be
offered contracts which they would accept, and will have to fall
back on the outside option. More interestingly, we are also able to solve
explicitly the case where the fallback option of the Agents depends on their
type in a concave manner. In this case, getting more efficient Agents can be
too costly, and the electricity provider may concentrate on the less
efficient but less expensive consumers. 

\medskip
Finally, a remark on the mathematics. We will be using $u-$convex analysis,
a tool introduced and developed elsewhere, notably by Carlier (see \cite{carlier2001general}), and which extends
classical convex analysis. 
Forgetting for ease of presentation about time
dependence, recall that the maximal utility Agent $x$ can obtain when a tariff $p(c)$ is set, is equal to
$\max_{c}\left\{ u( x,c) -p( c) \right\}.$ This maximum is denoted by $p^{\star }( x) $, 
note that it depends on the entire price
schedule and it can be computed for each Agent $x$. In this way, we
associate with each function $p( c) $ a function $p^{\star }(
x) $, which the economist knows as the indirect utility associated
with $p$ and which the mathematician knows as the $u-$transform of $p$. Conversely, if the indirect utility $p^{\star}$ is
known, the price schedule can be derived by the same formula $p( c) =\max_{x}\left\{ u( x,c) -p^{\star }(
x) \right\}.$ In the bilinear case, when $u( x,c) =xc$, we get the usual
Fenchel formulas of convex analysis.

\medskip
The rest of the paper is organised as follows. Section \ref{sec:main} sets the model and the main results, i.e. the expression of the tariffs for industrial and residential customers. In Section \ref{sec:interpret} we provide economic interpretation of numerical results. In Section \ref{sec:model}, we provide a rigorous definition of the model, from the mathematical point of view. Sections \ref{section-constant} and \ref{sec:general} provide the main results for constant and concave increasing reservation utility; proofs are left for appendixes. Finally, the conclusions are given in Section \ref{sec:conc}.

\section{Main results}\label{sec:main}

The model we propose is set up on Principal--Agent relationship where the Principal is an electricity provider and the Agents are consumers. Since the electricity consumption is observed by the Principal, there is no moral hazard. On the other hand, adverse selection is in force since the Agents' willingness to pay for electricity is not known by the Principal. This taste for electricity represents how much Agents value a given volume of electricity in terms of usefulness. For an industrial Agent, this would represent the benefit he gets by running his industrial process with this given volume of electricity. For a residential Agent, this would represent the comfort he gets by using this given volume of electricity to perform domestic tasks. 
Of course, this depends on the Agents' equipment, referred to as his type $x$. As classically assumed in adverse selection setting, even if the Principal does not know the exact type of a particular Agent, he knows the proportion of Agents' type among the population. This hypothesis is realistic as the electricity provider can make pre--marketing surveys, or can use historical data in order to acquire such information.

\subsection{Players' objectives and electricity particularity}
Both players have their particular objectives
\begin{itemize}[leftmargin=*]
	\item Agent's objective is to choose the level of electricity consumption $c$ at any time $t$, which maximizes his utility for electricity $u\left(t,x,c\right)$ with respect to his type $x$, minus the tariff $p(t,c)$ that he needs to pay for the electricity
	\[
	\max_{c}\bigg\{ \int_0^T u\left(t,x,c\right) -p\left(t,c\right) \mathrm{d}t \bigg\}. 
	\]
	\item Principal's objective is to offer the tariffs which maximise his own profits: all payments she receives from consumers accepting the contract minus the costs for providing the total volume of electricity consumed by her clients. 

\end{itemize}

One particular feature of electricity production, is the fact that it suffers from decreasing returns to scale: its marginal price increases with the total aggregate consumption. We consider therefore a convex cost function for the Principal, as discussed in the Introduction.

\subsection{Notations and model assumptions}

\begin{figure}[!ht]
	\begin{center}
	\hspace*{-6mm}
	\begin{tabular}{cc}
		$\gamma \in (0,1)$ &  $\gamma <0$  \\
		\includegraphics[scale=0.25]{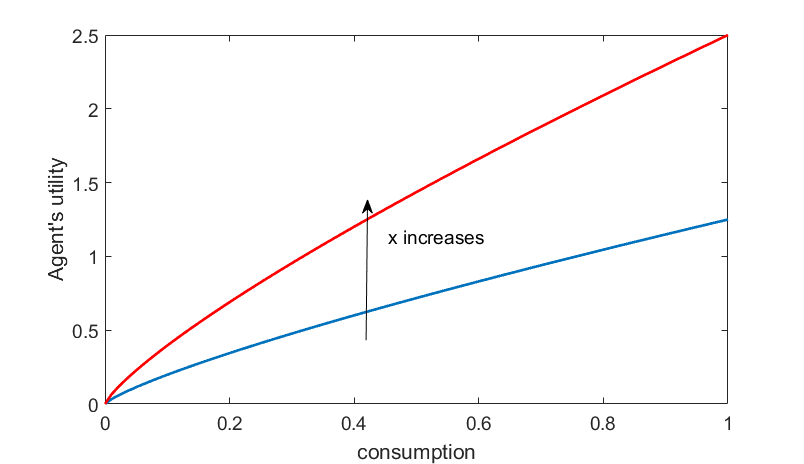} 	
		&
		\includegraphics[scale=0.25]{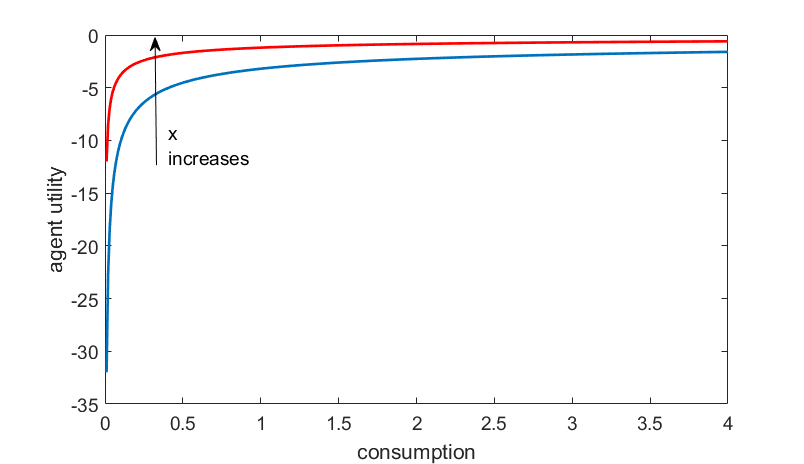} 	
		
	\end{tabular} 
	\caption{\label{utility} Agent's utility with respect to consumption for $\gamma \in (0,1)$ (left figure) and $\gamma <0$ (right figure).}
	\end{center}
	
\end{figure}

We consider constant relative risk aversion (CRRA) utility functions for the Agents
$$
u(t,x,c) = g_\gamma(x)\phi(t)\frac{c^\gamma}{\gamma},
$$
where $\phi(t)$ represents the Agents' time preference for electricity. This factor is common to every Agent and typically represents the preference to have electricity during the daytime than during the middle of the night. We suppose that $\gamma<1$  and we consider two different cases: $\gamma\in (0,1)$ and $\gamma<0$. The function $g_\gamma$ represents the willingness of the Agents to pay for their consumption depending on their type $x$, and we take typically $g_\gamma(x)=x$, if $\gamma\in (0,1)$ and $g_\gamma(x)=1-x$, if $ \gamma<0$. Graphic illustrations of the utility function are shown in Figure \ref{utility}. The case $\gamma \in (0,1)$ corresponds to the modelling of industrial Agents, whose utility grows to infinity if they can have infinite volume of electricity: they can always make their industrial capacities grow and generate more benefits whenever they have extra electricity. On the contrary, they can stop producing if they could not get any electricity or substitute it by another energy, which corresponds to a zero utility whenever $c=0$. The case $\gamma <0$ illustrates the residential Agents utility for whom electricity is a staple product: they can not avoid consuming electricity (they would get $-\infty$ utility). They also face a saturation effect: above a high volume of electricity, they do not gain much satisfaction with an extra quantity, because all their electrical needs are already fulfilled.

The density function of the types is known by the Principal and denoted by $f$. The indirect utility $P^\star(x)$ is the best level of utility that an Agent of type $x$ can obtain by signing the contract, during the period $[0,T]$
$$
P^\star(x):= \int_{0}^{T}p^\star(t,x) \mathrm{d}t = \underset{c}{\sup}\bigg\{\int_{0}^{T} u(t,x,c)-p(t,c)  \mathrm{d}t\bigg\}.
$$
The tariffs designed by the Principal need to respect the individual rationality of the Agents. Indeed, Agents are not forced to accept the contract offered by the Principal, since they can pick alternative electricity providers, offering better conditions. This is taken into account in the model via a reservation utility $H$ which represents the minimum level of satisfaction that an Agent needs to achieve in order to accept the contract. 

\medskip
Agent of type $x$ signs the contract with the Principal if and only if $P^\star(x)\ge H(x)$. We consider two cases for $H(x)$, either a constant function or an increasing concave one verifying a certain condition.\footnote{Namely, that the elasticity of reservation utility is smaller than the elasticity of willingness to pay for consumption, i.e. $ g_\gamma/g_\gamma^\prime \le H/H^\prime$.} The increasing property of $H$ indicates that competitors target principally the more efficient Agents. Indeed, 
a high value of $H$ means that the customer is harder to satisfy because he has better external alternatives, which we suppose to be competitors offers. 

\medskip
Finally, we denote by $X^\star$ the set of Agents who end up signing the contract 
$$X^\star(p):  = \left\{ x\in [0,1],\ P^\star (x) \geq H(x) \right\}.$$
As mentioned before, the cost of production depends on the set $X^\star(p)$.\footnote{In general the production costs depends on the consumption of all the population, so we are 
 implicitly assuming that aggregate consumption of clients who select the provider is correlated to the consumption of all the other consumers. 
This is justified since there are strong common preferences and behaviours among consumers, for example consumptions are higher during daytime than during the night.} In some explicit examples, we consider a convex cost of power production $
 K(t,\tilde{c}) = k(t)\frac{\tilde{c}^n}{n},
$ where ${\tilde{c}}$ refers to the total consumption of clients. The term $k(t)$ is positive and indicates the time dependence of electricity production costs (for example photovoltaic production occurs only at day and wind is blowing more during winter). The power $n>1$ reflects the production fleet composition; the fleet has expensive peak power plants when $n$ is high. 

\subsection{Optimal tariffs}
In the setting we previously described, the Principal--Agent problem can be explicitly solved. We describe in this part the form of the solution, and leave the rigorous mathematical proofs to the remaining sections of the paper. 
In order to be admissible, a tariff $p$ should verify the individual rationality and incentive compatibility conditions. This is denoted by $p\in\mathcal{P}$. Let us write formally the objective function of the Agents $U_A$ and Principal $U_P$
$$
	U_A(p,x) :=\sup_{c} \int_{0}^{T} \left( u(t,x,c(t)) - p(t,c(t))\right)\mathrm{d}t= \int_{0}^{T}  p^\star(t,x)\mathrm{d}t.
$$ 
The solution to this optimization problem is the optimal consumption of the Agents of type $x$ and it is denoted by $c^\star(\cdot,x)$. The map $c^\star$ becomes part of the problem of the Principal, as follows
$$	
U_P :=\sup_{p\in\mathcal{P}} \int_{0}^{T} \bigg[ \int_{X^\star(p^\star) } p(t,c^\star(t,x)) f(x) \mathrm{d}x  - K\bigg(t,\int_{X^\star (p^\star)} c^\star(t,x) f(x)\mathrm{d}x\bigg)  \bigg] \mathrm{d}t.
$$
We compute that, whatever $\gamma$ or $H$ (constant or concave), the optimal tariff is a function of three components at most, namely a constant part $p_3$, a proportional part $p_2$ of the consumed power $c$, and a proportional part $p_1$ of $c^\gamma$
$$
p(t,c)=p_1(t)c^\gamma + p_2(t)c + p_3(t).
$$
This tariff is always a concave increasing function of the consumed power $c$. An important observation is that this tariff is quite simple and close to current tariff structures proposed by electricity providers. Indeed, they are commonly split into a fixed charge in \euro, a volumetric charge in \euro/MWh, and possibly a demand charge in \euro/MW. The fixed and the volumetric charges can depend on the maximum subscribed power which is another way to price the demand charge. The optimal tariffs offered by the Principal are summarised in the following tables, where the explicit expressions for the functions  $(p_{i,\gamma})_i$, $(p_i^j)_{i,j}$ and $\hat{c}_i^\gamma(t)$ are respectively provided in Theorem \ref{th:mainex} and Theorem \ref{t:exam2} hereafter.

\begin{table}[ht!]
\small
\begin{center}
	\begin{tabular}{|c|c|c|}
		\hline
		\text{Selected Agents}&  \text{$H$ constant} &     \text{$H$ concave non decreasing}  \\
		\hline
		$\left[a_0,1\right]$ \text{\scriptsize most efficient Agents} & $p_{2,\gamma}(t)c + p_{3,\gamma}(t)$ &  $ p_ {2,\gamma}^1(t)c + p_{3,\gamma}^1(t)$, for $c<\hat{c}_1^\gamma(t)$ \\
\hline
		$\left[b_0,a_0\right]$ \text{\scriptsize intermediate Agents} & not picked &  not picked\\
		\hline
		$\left[0,b_0\right]$ \text{\scriptsize least efficient Agents} & not picked & $p_{1,\gamma}^3(t)c^\gamma + p_{2,\gamma}^3(t)c + p_{3,\gamma}^3(t)$, for $\hat{c}_2^\gamma(t) <c$  \\
		\hline
	\end{tabular}
	
	\vspace{1em}
	\begin{tabular}{|c|c|c|}
		\hline
		\text{Selected Agents}&  \text{$H$ constant} &     \text{$H$ concave non decreasing}  \\
		\hline 
		$\left[a_0,1\right]$ \text{\scriptsize most efficient Agents} & $p_{1,\gamma}(t)c^\gamma+ p_{2,\gamma}(t)c+p_{3,\gamma}(t)$ &  $ p_ {1,\gamma}^1(t)c^\gamma + p_{2,\gamma}^1(t)c+p^1_{3,\gamma}(t)$, for $c>\hat{c}_2^\gamma(t)$ \\
\hline
$\left[b_0,a_0\right]$ \text{\scriptsize intermediate Agents} & not picked &  not picked\\
		\hline
		$\left[0,b_0\right]$ \text{\scriptsize least efficient Agents} & not picked & $p^3_{2,\gamma}(t)c + p^3_{3,\gamma}(t)$, for $c<\hat{c}_1^\gamma(t) $  \\
		\hline
	\end{tabular}

	\end{center}
\caption{Optimal tariff of residential consumers $\gamma<0$ (top), and industrial consumers $\gamma\in\left(0,1\right)$ (bottom).}
	\end{table}


\medskip
Let us interpret the optimal tariffs and connect the three components to electricity pricing standard issues. As already said, $p_3$ represents the fixed charge. The volumetric charge is the combination of a standard term $p_2(t)c$ plus $p_1(t)c^\gamma$, where the latter is a way to charge more high demand consumers (indeed it only appears when $c$ is high enough). Finally, no explicit demand charge appears but the coefficients $(p_i)_{1\leq i\leq 3}$ depend on the maximum subscribed power $(\hat{c}_i^\gamma)_{1\leq i\leq 3}$, which limits the instantaneous power use and allows to charge more high power consumers. Let us point out that this method of electricity pricing is implementable in practice, thanks to the recent spread and development of smart meters, which enables a precise metering of electricity consumption, and a dynamic management of maximum power\footnote{See for instance the USmartConsumer report \cite{usmartconsumer2016european}, which states that at the end of 2016, 30\% of overall European electricity meters were equipped with smart technology }. In addition, the peak/off--peak issues are handled by the temporal structure of the tariff and high power consumption within peak period will be overcharged compared to off--peak period. Let us also mention that the proportional part $p_1(t)$ to $c^\gamma$ only depends on the Agent's utility parameters. Therefore this part should be common to any Principal, whatever her cost of production is, or the reservation utility of the consumers. 

\medskip
The selected Agents are the most efficient (highest utility for a given consumption) when $H$ is constant, which is a classical result. But when $H$ is concave, the Principal can also select Agents among the least efficient ones. Indeed, in this case reaching most efficient Agents is costly (
they require a high amount of electricity
) and it happens that getting less efficient Agents can be profitable as they are more easily satisfied (
they require a lower volume of electricity at a higher price per unit, compared to efficient Agents). This type of feature seems to be uncommon in the Principal--Agent literature and can only be found, as far as we know, in models with countervailing incentives (see \cite{lewis1989countervailing}).

\section{Economic interpretations and numerical results}	\label{sec:interpret}

\begin{figure}[ht!]
	\begin{center}
		\hspace*{-6mm}
			\includegraphics[scale=0.35]{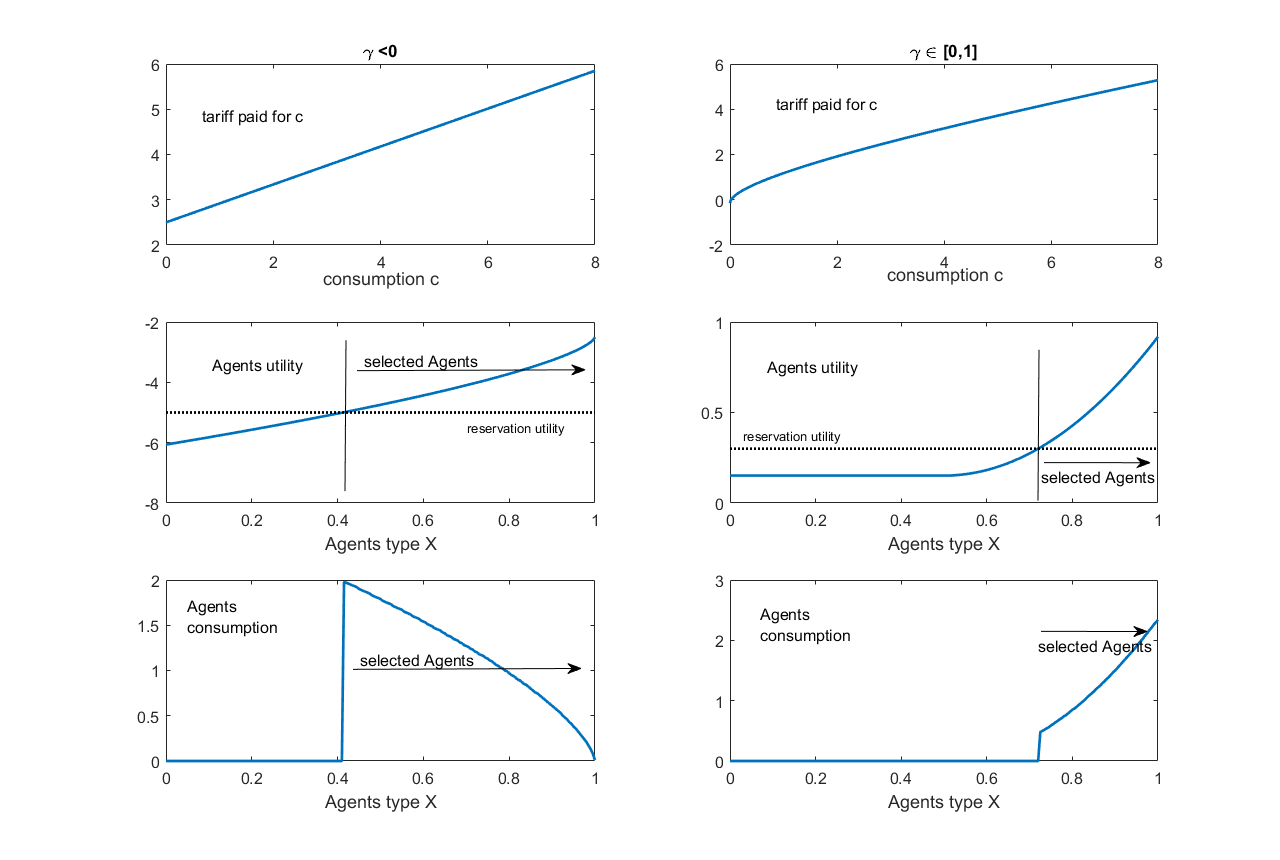} 	
	\small
		\caption{\label{Hconstant} Tariffs paid against consumption (upper graphs), Agents' utility against type (middle graphs), selected consumption against Agents' type (lower graphs); $H$ constant.}
	\end{center}
\end{figure}

\textbf{Examples when $H$ is constant}. For a constant reservation utility, the most efficient Agents are selected. We present numerical illustrations in Figure \ref{Hconstant}. The tariff structure is linear in consumption when $\gamma<0$ and is concave when $\gamma\in(0,1)$, which is represented in the upper graphics of Figure \ref{Hconstant}. Middle graphics represent the utility Agents can obtain by signing the contract, depending on their type. If this utility level is smaller than their reservation utility (represented by the dashed line) they do not enter the contract and their consumption is null, as represented in the lower graphics. These utility representations also illustrate a classical result of informational rent: the most efficient Agents obtain a tariff lower to what they are willing to pay, whereas the least efficient ones need to pay as much as they are able to, or are excluded.

\begin{center}
\begin{figure}[!ht]
	\begin{center}
		\hspace*{-6mm}
		\includegraphics[scale=0.35]{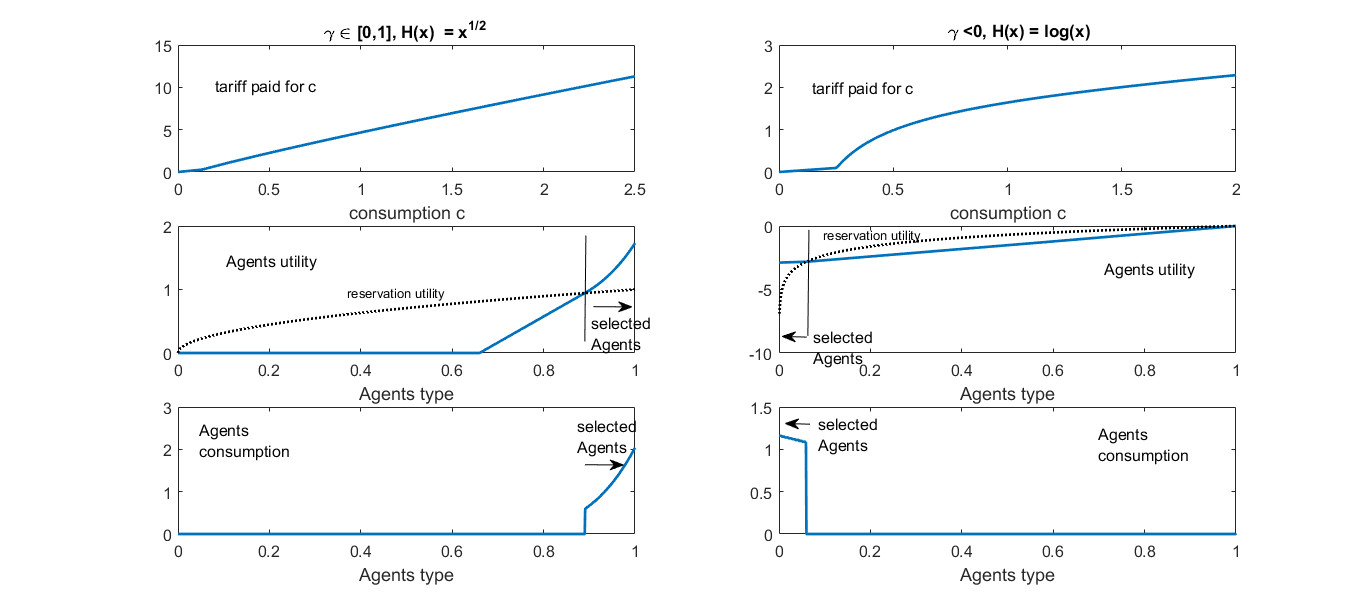} 
		\caption{\label{Hconcave} Tariffs against consumption (upper graphs), Agents's utility against type (middle graphs), and selected consumption against types (lower graphs); $H$ concave.}
	\end{center}
	
\end{figure}
	\end{center}

\textbf{Examples when $H$ is concave.} For a concave reservation utility, not only most efficient Agents are selected. We provide numerical illustrations where either the most or the least efficient Agents are selected on Figure \ref{Hconcave}. First, let us analyse the example when $H(x)=\sqrt{x}$ and $\gamma \in (0,1)$, which corresponds to the left column. In this example, only the most efficient Agents sign the contract as they are the only ones obtaining a higher utility than their reservation one. As presented in the previous section, the tariff structure is the combination of three functions of consumption (upper graphics), but Agents who sign the contract only choose consumption such that $\hat{c}_2^\gamma<c$ which corresponds to the concave tariff part $p_1^3(t)c^\gamma + p_2^3(t)c + p_3^3(t)$.

\medskip
When $H(x)=\log{(x)}$ and $\gamma <0$ (right column of Figure \ref{Hconcave}), only the least efficient Agents take the contract. Indeed, the concavity of the reservation utility makes it profitable for the Principal to select these Agents, rather than the most efficient ones. The tariff structure is again the combination of three functions of consumption (upper graphics) but Agents who sign the contract in this example only take consumption such that $\hat{c}_2^\gamma<c$ (of course $\hat{c}_2$ is different from the one in the previous example because we consider a different $H$). This again corresponds to the concave tariff part $p_1^3(t)c^\gamma + p_2^3(t)c + p_3^3(t)$.

\begin{figure}[H]
	\centering
	\includegraphics[scale=0.3]{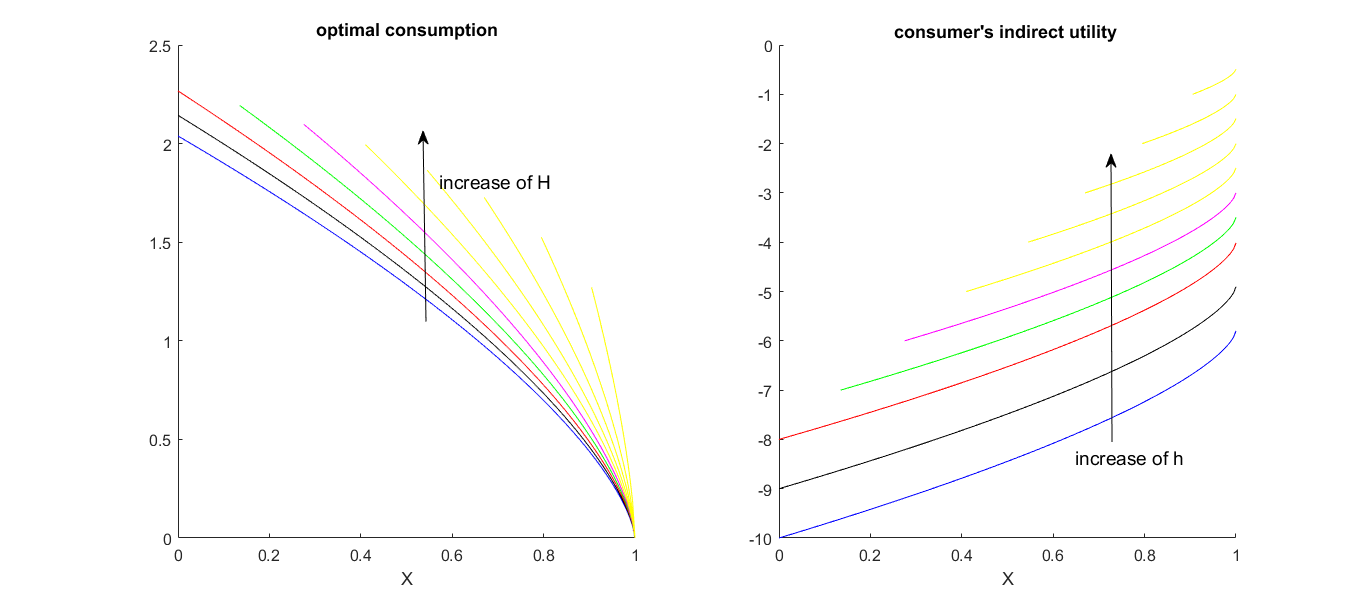}
	\caption{\label{Hincrease} Evolution of Agents' utility (left) and consumption (right) against type; $H$ increases and $\gamma<0$.}
\end{figure}

\begin{figure}[H]
	\centering
	\includegraphics[scale=0.3]{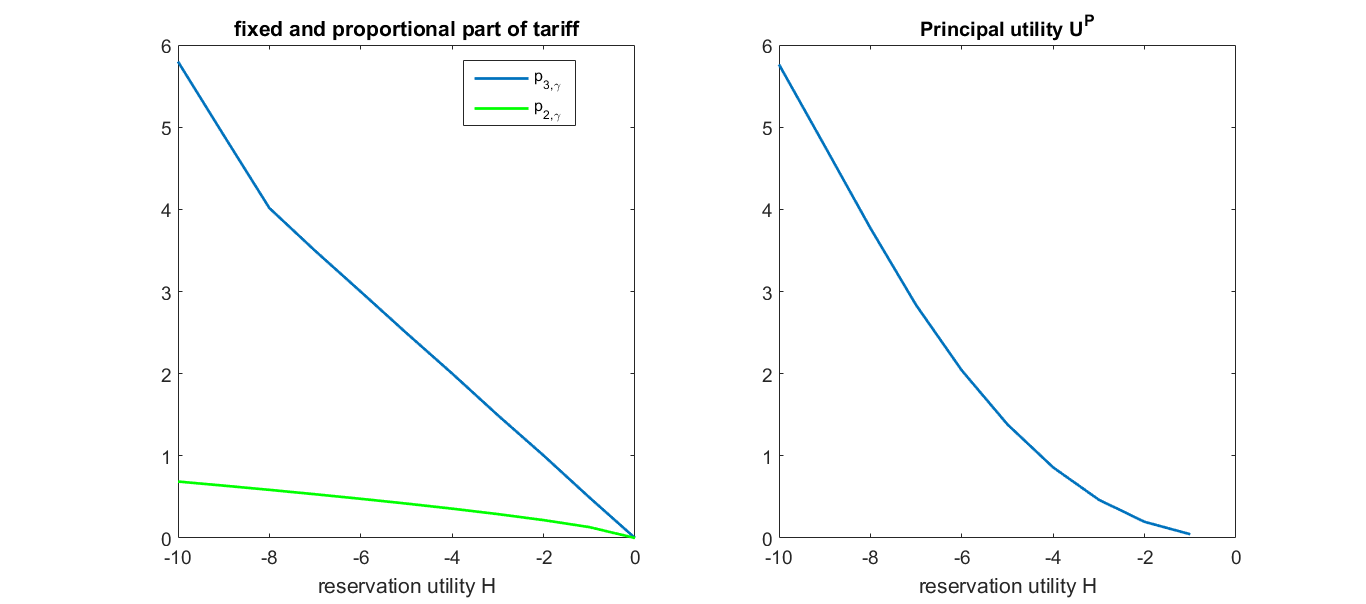}
	\caption{\label{Hincrease_tariff} Evolution of tariff's components (left) and Principal utility $U^P$ (right); $H$ increases and $\gamma<0$. }
\end{figure}

\textbf{Impact of competition when $H$ is constant.} For a bigger function $H$, because competition is more intense, the Principal adapts her tariff in order to remain competitive. In that case, the Principal mainly decreases the constant part $p_{3,\gamma}$ of the tariff in order to attract consumers (see the left graphic of Figure \ref{Hincrease_tariff} when $\gamma<0$). The consumers selecting this new tariff obtain better conditions and as such consume more power, because it is cheaper, see the same example on the left graphic of Figure \ref{Hincrease} when $\gamma<0$. Therefore, when the Principal decreases his tariff, he does not decrease it enough in order to keep the same quantity of consumers: he accepts to retain less of them, but the selected ones do consume more, as represented on the right graphic of Figure \ref{Hincrease}. Nevertheless, the utility of the Principal decreases with competition, see the right part of Figure \ref{Hincrease_tariff}. At the extreme, the Principal even offers no tariff whenever $H$ is too high. Observe that in this example, the fixed part of the tariff represents more than a half of the total cost of electricity for the consumers.

\begin{figure}[H]
	\centering
	\includegraphics[scale=0.3]{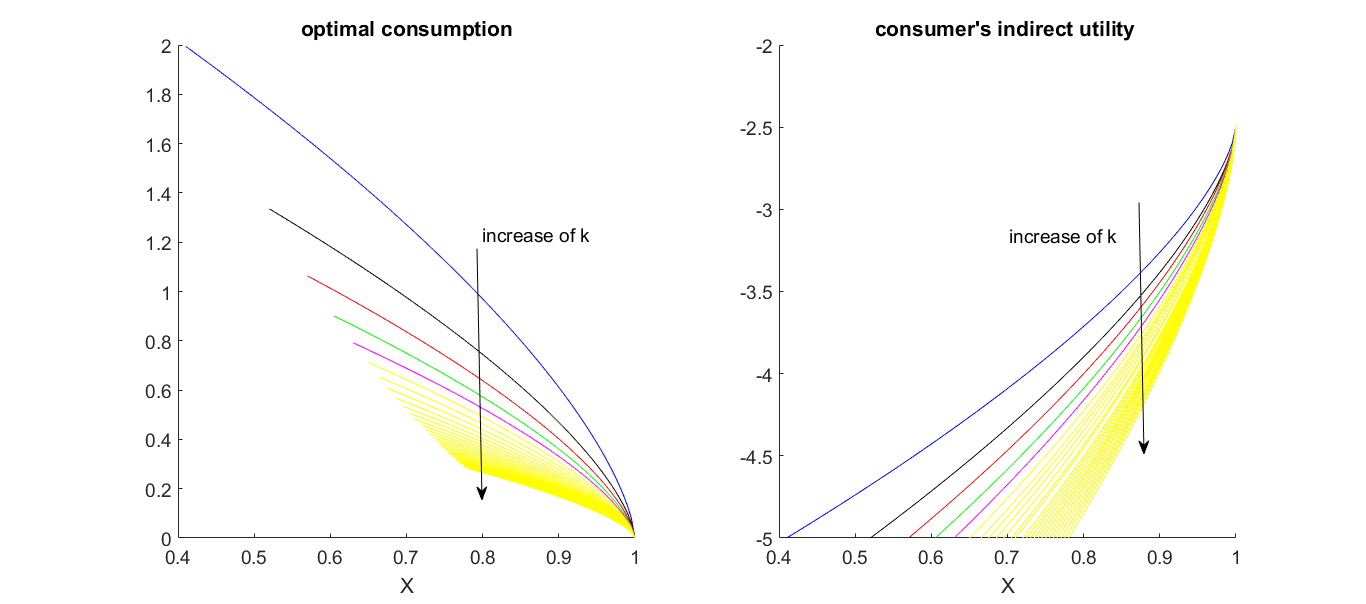}
	\caption{\label{Kincrease} Agents' utility (left) and consumption (right) against type; $k$ increases and $\gamma<0$.}
\end{figure}

\begin{figure}[H]
	\centering
	\includegraphics[scale=0.3]{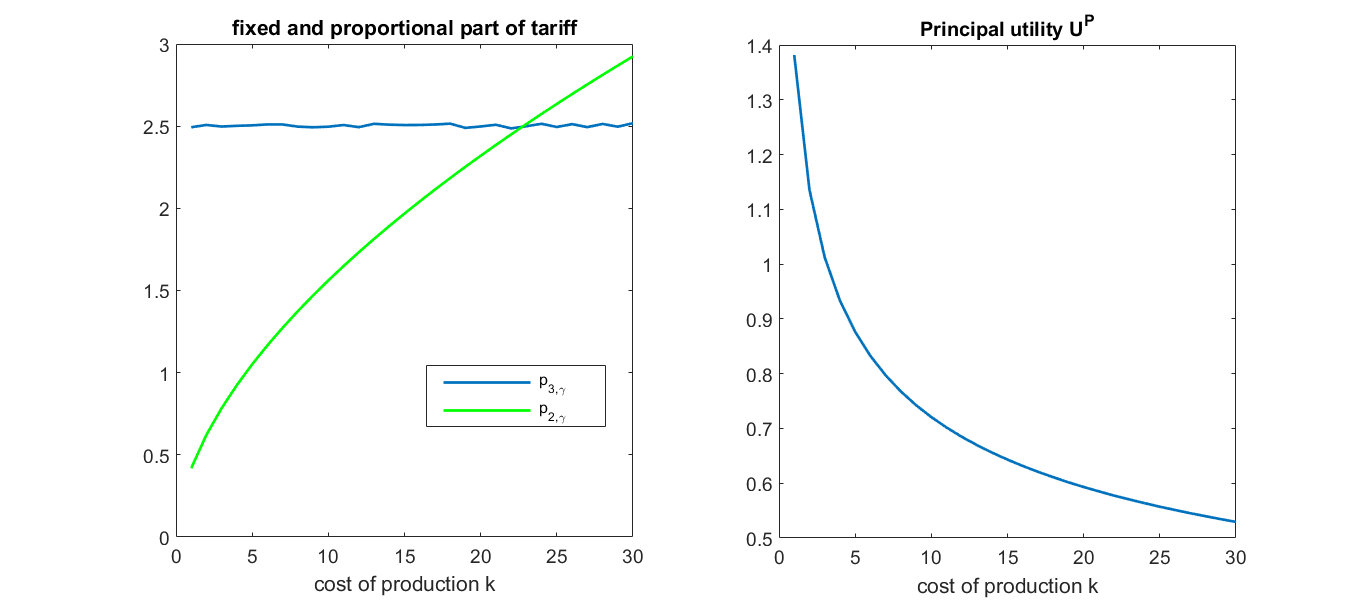}
	\caption{\label{Kincrease_tariff} Tariff's components (left) and Principal utility $U^P$ (right) against  $k$; $\gamma<0$.}
\end{figure}

\textbf{Impact of cost of production when $H$ is constant.} For an increase of the cost of production $k$, the Principal also adapts his tariff in order to reflect this cost modification. In that case, the Principal mainly increases the proportional part $p_{2,\gamma}$ of its tariff, see the example for constant $H$ and $\gamma<0$ on the left part of Figure \ref{Kincrease_tariff}. Consumers who select this new tariff are offered worse conditions, and as such consume less power, because it is more expensive, see the same example on the left graphic of Figure \ref{Kincrease}. In addition, less consumers select the contract. Therefore, the utility of the Principal decreases with the cost of production, as illustrated in the right part of Figure of \ref{Kincrease_tariff}. On the contrary, some production technologies like renewable have no variable production costs and only fixed costs which corresponds mainly to the investment and the maintenance costs. If we imagine a system with a very large share of these types of technology, the variable cost $k$ could be very low: energy is not expensive only the installation of equipment is costly. In that case, the proportional part $p_2$ goes to zero. This means that for residential consumers ($\gamma <0$), the electricity tariff would reduce only to a fixed charge (see illustration on figure \ref{Kincrease_tariff}) and consumers would pay the same whatever their consumption. But their consumption is "naturally" limited by the saturation for electricity expressed in their utility function. For industrial Agents ($\gamma\in \left[0,1\right]$), the electricity tariff would reduce to a fixed charge $p_3$ and to a proportional part $p_1$ of $c^\gamma$. Therefore, this is this last term which enables to limit the consumption for industrial Agents. 

Let us also point out that we can simulate the impact of the convexity of the cost of production $K$ and this shows that the structure of the optimal tariff is unchanged when $n\rightarrow 1$.

\section{Model specification} \label{sec:model}

We now turn to a more precise exposition of the model and try to present it in a rather general setting. We will state the main ideas and methodologies in the highest level of generality. Specific assumptions on the shape on the utility, type distribution or cost functions will only be introduced later, in order to present our results in a more explicit fashion. 

\medskip
In this model, the Principal is a power company, whose purpose is to offer to its clients a collection of tariffs which maximise its profits. The time horizon $T>0$ is fixed. The following notations are used throughout the article.
\begin{itemize}
	\item $\mathcal C$ represents the admissible levels of consumption for the Agents, and it is either equal to $\R_+$ or $\R_+^*$ depending on the utility function of the Agents.
	\item $p:[0,T]\times \mathcal C\longrightarrow \R_+$ is the tariff proposed by the Principal, such that $p(t,c)$ represents the instantaneous price of electricity at time $t$ corresponding to a level of consumption $c$.
	\item $K:[0,T]\times\mathcal C\longrightarrow\R_+$ is the cost of production of electricity for the Principal, such that $K(t,c)$ represents the cost at time $t$ for an aggregate level of production $c$. We assume that $K$ is continuous in $t$, increasing, continuously differentiable and strictly convex in $c$.
	\item $x$ is the Agent's type, assumed to take values in some subset $X$ of $\mathbb R$.
	\item $c:[0,T]\times X\longrightarrow \mathcal C$ is the consumption function, such that $c(t,x)$ represents the consumption of electricity by an Agent of type $x$ at time $t$.
	\item $u:[0,T]\times X\times\mathcal C\longrightarrow \R$ is the utility function of the Agents, such that $u(t,x,c)$ represents the utility obtained by an Agent of type $x$ at time $t$ when he consumes $c$. We assume that the map $c\longmapsto u(t,x,c)$ is non-decreasing and concave for every $(t,x)\in[0,T]\times X$. Moreover, the map $u$ is assumed to be jointly continuous, such that $x\longmapsto u(t,x,c)$ is non-decreasing and differentiable Lebesgue almost everywhere for every $(t,c)\in[0,T]\times\Cc$, and such that $c\longmapsto  \frac{\partial u}{\partial x}(t,x,c)$ is invertible. Finally, we assume that if $\Cc=\R_+$, the value $u(t,x,0)\in\R_+$ is independent of $x$, and if $\Cc=\R_+^*$ that $\lim_{c\to 0}u(t,x,c)=-\infty$, for every $(t,x)\in[0,T]\times X$. In other words, all the Agents have the same utility when they do not consume electricity.
	\item $f:X\longrightarrow \R_+$ is the distribution of the Agent's type over the population. As is customary in adverse selection problems, $f$ is supposed to be known by the Principal.
\end{itemize} 
%
%
\subsection{Agent's problem} \label{sec:Agent}

Let us start by defining the consumption strategies that the Agents are allowed to use. A consumption strategy $c$ will be said to be admissible, which we denote by $c\in\mathfrak C$, if it is a Borel measurable map from $[0,T]$ to $\Cc$. Given a tariff $p$, that is a map from $[0,T]\times\R_+$ to $\R$, proposed by the Principal, an Agent of type $x\in X$ determines his consumption by solving the following problem
\begin{equation}\label{prob:Agent}
U_A(p,x):=\sup_{c\in\mathfrak C} \int_{0}^{T} \big( u(t,x,c(t)) - p(t,c(t))\big)\mathrm{d}t.
\end{equation}

The tariff that the Principal can offer to the Agents has to satisfy the incentive compatibility (IC) and the individual rationality (IR) conditions. In our setting, there is no moral hazard, so that the incentive compatibility condition is automatically satisfied. Furthermore, the (IR) condition can be expressed through the set $X(p)$ of the types of Agents which accept the contract $p$, which can be defined as
$$
X(p):  = \left\{ x\in X,\; U_A(p,x) \geq H(x) \right\},
$$ 
with a continuous and non-decreasing function H which represents the reservation utility of the Agents of different types, that is to say the utility that the Agents can hope to obtain by subscribing their power contract with a competitor. Agents for which the map $U_A(p,\cdot)$ is smaller than $H$ will not accept the contract offered by the Principal.

\medskip
We are now ready to give our definition of admissible tariffs, which uses vocabulary from $u-$convex analysis. We have regrouped all the pertinent results and definitions in Appendix \ref{Appendix A}, for readers not familiar with this theory.
\begin{definition}\label{def:admissible tariff}
	A tariff $p:[0,T]\times\Cc\longrightarrow \R$ will be said to be admissible, denoted by $p\in\Pc$, if it satisfies
	\begin{itemize}
		\item[$(i)$] For any $(t,x)\in[0,T]\times X(p)$ the set $\partial^\star p^\star(t,x)$ is non--empty.
		\item[$(ii)$] The map $x\longmapsto p^\star(t,x)$ is continuous on $X$, differentiable Lebesgue almost everywhere, for every $t\in[0,T]$, and satisfies
		$$\int_0^T\int_X\left|\frac{\partial p^\star}{\partial x}\right|(t,x)\mathrm{d}x\mathrm{d}t<+\infty.$$
		\item[$(iii)$] If one defines the map $c^\star:[0,T]\times[0,1]\longrightarrow \mathbb R_+$ by
		\begin{equation}
		c^\star (t,x) = \left(\frac{\partial u}{\partial x}(t,x,\cdot)\right)^{(-1)} \left( \frac{\partial p^\star }{\partial x}(t,x) \right),
		\label{optim_conso}
		\end{equation}
		then the restriction of $p$ to $\{(t,c)\in[0,T]\times\Cc,\; \exists x\in X(p),\; c=c^\star(t,x)\}$ is $u-$convex.
	\end{itemize}
\end{definition}   
Let us comment on the above definition. First of all, the regularity assumptions are mainly technical. In principle, we would like to deal with tariffs $p$ which are $u-$convex since they are completely characterised by their $u-$transform $p^\star$, which will be our object of study later on. However, we can ask for less than that. The main point here is that since only the clients with type in $X(p)$ are going to accept the contract, the Principal will only have to face consumption levels chosen by these clients. Besides, as we are going to prove in the next proposition, this optimal consumption is exactly $c^\star(t,x)$. Therefore, any consumption $c\in\Cc$ which does not belong to the pre--image of $X(p)$ will never have to be considered by the Principal. In particular, there is a degree of freedom when defining the value of the tariff $p$ there. Indeed, if clients of some type $x$ reject the contract $p$, they will reject any contract with a higher price. This is the reason why we do not impose the admissible tariffs to be $u-$convex on $\Cc$ but only on the corresponding pre--image of the set $X(p)$.
%
%
%

\medskip
Our main result in this section is
\begin{proposition} \label{prop:agent-problem}
	For every $p\in\Pc$ and for almost every $x\in X(p)$, we have
	\[U_A(p,x)=\int_0^Tp^\star(t,x)\mathrm{d}t,\]
	and the optimal consumption of Agents of type $x$ at any time $t\in[0,T]$ is given by $c^\star (t,x)$ defined in \eqref{optim_conso}. In particular, $X(p)$ can be defined through $p^\star$ only as follows
	\[X(p)=X^\star(p^\star):=\left\{ x\in X,\; P^\star(x):=\int_0^Tp^\star(t,x)\mathrm{d}t \geq H(x) \right\}.\]
\end{proposition}

\subsection{The Principal's problem}
%
%
%
The Principal sets a tariff $p\in\Pc$ as a solution to her maximisation problem
\begin{equation}
U_P:=	\sup_{p\in\Pc} \int_{0}^{T} \bigg[ \int_{X(p) } p(t,c^\star(t,x)) f(x) \mathrm{d}x  - K\bigg(t,\int_{X (p)} c^\star(t,x) f(x)\mathrm{d}x\bigg)  \bigg] \mathrm{d}t.
\label{Principla_obj}
\end{equation}

\medskip
Using the results of Section \ref{sec:Agent}, we can rewrite this problem in terms of $p^\star$ only as
\begin{align}
\nonumber U_P=&\	\sup_{p\in\Pc} \int_{0}^{T} \bigg[ \int_{X^\star(p^\star) } \bigg(u\bigg(t,x,\left(\frac{\partial u}{\partial x}(t,x,\cdot)\right)^{(-1)} \bigg( \frac{\partial p^\star }{\partial x}(t,x) \bigg) 
\bigg)-p^\star(t,x)\bigg) f(x) \mathrm{d}x\\
& \hspace{4em}- K\bigg(t,\int_{X^\star (p^\star)} \left(\frac{\partial u}{\partial x}(t,x,\cdot)\right)^{(-1)} \left( \frac{\partial p^\star }{\partial x}(t,x) \right) 
f(x)\mathrm{d}x\bigg)  \bigg] \mathrm{d}t,
\label{Principla_obj2}
\end{align}

\medskip
Now notice from \eqref{optim_conso} that $p^\star$ is actually non-decreasing in $x$ (since $x\longmapsto u(t,x,c)$ is non-decreasing for every $(t,c)\in[0,T]\times \Cc$). Let us then consider the space $C^{+}$ of maps $g$, such that for every $t\in[0,T]$, $x\longmapsto g(t,x)$ is continuous and non--decreasing with
$$\int_0^T\int_X\abs{\frac{\partial g}{\partial x}(t,y)}\mathrm{d}y\mathrm{d}t<+\infty.$$

We shall actually consider a relaxation of the problem of the Principal $\widetilde U_P\geq U_P$, defined by
\begin{align}
\nonumber \widetilde U_P:=&	\sup_{p^\star\in C^+} \int_{0}^{T} \bigg[ \int_{X^\star(p^\star) } \bigg(u\bigg(t,x,\left(\frac{\partial u}{\partial x}(t,x,\cdot)\right)^{(-1)} \bigg( \frac{\partial p^\star }{\partial x}(t,x) \bigg) 
\bigg)-p^\star(t,x)\bigg) f(x) \mathrm{d}x\\
& \hspace{4em}- K\bigg(t,\int_{X^\star (p^\star)} \left(\frac{\partial u}{\partial x}(t,x,\cdot)\right)^{(-1)} \left( \frac{\partial p^\star }{\partial x}(t,x) \right) 
f(x)\mathrm{d}x\bigg)  \bigg] \mathrm{d}t,
\label{Principla_obj3}
\end{align}
where we have forgotten the implicit link existing between $p$ and $p^\star$, which explains why we have in general $\widetilde U_P\geq U_P$. We will see in the frameworks described below that we can give conditions under which the two problems are indeed equal. The main advantage of $\widetilde U_P$ is that it no longer contains the condition that $p^\star$ has to be $u-$convex, a constraint that is not easy to consider in full generality.

\medskip
We emphasise that problem $\widetilde U_P$ is well defined, since the elements of $C^+$ are non--decreasing with respect to $x$ and thus differentiable Lebesgue almost everywhere. Our aim now will be to compute $\widetilde U_P$. However, the present framework is far too general to hope obtaining explicit solutions, which are of the utmost interest in our electricity pricing model, so for the rest of the paper we will concentrate our attention on the case of Agents with power--type CRRA utilities. 

\subsection{Agents with CRRA utilities}

For the sake of tractability, we shall use the following standing assumptions

\begin{assumption}\label{assump:power}
	$(i)$ $X=[0,1]$.
	
	\medskip
	$(ii)$ We have for every $(t,x,c)\in[0,T]\times X\times \Cc$
	\begin{equation*}
	u(t,x,c) = g_\gamma(x)\phi(t) \frac{c^\gamma}{\gamma},
	\end{equation*}
	for some $\gamma\in(-\infty,0)\cup(0,1)$, some map $g_\gamma:X\longrightarrow \R_+$ which is continuous, increasing if $\gamma\in(0,1)$, decreasing if $\gamma\in(-\infty,0)$, and for some continuous map $\phi:[0,T]\longrightarrow \R_+^\star$. 
\end{assumption}
Let us comment on this modelling choice for the utility function. The term $g_\gamma(x)$ represents the willingness of the Agents to pay for their consumption, {\sl i.e.} their need for energy depends on their type. The term $ \phi$ is common to every type of Agents and represents the fact that (almost) everyone is eager to consume at the same time (for example during the day rather than at night). Furthermore, we consider both the cases $\gamma\in(0,1)$, which would be the classical power utility function, as well as the case $\gamma<0$, which corresponds to a situation where Agents actually cannot avoid consuming electricity, as it would provide them a utility equal to $-\infty$, which may be seen as more realistic. As discussed previously, taking $\gamma\in(0,1)$ identifies to considering industrial Agents, and $\gamma<0$ more typically refers to residential Agents.

\medskip
Equation \eqref{optim_conso} now can be written as
\begin{equation} \label{optim_conso_benchmark2}
c^\star(t,x) = \left(\frac{\gamma}{\phi \left(t\right)g'_\gamma(x)} \frac{\partial p^\star }{\partial x} \left(t,x\right)\right)^{\frac{1}{\gamma}}.
\end{equation}

By inserting the previous expression in equation \eqref{Principla_obj} and using that 
$$
p(t,c^\star (t,x))=  g_\gamma(x) \phi(t) \frac{c^\star (t,x)^\gamma}{\gamma} - p^\star (t,x),
$$
the relaxed problem to solve can now be expressed as
\begin{align}\label{Principla_obj_benchmark2}
 \widetilde U_P=&\sup_{p^\star\in C^+ } \int_{0}^{T}  \bigg[ \int_{X^\star(p^\star) }  \left( \frac{g_\gamma(x)}{g'_\gamma(x)}\frac{\partial p^\star }{\partial x} (t,x)  - p^\star (t,x)\right)   f(x) \mathrm{d}x  - K\bigg(t,  \int_{X^\star(p^\star) }  \left( \frac{\gamma}{\phi (t)g'_\gamma(x)} \frac{\partial p^\star }{\partial x} (t,x) \right)^{\frac{1}{\gamma}} f(x)\mathrm{d}x\bigg)  \bigg] \mathrm{d}t.
\end{align}

\medskip
The aim of next sections is to solve problem $\widetilde U_P$. In Section \ref{section-constant}, we will focus on the particular case where the reservation utility $H$ is constant, and shall consider the more general case where it may depend on the Agents' type in Section \ref{sec:general}.

\section{Constant reservation utility}\label{section-constant}

We consider in this section a further simplification, related to the reservation utility of the Agents, which we suppose to be independent of their type. This assumption will be relaxed in Section \ref{sec:general}.
\begin{assumption}\label{assump:hconst}
	The reservation utility $H$ is actually independent of $x$, that is
	\[H(x)=:H,\ \text{for every $x\in[0,1]$}.\]
\end{assumption}
When the function $H$ is constant, the set of Agents who accept any tariff becomes an interval. Indeed, under Assumptions \ref{assump:power} and \ref{assump:hconst}, the (IR) condition reduces to
\[
X^\star(p^\star)  = \left\{ x\in[0,1],\  \int_{0}^{T} p^\star(t,x)\mathrm{d}t \geq H \right\}.
\]
Since $p^\star $ is non--decreasing in $x$, we have for any $x_0\in[0,1]$ that
\[
\int_0^T p^\star (t,x_0)\mathrm{d}x \geq H ~\Longrightarrow~  \int_0^T p^\star (t,x)\mathrm{d}x \geq H, \ \forall x\geq x_0.
\]
Therefore, the set $X^\star(p^\star) $ has necessarily the form
\[
X^\star(p^\star) =[x_0,1],
\]
where $x_0\in[0,1]$ needs to be determined and verifies, by continuity, that $P^\star(x_0)=H$. This means that the Principal will only select the Agents of high type, greater than some value $x_0$. 

\medskip
We can now look at an equivalent formulation of the problem of the Principal, in which we first choose $x_0$, and then optimise over the tariffs for which $X^\star(p^\star)$ is exactly $[x_0,1]$. The problem \eqref{Principla_obj_benchmark2} can therefore be written as
\begin{align}\label{Principla_obj_benchmark2_simple}
\nonumber&\widetilde U_P=\sup_{x_0 \in[0,1]}\ \sup_{p^\star\in C^+(x_0)}\ \int_{0}^{T} \bigg[ \int_{x_0}^1  \left( \frac{g_\gamma(x)}{g'_\gamma(x)}\frac{\partial p^\star }{\partial x} (t,x)    - p^\star (t,x)\right)f(x)  \mathrm{d}x  \\ 
&\hspace{10em}-  K\bigg( t,\int_{x_0}^1  \left( \frac{\gamma}{\phi(t)g'_\gamma(x)} \frac{\partial p^\star }{\partial x} (t,x) \right)^{\frac{1}{\gamma}}f(x)\mathrm{d}x \bigg)  \bigg] \mathrm{d}t,
\end{align}
with 
$$
C^+(x_0):=\left\{p^\star\in C^+,\ \int_0^T p^\star(t,x_0)\mathrm{d}t = H\right\} =\left\{p^\star\in C^+,\ X^\star(p^\star)=[x_0,1]\right\}.
$$
Formulation \eqref{Principla_obj_benchmark2_simple} is convenient because it allows to reduce problem $\widetilde U_P$ to a one--dimensional one. To do so, we first need to find the best map $p^\star$ in the set $C^+(x_0)$, for each $x_0$. 

\medskip
Finally, we will rewrite the relaxed problem of the Principal once more, in order to focus only in the derivative of the $u-$transform $p^\star$. Denote by $F$ the cumulative distribution function of the types of Agents. By integration by parts we have for every $x_0\in[0,1]$ and every $p^\star\in C^+(x_0)$
{\small
	\begin{align*}
	& \int_{0}^{T} \bigg[ \int_{x_0}^1  \left( \frac{g_\gamma(x)}{g'_\gamma(x)} \frac{\partial p^\star }{\partial x} (t,x)    - p^\star (t,x)\right)f(x)  \mathrm{d}x  - K\bigg( t,\int_{x_0}^1  \left( \frac{\gamma}{\phi(t)g'_\gamma(x)} \frac{\partial p^\star }{\partial x} (t,x) \right)^{\frac{1}{\gamma}}f(x)\mathrm{d}x \bigg)  \bigg] \mathrm{d}t\\
	&=\int_{0}^{T} \bigg[ \int_{x_0}^1  \left( \frac{g_\gamma(x)}{g'_\gamma(x)}f(x)+F(x)-1\right) \frac{\partial p^\star }{\partial x} (t,x)   \mathrm{d}x  - K\bigg( t,\int_{x_0}^1  \left( \frac{\gamma}{\phi(t)g'_\gamma(x)} \frac{\partial p^\star }{\partial x} (t,x) \right)^{\frac{1}{\gamma}} f(x)\mathrm{d}x \bigg)\bigg]\mathrm{d}t\\
	&\hspace{0.9em}  +(F(x_0)-1)\int_0^Tp^\star (t,x_0)\mathrm{d}t.
	\end{align*}
}
We therefore end up with the maximization problem
{\small \begin{align} \label{solvableproblem}
	&\nonumber\widetilde U_P=\sup_{x_0 \in[0,1]}\ \sup_{p^\star\in C^+(x_0)}\ \int_{0}^{T} \bigg[ \int_{x_0}^1  \frac{\left( g_\gamma(x)f(x)+g'_\gamma(x)F(x)-g'_\gamma(x) \right)}{g'_\gamma(x)} \frac{\partial p^\star}{\partial x}(t,x)  \mathrm{d}x  \\ 
	&\hspace{10em}- K\bigg( t,\int_{x_0}^1  \left( \frac{\gamma}{\phi(t)g'_\gamma(x)} \frac{\partial p^\star }{\partial x} (t,x) \right)^{\frac{1}{\gamma}}f(x)\mathrm{d}x \bigg) \bigg] \mathrm{d}t +(F(x_0)-1)H.
	\end{align} }

\medskip
We can now state our main result of this section, providing the solution to the relaxed problem $\widetilde U_P$ and conditions under which we can recover the solution to the original problem $U_P$. For ease of presentation, we introduce the following function
\[
\ell(x_0):=\int_{x_0}^1  \left( \frac{ \left[ g_\gamma(x)f(x)+g'_\gamma(x)F(x)-g'_\gamma(x) \right]^+ }{f^\gamma(x)}\right)^{\frac{1}{1-\gamma}}   \mathrm{d}x,\ x_0\in[0,1].
\]

\begin{theorem}\label{th:main}
	Let Assumptions \ref{assump:power} and \ref{assump:hconst} hold. We have
	
	\medskip
	$(i)$ The maximum in \eqref{Principla_obj_benchmark2_simple} is attained for the maps 
	\begin{align*}
	p^\star(t,x)=p^\star(t,x_0^\star)+\int_{x^\star_0}^x\frac{g'_\gamma(y)}{\gamma}\Bigg(\frac{\phi(t)^{\frac{1}{\gamma}}\big[ g_\gamma(y)f(y)+g'_\gamma(y)F(y)-g'_\gamma(y)\big]^+ }{f(y)\frac{\partial K}{\partial c}(t,A(t,x_0^\star))}\Bigg)^{\frac{\gamma}{1-\gamma}}\mathrm{d}y,\ x\in[0,1],
	\end{align*}
	where $A(t,x_0)$ is defined
	\[
	A(t,x_0) := \int_{x_0}^1  \left( \frac{\gamma}{\phi(t)g'_\gamma(x)} \frac{\partial p^\star }{\partial x} (t,x) \right)^{\frac{1}{\gamma}}  f(x) \mathrm{d}x.
	\]
	and $x_0^\star$ is any maximiser of the map 
	\begin{align*}
	[0,1]\ni x_0 \longmapsto &\int_{0}^{T}\left( \frac{\partial K}{\partial c}(t,A(t,x_0))m(t,x_0)- K( t,\gamma m(t,x_0))\right)\mathrm{d}t +(F(x_0)-1)H,
	\end{align*}
	with
	\[m(t,x_0):=\frac{\phi^{\frac{1}{1-\gamma}}(t)\ell(x_0)}{\gamma\displaystyle\left(\frac{\partial K}{\partial c}(t,A(t,x_0))\right)^{\frac{1}{1-\gamma}}}, \; (t,x_0)\in[0,T]\times[0,1],\]
	and $t\longmapsto p^\star(t,x_0^\star)$ is any map such that
	\[\int_0^Tp^\star(t,x_0^\star)\mathrm{d}t=H.\]
	For instance, one can choose $p^\star(t,x_0^\star):=H/T,$ $t\in[0,T]$. 
	
	\medskip
	$(ii)$ Define $p$ for any $(t,c)\in[0,T]\times\R_+$ by
	\[p(t,c)=\underset{x\in [0,1]}{\sup}\left\{g_\gamma(x)\phi(t)\frac{c^\gamma}{\gamma}-p^\star(t,x)\right\}.\]
	If the map defined on $[0,1]$ by
	\[ x\longmapsto g_\gamma'(x)\left(\frac{\left[ {g_\gamma(x)}f(x)+{g_\gamma'(x)}F(x)-g_\gamma'(x)\right]^+ }{f(x)}\right)^{\frac{\gamma}{1-\gamma}},\]
	is non--decreasing, then $p^\star$ is $u-$convex, and $p$ is the optimal tariff for the problem \eqref{Principla_obj}. Furthermore, the Principal only signs contracts with the Agents of type $x\in[x_0^\star,1]$. 
	
	\medskip
	$(iii)$ Finally, in the case $\gamma\in(0,1)$, if $f$ is non--increasing and the map
	\[
	\beta:  x\longmapsto  \frac{ (g_\gamma(x)f(x)+g'_\gamma(x)F(x)-g'_\gamma(x) )}{f^\gamma(x)},
	\]
	is increasing over the set $L:=\left\{ x\in[0,1],\ \beta(x) > 0  \right\}$, then $x_0^\star$ is unique and is characterised by the equation
	\[
	\left(\frac{1-\gamma}{\gamma}\right)\frac{\phi(t)^\frac{1}{1-\gamma}\beta(x_0^\star)}{\left(\frac{\partial K}{\partial c}(t,A(t,x_0^\star))\right)^{\frac{\gamma}{1-\gamma}}} = f(x_0^\star)H.
	\]
	The same result holds in the case $\gamma\in(-\infty, 0)$ if $f$ is non--decreasing and $\beta$ is decreasing over $L$.
	
\end{theorem}

\subsection{An explicit example}\label{sec:example}
We insist on the fact that the tariff $p$ defined in Theorem \ref{th:main} is $u-$convex by definition, and it is finite since it is written as a supremum of a continuous function over a compact set. In order to verify that $p\in\Pc$, one therefore only needs to make sure that $p^\star$ is indeed the $u-$transform of $p$ (which is the case if $p^\star$ is $u-$convex) and satisfies the other required properties. We will consider here a simplified framework where all the computations can be done almost explicitly.

\begin{assumption}\label{assump:k}
	The cost function $K$ is given, for some $n>1$, by
	\[
	K(t,c):=k(t)\frac{c^n}{n},\  (t,c)\in[0,T]\times\R_+,
	\]
	for some map $k:[0,T]\longrightarrow \R_+^\star$. Moreover, the distribution of the type of Agents is uniform, that is $f(x)=1$, and we impose $g_\gamma(x):=x{\bf 1}_{\gamma\in(0,1)}+(1-x){\bf 1}_{\gamma<0}$, for every $x\in[0,1]$.
\end{assumption}

Under Assumption \ref{assump:k}, we then have
\[
A(t,x_0)=\left(\frac{\phi(t)}{k(t)}\right)^{\frac{1}{n-\gamma}}\ell^{\frac{1-\gamma}{n-\gamma}}(x_0),
\]
and the maximisation problem becomes
\[
\widetilde U_P=\sup_{ x_0\in[0,1]}\left\{ \left(  \frac{1}{\gamma} - \frac{1}{n}\right)\int_0^T\left(\frac{\phi(t)^n}{k(t)^\gamma}\right)^{\frac{1}{n-\gamma}}\mathrm{d}t\ \ell(x_0)^{\frac{n(1-\gamma)}{n-\gamma}}+(x_0-1)H\right\}.
\]
Define
\[
B_\gamma(T):= \left(  \frac{1}{\gamma} - \frac{1}{n}\right)\int_0^T\left(\frac{\phi(t)^n}{k(t)^\gamma}\right)^{\frac{1}{n-\gamma}}\mathrm{d}t,~ \Phi(x_0):=B_\gamma(T)\ell(x_0)^{\frac{n(1-\gamma)}{n-\gamma}}+(x_0-1)H,
\]
where we emphasise that since $n>1$, when $\gamma\in(0,1)$, we easily have that $B_\gamma(T)>0$, while $B_\gamma(T)<0$ when $\gamma<0$. Furthermore, we remind the reader that when $\gamma>0$, the reservation utility of the Agents is necessarily non--negative, while it has to be negative when $\gamma<0$, since the utility function itself is negative.

\medskip
In this setting, the sufficient conditions for the $u-$convexity of the solution $p^\star$ to the relaxed problem are satisfied. Our result therefore rewrites in this case
\begin{theorem}\label{th:mainex}
	Let Assumptions \ref{assump:power}, \ref{assump:hconst} and \ref{assump:k} hold. 
	
	\medskip
	$(i)$ If $\gamma\in(0,1)$ then, the optimal tariff $p\in\Pc$ is given for any $(t,c)\in[0,T]\times\R_+$ by
	
	
	\[p(t,c)=\displaystyle\phi(t)\frac{c^\gamma}{2\gamma}+\left(\left(\frac{\phi(t)}{2}\right)^{\frac{1}{1-\gamma}}\frac{1-\gamma}{\gamma M(t)}\right)^{\frac{1-\gamma}{\gamma}}c-\frac HT+M(t)(2x_0^\star -1)^{\frac{1}{1-\gamma}},
	\]
	where
	\[M(t)=\frac{1-\gamma}{2\gamma} \left( \frac{2(2-\gamma)}{1-\gamma}\right)^{\frac{\gamma(n-1)}{n-\gamma}}\left(\frac{\phi^n(t)}{k^\gamma(t)}\right)^{\frac{1}{n-\gamma}}\left(1-(2x_0^\star -1)^\frac{2-\gamma}{1-\gamma}\right)^{-\frac{\gamma(n-1)}{n-\gamma}},\]
	and where $x_0^\star$ is the unique solution in $(1/2,1)$ of the equation
	\[H=2nA_\gamma(T)\frac{2-\gamma}{n-\gamma}(2x_0^\star-1)^{\frac{1}{1-\gamma}}\left(1-(2x_0^\star-1)^{\frac{2-\gamma}{1-\gamma}}\right)^{-\frac{\gamma(n-1)}{n-\gamma}}.\]
	Furthermore, only the Agents of type $x\geq x_0^\star$ will accept the contract.
	
	\medskip
	$(ii)$ If $\gamma<0$, then the optimal tariff $p\in\Pc$ is given for any $(t,c)\in[0,T]\times\R_+$ by
	\begin{align*}
	p(t,c)=	\displaystyle -\gamma c\left(-\frac{\phi(t)}{\gamma}\right)^{\frac1\gamma}\left(\frac{1-\gamma}{\widehat M(t)}\right)^{\frac{1-\gamma}{\gamma}}- \frac HT-\widehat M(t)(1-\widehat x_0^\star)^{\frac{1}{1-\gamma}},
	\end{align*}
	where
	\[\widehat M(t)=-\frac{1-\gamma}{\gamma} \left( \frac{2-\gamma}{1-\gamma}\right)^{\frac{\gamma(n-1)}{n-\gamma}}\left(\frac{2^\gamma\phi^n(t)}{k^\gamma(t)}\right)^{\frac{1}{n-\gamma}}(1-\widehat x_0^\star)^{-\frac{\gamma(2-\gamma)(n-1)}{(n-\gamma)(1-\gamma)}},\]
	and where
	\[\widehat x_0^\star:=\bigg(1-\left(\frac{n-\gamma}{n(1-\gamma)B_\gamma(T)}H\right)^{\frac{n-\gamma}{n(1-\gamma)+\gamma}}\left(\frac{2-\gamma}{1-\gamma}\right)^{\frac{-\gamma(n-1)}{n(1-\gamma)+\gamma}} 2^\frac{-n}{n(1-\gamma)+\gamma}\bigg)^+.\]
	Furthermore, only the Agents of type $x\geq \widehat x_0^\star$ will accept the contract.
\end{theorem}

\section{Type--dependent reservation utilities}\label{sec:general}
In this section, we study the case where the reservation utility $H$ is a general continuous and non--decreasing function of the type $x\in[0,1]$. This case strongly differs from the previous section as we can no longer guarantee that the set of Agents signing the contract with the Principal is in general an interval. The best we can say is that, under appropriate conditions on $H$, the set $X^\star(p)$ is indistinguishable from a countable union of open intervals. However, this fact does not allow us to reduce the dimensionality of the problem of the Principal and even the existence of a solution to it is not guaranteed. For this reason, we need to impose some additional structure to the set of admissible tariffs in order to obtain a well--posed problem. Specifically, we will consider a new set of admissible tariffs which is contained in a reflexive Banach space and we will use classical results from functional analysis to prove the existence of solutions to the Principal's problem. With that purpose in mind, we introduce the following Sobolev--like spaces.

\begin{definition}\label{def:sobolev-space}
	For any $\ell\geq 1$ and any open subset $\Oc$ of $X$, we denote by $W_x^{1,\ell}(\Oc)$ the space of maps $q:[0,T]\times\Oc\longrightarrow\R$ for which there exists a null set $\Nc(q)\subset[0,T]$ $($for the Lebesgue measure$)$ satisfying that for every $t\in[0,T]\setminus\Nc(q)$ the map $x\longmapsto q(t,x)$ belongs to $W^{1,\ell}(\Oc)$\footnote{That is to say the usual Sobolev space of maps admitting a weak first order derivative.} and such that
	\[
	\|q\|_{\ell,\Oc} := \left(\int_0^T \int_{\Oc} |q(t,x)|^\ell \mathrm{d}x\mathrm{d}t\right)^\frac1\ell + \left(\int_0^T \int_{\Oc} \bigg| \frac{\partial q}{\partial x}(t,x) \bigg|^\ell \mathrm{d}x\mathrm{d}t\right)^\frac1\ell < \infty.
	\]
\end{definition} 

\begin{remark}
	For the rest of the paper, for every map $q$ belonging to some space $W_x^{1,\ell}(\Oc)$, the set $\Nc(q)$ will make reference to the one mentioned in Definition \ref{def:sobolev-space}.
\end{remark}

For all the analysis of this section, we fix a number $m>1$ such that $m\gamma<1$. We are now ready to give our new definition of admissible tariffs.
\begin{definition}
	A tariff $p:[0,T]\times\R_+\longrightarrow \R$ is said to be admissible $($in the case when $H$ is not constant$)$, denoted by $p\in \widehat{\mathcal{P}}$, if in addition to Definition \ref{def:admissible tariff}, it satisfies that $p^\star\in W_x^{1,m}(\overset{o}{X})$. 
\end{definition}  
In this new setting, the Principal offers a tariff $p\in\widehat{\mathcal{P}}$ which solves her maximisation problem
\begin{equation}
\widehat{U}_P:=	\sup_{p\in\widehat{\mathcal{P}}} \int_{0}^{T} \bigg[ \int_{X^\star(p^\star) } p(t,c^\star(t,x)) f(x) \mathrm{d}x  - K\bigg(t,\int_{X^\star (p^\star)} c^\star(t,x) f(x)\mathrm{d}x\bigg)  \bigg] \mathrm{d}t.
\label{Principal_obj_sobolev}
\end{equation}
Following the previous sections, we will consider the problem $\overline U_P\geq \widehat{U}_P$, in which we drop the $u-$convexity property, defined by
\begin{align}
 \overline U_P=&\sup_{p^\star\in \widehat{C}^+ } \int_{0}^{T}  \bigg[ \int_{X^\star(p^\star) }  \left( \frac{g_\gamma(x)}{g'_\gamma(x)}\frac{\partial p^\star }{\partial x} (t,x)  - p^\star (t,x)\right)   f(x) \mathrm{d}x  - K\bigg(t,  \int_{X^\star(p^\star) }  \left( \frac{\gamma}{\phi (t)g'_\gamma(x)} \frac{\partial p^\star }{\partial x} (t,x) \right)^{\frac{1}{\gamma}} f(x)\mathrm{d}x\bigg)  \bigg] \mathrm{d}t,
\label{Principla_obj2_sobolev}
\end{align}
where $\widehat{C}^+= C^+\cap W_x^{1,m}(\overset{o}{X})$. We aim at solving the relaxed problem $\overline U_P$ and give conditions under which its solution coincides with the solution to $\widehat{U}_P$.  To facilitate the reading, we provide the main steps we will follow to achieve our results.

\begin{itemize}
	\item The structure of $X^\star$ is determined. Whenever $H$ is non--decreasing, we prove that $X^\star$ is a countable union of intervals. If $H$ is in addition concave, then $X^\star$ is of the form $\left[0,b_0\right] \cup \left[a_0,1\right]$, meaning that the Principal selects the most and least efficient Agents. 
	\item In the concave case, using the structure of $X^\star$, we rewrite the objective of the Principal as a finite--dimensional function of $a_0$ and $b_0$. The maximum of this function can be determined by standard optimisation techniques.  
	\item We finally verify whether the solution $p^\star$ for problem $\overline U_P$ satisfies the conditions of the initial problem $\widehat{U}_P$, \textit{i.e.}, it is admissible. By doing so, we can conclude that the solutions of both problems coincide and it is given by $p^\star$. 
\end{itemize}

Thus, we move to the reservation utility function $H$, recalling that it determines the structure of the set $X^\star(p^\star)$. In order to avoid complex forms of this set we make the following assumption on $g$, $H$ and $f$.

\begin{assumption}\label{assump:ghf}
	The functions $g$ and $H$ are such that for every $x\in [0,1]$
	\begin{equation} \label{hg}
	\frac{g_\gamma(x)}{g_\gamma'(x)} \leq \frac{H(x)}{H'(x)}.
	\end{equation}
	Moreover, the following maps
	\begin{align*}
	& v_1(x):=g_\gamma'(x)\bigg(\frac{\left[ {g_\gamma(x)}f(x)+{g_\gamma'(x)}F(x)\right]^+ }{f(x)}\bigg)^{\frac{\gamma}{1-\gamma}}, \;  v_2(x):=g_\gamma'(x)\bigg(\frac{\left[ {g_\gamma(x)}f(x)+{g_\gamma'(x)}F(x)-g_\gamma'(x)\right]^+ }{f(x)}\bigg)^{\frac{\gamma}{1-\gamma}},
	\end{align*}
	are non--decreasing on $[0,1]$.
\end{assumption}

\begin{remark}
	Condition \eqref{hg} is equivalent to the elasticity of reservation utility being less than the elasticity of willingness to pay for consumption. For instance, in the case $\gamma\in(0,1)$, it is automatically satisfied when $H$ is constant, and if $g_\gamma(x)=x$ then \eqref{hg} reduces to $H$ being concave. For $g_\gamma(x)=x$, $v_1$ and $v_2$ are increasing if the distribution of Agents is uniform and in some cases of the Beta distribution.\footnote{For the Beta distribution, $f(x)=C_{\alpha,\beta}x^{\alpha-1}(1-x)^{\beta-1}$, and the condition is satisfied if for instance either $\alpha=1$ or $\beta=1$
} Similarly, in the case $\gamma<0$, \eqref{hg} holds if $g_\gamma(x)=1-x$ and $H(x)=x^\alpha$ with $\alpha>1$. On the other hand, if for instance, $f(x)=1$ and $g(x)=x^\alpha$ with $\alpha\in(0,1]$, then $v_1$, $v_2$ are increasing when $\alpha \geq 1-\gamma$. 
\end{remark}

The following proposition shows that when Condition \eqref{hg} holds, it is actually never optimal for the Principal to propose a tariff for which the utility of the Agents is exactly their reservation utility on a set with positive Lebesgue measure. 
\begin{proposition} \label{prop:P=H}
	Let Assumptions \ref{assump:power} and \ref{assump:ghf} hold, and let $p^\star \in \widehat{C}^{+}$ be any function such that the set
	\[
	Y^\star(p^\star) : = \left\{ x\in[0,1],~  P^\star(x)=H(x) \right\},
	\]
	has positive Lebesgue measure. Then $p^\star $ is not optimal for problem \eqref{Principla_obj2_sobolev}.
\end{proposition}

We can now split the problem of the Principal into subintervals. Thanks to the previous proposition, we can consider without loss of generality functions $p^\star \in \widehat{C}^{+}$ such that the Lebesgue measure of $Y^\star(p^\star) $ is zero. For these functions, we define the set
\[
\widehat{X}^\star (p^\star):=X^\star(p^\star) \setminus Y^\star(p^\star)  = \left\{ x\in[0,1],\ P^\star (x) > H(x) \right\},
\]
which by continuity is an open subset of $[0,1]$. As $X^\star$ and $p^\star$ are continuous, we can replace all the integrals over $X^\star(p^\star)$ by integrals over $\widehat{X}^\star(p^\star)$ and we can write the latter set as a countable union of open disjoint intervals, that is
\[
\widehat{X}^\star (p^\star) := [0,b_0) \cup \bigcup_{n\geq 1} (a_n,b_n) \cup  (a_0,1],
\]
for some $a_0\in(0,1]$, $b_0\in[0,1)$, and $0<a_n<b_n<1$, $\forall n\geq 0$. We denote $a:=(a_n)_{n\geq 0}$, $b:=(b_n)_{n\geq 0}$ and define ${\mathcal{A}}$ as the set of such that pairs $(a,b)$. Similar to the previous section, we will use the characterisation of the set $X^\star(p^\star)$ to reformulate the problem of the Principal. For each possible set $X^\star(p^\star)$, we will solve the sub--problem in which the set of tariffs reduces to the ones for which the set of Agents accepting the contract is exactly $X^\star(p^\star)$. For any $(a,b)\in\Ac$, we define the set
\[
X^\star(a,b) = [0,b_0) \cup \bigcup_{n\geq 1} (a_n,b_n) \cup  (a_0,1].
\]
We can therefore write
\begin{align}\label{problem-ab-noibp}
\nonumber &\overline U_P=\sup_{(a,b)\in \mathcal{A}} \sup_{p^\star\in C^+(a,b)} \int_{0}^{T}  \bigg[ \int_{X^\star(a,b)}  \left( \frac{g(x)}{g'(x)}\frac{\partial p^\star }{\partial x} (t,x)  - p^\star (t,x)\right)   f(x) \mathrm{d}x \\ 
&\hspace{9em} - K\bigg(t,  \int_{X^\star(a,b)}  \left( \frac{\gamma}{\phi (t)g'(x)} \frac{\partial p^\star }{\partial x} (t,x) \right)^{\frac{1}{\gamma}} f(x)\mathrm{d}x\bigg)  \bigg] \mathrm{d}t,
\end{align}
where $C^+(a,b)$ is given by all the maps $p^\star\in \widehat{C}^+$ such that $\widehat{X}^\star (p^\star) = X^\star(a,b)$. 
\begin{remark}\label{remark-an-binding}
	The case $b_0=0$ stands for $P^\star(b_0) < H(b_0)$ and the case $a_0=1$ stands for $P^\star(a_0) < H(a_0)$. By continuity we have $P^\star(a_n) = H(a_n)$ and $P^\star(b_n) = H(b_n)$ for every $n\geq 1$.
\end{remark}
For fixed $(a,b)\in\Ac$, define the operator $\Psi_{(a,b)}:C^+(a,b)\longrightarrow\R$ by		
\begin{align*}
\nonumber \Psi_{(a,b)}(p^\star):=& \int_{0}^{T}  \bigg[ \int_{X^\star(a,b)}  \left( \frac{g_\gamma(x)}{g_\gamma'(x)}\frac{\partial p^\star }{\partial x} (t,x)  - p^\star (t,x)\right)   f(x) \mathrm{d}x \\ 
&\hspace{6em} - K\bigg(t,  \int_{X^\star(a,b)}  \left( \frac{\gamma}{\phi (t)g_\gamma'(x)} \frac{\partial p^\star }{\partial x} (t,x) \right)^{\frac{1}{\gamma}} f(x)\mathrm{d}x\bigg)  \bigg] \mathrm{d}t.
\end{align*}
As previously explained, we focus on the problem in which the set of Agents signing the contract is fixed. Define 
\[
(P_{a,b})~~~\sup_{p^\star\in C^+(a,b)} \Psi_{(a,b)}(p^\star).
\]
Our first result gives the existence of a solution to the above infinite-dimensional optimisation problem, and requires the following assumption, which involves mainly the cost function. It is required in order to obtain nice coercivity properties
\begin{assumption}\label{assump:cost}
	The cost function $K$ satisfies the following growth condition
	\[
	K(t,c) \geq k(t) c^n, ~\forall c\in\Cc,
	\]
	where the map $k:[0,T]\longrightarrow\R_+$ is bounded from below by some constant $\underline{k}>0$ and $n\geq1$. Moreover, we have that
	\[
	I:=\inf\left\{ k(t)\left( \frac{\gamma f(x)}{\phi(t)g_\gamma'(x)} \right)^n,\ (t,x)\in[0,T]\times[0,1]\right\}>0.\]
\end{assumption}
We now have
\begin{proposition}\label{prop:existence}
	Let Assumptions \ref{assump:power}, \ref{assump:ghf} and \ref{assump:cost} hold. For every $(a,b)\in\Ac$, the optimisation problem $(P_{a,b})$ has at least one solution.
\end{proposition}

Next, we obtain necessary optimality conditions for problem $(P_{a,b})$. Recalling from Remark \ref{remark-an-binding} that the (IR) condition is binding at each $a_n,b_n$, by integration by parts we can rewrite $\Psi_{(a,b)}$ as
\begin{align}
\nonumber \Psi_{(a,b)}(p^\star) =  &  \int_{0}^{T}  \Bigg(\int_{0}^{b_0}  \frac{\left( g_\gamma(x)f(x) + g_\gamma'(x)F(x) \right)}{g_\gamma'(x)} \frac{\partial p^\star}{\partial x}(t,x) \mathrm{d}x+ \int_{a_0}^1  \frac{\left( g_\gamma(x)f(x) + g_\gamma'(x)F(x) - g_\gamma'(x) \right)}{g_\gamma'(x)} \frac{\partial p^\star}{\partial x}(t,x) \mathrm{d}x \Bigg)\mathrm{d}t \\
&\nonumber + \sum_{n=1}^\infty \int_0^T \int_{a_n}^{b_n} \frac{\left( g_\gamma(x)f(x) + g_\gamma'(x)F(x) \right)}{g'(x)} \frac{\partial p^\star}{\partial x}(t,x) \mathrm{d}x \mathrm{d}t   \\
&\nonumber  - K\bigg(t,  \int_{(0,b_0) \cup \bigcup_{n\geq 1} (a_n,b_n) \cup  (a_0,1)}  \left( \frac{\gamma}{\phi (t)g_\gamma'(x)} \frac{\partial p^\star }{\partial x} (t,x) \right)^{\frac{1}{\gamma}} f(x)\mathrm{d}x\bigg)  \mathrm{d}t\\ 
&  +\sum_{n=1}^\infty F(a_n) H(a_n)- \sum_{n=1}^\infty F(b_n) H(b_n)-F(b_0)H(b_0) + (F(a_0)-1)H(a_0).
\end{align}
To simplify notations, denote 
\[
A(t,a,b):=\int_{X^\star(a,b)}  \left( \frac{\gamma}{\phi (t)g_\gamma'(x)} \frac{\partial p^\star }{\partial x} (t,x) \right)^{\frac{1}{\gamma}} f(x)\mathrm{d}x.
\]

\begin{theorem}\label{thr-optimality-condition-concave-case}
	Let Assumptions \ref{assump:power}, \ref{assump:ghf} and \ref{assump:cost} hold and let $p^\star$ be a solution of $(P_{a,b})$. Consider an interval $I=(x_\ell,x_r)\subseteq X^\star(a,b)$ such that $P^\star(x)>H(x)$ for every $x\in I$, $P^\star(x_\ell)=H(x_l)$ and $P^\star(x_r)=H(x_r)$. Then there exists a null set $\Nc\subset[0,T]$ and a constant $\mu_t$ for every $t\in[0,T]\setminus\Nc$ such that the following optimality condition is satisfied
	\begin{enumerate}
		\item[$(i)$] In the case $I\subseteq(a_0,1)$, for every $x\in I$ we have
		\begin{equation}\label{eq-optimality-concave-h-a0}
		\frac{\partial p^\star}{\partial x}(t,x) = \displaystyle\left(\frac{\phi(t)^{\frac1\gamma}\left[ g_\gamma(x)f(x)+g'_\gamma(x)F(x)-g'_\gamma(x)+g'_\gamma(x)\mu_t\right]^+}{f(x)\frac{\partial K}{\partial c}\left(t,A(t,a,b)\right)}\right)^{\frac{\gamma}{1-\gamma}} \frac{g'_\gamma(x)}{\gamma}.
		\end{equation}
		\item[$(ii)$] In the case $I\subseteq(0,b_0)\cup_{n\geq 1}(a_n,b_n)$, for every $x\in I$ we have
		\begin{equation}\label{eq-optimality-concave-h-an}
		\frac{\partial p^\star}{\partial x}(t,x) = \displaystyle\left(\frac{\phi(t)^{\frac1\gamma}\left[ g_\gamma(x)f(x)+g'_\gamma(x)F(x)+g'_\gamma(x)\mu_t\right]^+}{f(x)\frac{\partial K}{\partial c}\left(t,A(t,a,b)\right)}\right)^{\frac{\gamma}{1-\gamma}} \frac{g'_\gamma(x)}{\gamma}.
		\end{equation}
	\end{enumerate}
\end{theorem}

The proof of Theorem \ref{thr-optimality-condition-concave-case} consists in several technical propositions which are given and proved in Appendix \ref{app:general} below. 

\medskip
Even after solving the sub--problems $(P_{a,b})$, the main difficulty for moving back to the relaxed problem of the Principal is the infinite dimensionality of the set $\Ac$. However, equations \eqref{eq-optimality-concave-h-a0} and \eqref{eq-optimality-concave-h-an} give us some insight on the behaviour of the optimal tariff in the set of Agents accepting the contract, and how it can be used to obtain a finite--dimensional formulation of Problem $\overline U_P$.

\medskip
Under Assumption \ref{assump:ghf}, the optimal tariffs are convex over intervals where the (IR) condition is not binding. This is the last result of this section and it is a direct consequence of the fact that the functions $v_1$ and $v_2$ are non--decreasing, which makes the derivative of $p^\star$ non--decreasing as well. The convexity property will allow us to completely solve the Principal's problem in the next subsection, when the function $H$ is strictly concave, because $\Ac$ will reduce to a 2--dimensional set. Other cases in which $\Ac$ is also finite--dimensional, and therefore the problem of the Principal can be easily solved, are mentioned in the last subsection.  

\begin{proposition}\label{optimalconvexity}
	Let Assumptions \ref{assump:power}, \ref{assump:ghf} and \ref{assump:cost} hold. Let $p^\star$ be a solution to problem $(P_{a,b})$. Then $P^\star$ is convex on every interval over which $P^\star$ is strictly greater than $H$. 
\end{proposition}

\subsection{Strictly concave reservation utility}

In this section we assume the reservation utility function of the Agents is strictly concave.

\begin{assumption}\label{assump:Hconc}
	The map $x\longmapsto H(x)$ is strictly concave and non--decreasing.
\end{assumption}

The main consequence of Assumption \ref{assump:Hconc} is the following simple result, which shows that we can always restrict our attention to sets $X^\star$ where only the most effective and less effective Agents sign the contract.

\begin{proposition}\label{prop-an-bn}
	Let Assumptions \ref{assump:power}, \ref{assump:ghf}, \ref{assump:cost} and \ref{assump:Hconc} hold. Let $(a,b)\in\Ac$ be such that $0<a_{n_0}<b_{n_0}<1$ for some $n_0\geq 1$. Then the solution to problem $(P_{a,b})$ is not optimal for problem \eqref{Principla_obj2_sobolev}.
\end{proposition}

\medskip
In the rest of this section, we start by deriving a general solution under some implicit assumptions, and then show that the latter can be verified in an example in the context of Assumption \ref{assump:k}.

\subsubsection{The general tariff}

\begin{figure}
	\centering
	\begin{tikzpicture}
	[/pgf/declare function={ H=2*sqrt(x+.5)+.5*x; p=.05*x*x+.005*x+4;}]
	\begin{axis}[
	axis lines=left,    
	domain=0:22, samples=100, 
	xmin=0, xmax=30, ymin=0, ymax=30, 
	no markers, 
	unit vector ratio*=1 1 1, 
	xtick=\empty, ytick=\empty, 
	cycle list={} 
	]
	\addplot [ line width=1pt, blue] {max(p,H)};
	\addplot [ line width=1pt, red] {H};
	\end{axis}
	
	\foreach \a/\b/\c [count=\i] in { .4/.8/{b}, 2.85/2.9/{a}} { 
		\begin{scope}
		\draw[dashed] (\a,0) -- (\a,\b); 
		\node [below] at (\a,0) {$\c_0$};
		\end{scope}   
	}
	\draw (4,5.5) node[left] {$P^\star (x)$};
	\draw (4.2,3.8) node[right] {$H(x)$};
	\end{tikzpicture}
	\caption{$\widehat X^\star(p^\star) $ for strictly concave $H$.}  \label{fighconcave}
\end{figure}
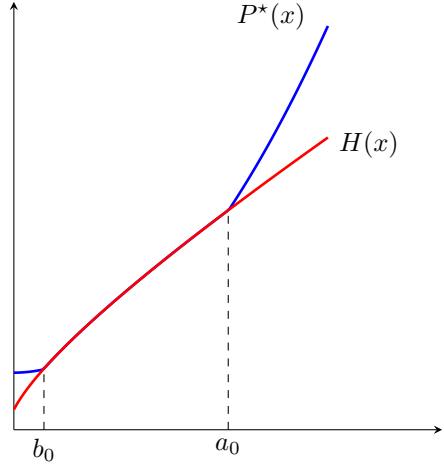

Proposition \ref{prop-an-bn} implies that the solution of \eqref{problem-ab-noibp} is attained at some $p^\star$ satisfying $\widehat{X}^\star(p^\star)  = [0,b_0) \cup (a_0,1].$ We expect then the optimal tariff to look like the curve in Figure \ref{fighconcave}. Let us then define the set
\[\Ac_2:=\left\{(a,b)\in[0,1]^2,\ b\leq a\right\}.\]
Theorem \ref{thr-optimality-condition-concave-case} gives us only partial information about the solution of the problem $(P_{a,b})$, for $(a,b)\in\Ac_2$. Now that we assume in addition that $H$ is strictly concave, we can actually precise the necessary optimality conditions with the following proposition, which tells us that the value of the constants $\mu_t$ is zero in the intervals of the form $[0,b)$ and $(a,1]$. 

\begin{proposition}\label{prop:mu_t=0}
	Let Assumptions \ref{assump:power}, \ref{assump:ghf}, \ref{assump:cost} and \ref{assump:Hconc} hold. Let $p^\star$ be a solution of $(P_{a,b})$, for $(a,b)\in\Ac_2$, and $I$ be as in Theorem \ref{thr-optimality-condition-concave-case}.  Then there exists a null set $\Nc'\subset[0,T]$ such that for every $t\in[0,T]\setminus\Nc'$ the optimality conditions from Theorem \ref{thr-optimality-condition-concave-case} hold with $\mu_t=0$.
\end{proposition} 

\medskip
Following the computations of Section \ref{section-constant}, we define
\begin{align*}
A(t,a_0,b_0):= & \ g_K^{(-1)}  \bigg( \phi(t)^{ \frac{1}{1-\gamma}} \int_0^{b_0}\bigg(\frac{ \left[ g(x)f(x)+g'(x)F(x)\right]^+ }{f^\gamma(x)}\bigg)^{\frac{1}{1-\gamma}}\mathrm{d}x   \\
& \hspace{4em} + \phi(t)^{ \frac{1}{1-\gamma}}\int_{a_0}^1\bigg(\frac{ \left[ g(x)f(x)+g'(x)F(x)-g'(x) \right]^+ }{f^\gamma(x)}\bigg)^{\frac{1}{1-\gamma}}\mathrm{d}x  \bigg).
\end{align*}
The aim of the next proposition is similar in spirit to that of Proposition \ref{prop-an-bn}, in the sense that it allows to exclude many specifications of $(a,b)\in\mathcal A_2$, for which the solution to the sub--problem $(P_{a,b})$ does not solve the general relaxed problem $\overline U_P$. 
\begin{proposition}\label{prop:great}
	Let Assumptions \ref{assump:power}, \ref{assump:ghf}, \ref{assump:cost} and \ref{assump:Hconc} hold. Let $p^\star$ be a solution of $(P_{a,b})$. If either 
	\[
	\Xi_\gamma(a_0,b_0):=\displaystyle\int_0^T \Bigg(\frac{\phi(t)^{\frac1\gamma}\left[ g_\gamma(a_0)f(a_0)+g'_\gamma(a_0)F(a_0)-g'_\gamma(a_0)\right]^+}{f(a_0)\frac{\partial K}{\partial c}\left(t,A(t,a_0,b_0)\right)}\Bigg)^{\frac{\gamma}{1-\gamma}} \frac{g'_\gamma(a_0)}{\gamma}\mathrm{d}t < H'(a_0), 
	\]
	or 
	\[
	\Psi_\gamma(a_0,b_0):=\displaystyle \int_0^T \Bigg(\frac{\phi(t)^{\frac1\gamma}\left[ g_\gamma(b_0)f(b_0)+g'_\gamma(b_0)F(b_0)\right]^+}{f(b_0)\frac{\partial K}{\partial c}\left(t,A(t,a_0,b_0)\right)}\Bigg)^{\frac{\gamma}{1-\gamma}} \frac{g'_\gamma(b_0)}{\gamma}\mathrm{d}t > H'(b_0),
	\]
	then the solution to problem $(P_{a,b})$ is not optimal for problem \eqref{Principla_obj2_sobolev}.
\end{proposition}

\medskip
Judging by the results of Proposition \ref{prop:great}, it is natural to define $\Ac_2^\prime$ as the set of all the pairs $(a,b)\in\Ac_2$ for which 
\begin{align*}
\Xi_\gamma(a_0,b_0) \geq H'(a_0), \; \Psi_\gamma(a_0,b_0)\leq H'(b_0).
\end{align*}

Thanks to Proposition \ref{prop:great}, we have thus reduced problem $\overline U_P$ to
\begin{equation*}
\overline U_P = \sup_{ (a_0,b_0)\in\Ac_2^\prime}  \int_{0}^{T} \left[ \frac{\phi(t)^\frac{1}{1-\gamma}\ell(a_0,b_0)}{\gamma \left(\frac{\partial K}{\partial c}(t,A(t,a_0,b_0))\right)^\frac{\gamma}{1-\gamma}} - K\left( t,\frac{\phi(t)^\frac{1}{1-\gamma}\ell(a_0,b_0)}{\left(\frac{\partial K}{\partial c}(t,A(t,a_0,b_0))\right)^\frac{1}{1-\gamma}}\right) \right]\mathrm{d}t + \theta(a_0,b_0),
\end{equation*}
where we abused notations and defined the corresponding functions 
\begin{align*}
\ell(a_0,b_0) :=  & \int_0^{b_0}\left(\frac{ \left[ g(x)f(x)+g'(x)F(x)\right]^+ }{f^\gamma(x)}\right)^{\frac{1}{1-\gamma}}\mathrm{d}x  + \int_{a_0}^1\left(\frac{ \left[ g(x)f(x)+g'(x)F(x)-g'(x) \right]^+ }{f^\gamma(x)}\right)^{\frac{1}{1-\gamma}}\mathrm{d}x,\\
\theta(a_0,b_0) := & -F(b_0)H(b_0) + (F(a_0)-1)H(a_0).
\end{align*}
Since all these maps are continuous on $[0,1]^2$, the supremum over the compact set above is attained at some $(a_0^\star,b_0^\star)\in\Ac_2^\prime$. We have therefore proved our main result of this section
\begin{theorem}\label{th:main2}
	Let Assumptions \ref{assump:power}, \ref{assump:ghf}, \ref{assump:cost} and \ref{assump:Hconc} hold. We have
	
	\medskip
	$(i)$ The maximum in \eqref{Principla_obj2_sobolev} is attained for the map 
	\begin{align*}
	p^\star(t,x)=\begin{cases}
	\displaystyle  \frac{H(b_0^\star)}{T} -\displaystyle\frac{\phi(t)^{\frac{1}{1-\gamma}} }{\gamma\left(\frac{\partial K}{\partial c}(t,A(t,a^\star_0,b^\star_0))\right)^{\frac{\gamma}{1-\gamma}}}\int_{x}^{b_0^\star} v_1(y)dy,\ \textrm{if } x\in [0,b_0^\star), \\[1.8em]
	\tilde p^\star(t,x), \text{ if }x\in[b_0^\star,a_0^\star],\\[0.8em]
	\displaystyle \frac{H(a_0^\star)}{T} +\frac{\phi(t)^{\frac{1}{1-\gamma}} }{\gamma\left(\frac{\partial K}{\partial c}(t,A(t,a^\star_0,b^\star_0))\right)^{\frac{\gamma}{1-\gamma}}}\int_{a_0^\star}^xv_2(y)dy,\  \textrm{if } x\in (a_0^\star,1],
	\end{cases}
	\end{align*}
	where $\tilde p^\star(t,x)$ is any continuous and non-decreasing map $($with respect to $x)$ such that 
	\[\int_0^T\tilde p^\star(t,b_0^\star)\mathrm{d}t=H(b_0^\star),\ \int_0^T\tilde p^\star(t,a_0^\star)\mathrm{d}t=H(a_0^\star),\ \int_0^T\tilde p^\star(t,x)\mathrm{d}t<H(x),\ \text{for all }x\in(b_0^\star,a_0^\star).\]

	$(ii)$ Define $p$, for any $(t,c)\in[0,T]\times\R_+$, by
	\[p(t,c):=\underset{x\in [0,1]}{\sup}\left\{g_\gamma(x)\phi(t)\frac{c^\gamma}{\gamma}-p^\star(t,x)\right\}.\]
	If $p^\star$ is $u-$convex on $X^\star(p^\star)$, then $p$ is the optimal tariff for the problem \eqref{Principal_obj_sobolev}. Furthermore, the Principal only signs contracts with the Agents of type $x\in[0,b_0^\star]\cup[a_0^\star,1]$.
\end{theorem}

\subsubsection{Power type cost function}
Exactly as in the case where $H$ was independent of $x$, the computations become much simpler as soon as Assumption \ref{assump:k} holds. Let us note $R_\gamma(a_0,b_0) = 1 + (2b_0)^\frac{2-\gamma}{1-\gamma} - \left((2a_0-1)^+\right)^\frac{2-\gamma}{1-\gamma}$ if $\gamma\in(0,1)$ and  $R_\gamma(a_0,b_0) = 1 -((1-2b_0)^+)^\frac{2-\gamma}{1-\gamma} + \left(2 - 2a_0\right)^\frac{2-\gamma}{1-\gamma}$ if $\gamma<0$. Then, the functions $\ell$ and $A$ are given, for any $(t,a_0,b_0)\in[0,T]\times\Ac_2^\prime$, by 
\begin{align*}
\ell_\gamma(a_0,b_0)  = \frac{1-\gamma}{2(2-\gamma)}R_\gamma(a_0,b_0),\ A(t,a_0,b_0)  = \left( \frac{\phi(t)}{k(t)}  \right)^\frac{1}{n-\gamma} \ell(a_0,b_0)^\frac{1-\gamma}{n-\gamma}.
\end{align*}
So, in order to obtain $(a_0^\star, b_0^\star)$ we have to solve
\begin{equation}\label{eq:prob}
\sup_{(a_0,b_0)\in\Ac_2^\prime} \left( \frac{1}{\gamma}-\frac{1}{n}\right) \int_0^T \left( \frac{\phi(t)^n}{k(t)^\gamma}  \right)^\frac{1}{n-\gamma}\mathrm{d}t ~\ell(a_0,b_0)^\frac{n(1-\gamma)}{n-\gamma} -b_0H(b_0) + (a_0-1)H(a_0).
\end{equation}
Let us now compute the associated tariff $p$ and check that $p$ indeed belongs to $\Pc$ and that its $u-$transform is $p^\star$. Fix some $t\in[0,T]$ and define
\begin{align*}
N_\gamma&:=
\displaystyle\frac{2^{\frac{\gamma}{1-\gamma}}(1-\gamma)}{\gamma} \left( \frac{2(2-\gamma)}{1-\gamma}\right)^{\frac{\gamma(n-1)}{n-\gamma}}\left(\frac{\phi^n(t)}{k^\gamma(t)}\right)^{\frac{1}{n-\gamma}}R_\gamma(a_0,b_0)^{-\frac{\gamma(n-1)}{n-\gamma}}.
\end{align*} 
Recall that by Proposition \ref{prop:great}, the following inequalities must be satisfied
\begin{itemize}
	\item[$(i)$] If $\gamma\in(0,1)$
	\begin{align}\label{eq:const}
	\frac{((2a_0^\star-1)^+)^{\frac{\gamma}{1-\gamma}}}{\gamma\ell(a_0^\star,b^\star_0)^{\frac{\gamma(n-1)}{n-\gamma}}}\int_0^T\left(\frac{\phi(t)^n}{k(t)^\gamma}\right)^{\frac1{n-\gamma}}\mathrm{d}t&\geq H^\prime(a_0^\star),\ \frac{\left(2b_0^\star\right)^{\frac{\gamma}{1-\gamma}}}{\gamma\ell(a_0^\star,b^\star_0)^{\frac{\gamma(n-1)}{n-\gamma}}}\int_0^T\left(\frac{\phi(t)^n}{k(t)^\gamma}\right)^{\frac1{n-\gamma}}\mathrm{d}t\leq H^\prime(b_0^\star).
	\end{align}
	\item[$(ii)$] If $\gamma<0$
	\begin{align}\label{eq:const2}
	-\frac{\left(2(1-a_0^\star)^+\right)^{\frac{\gamma}{1-\gamma}}}{\gamma\ell(a_0^\star,b^\star_0)^{\frac{\gamma(n-1)}{n-\gamma}}}\int_0^T\left(\frac{\phi(t)^n}{k(t)^\gamma}\right)^{\frac1{n-\gamma}}\mathrm{d}t&\geq H^\prime(a_0^\star),\
	-\frac{\left((1-2b_0^\star)^+\right)^{\frac{\gamma}{1-\gamma}}}{\gamma\ell(a_0^\star,b^\star_0)^{\frac{\gamma(n-1)}{n-\gamma}}}\int_0^T\left(\frac{\phi(t)^n}{k(t)^\gamma}\right)^{\frac1{n-\gamma}}\mathrm{d}t\leq H^\prime(b_0^\star).
	\end{align}
\end{itemize}
Notice in particular that when $\gamma\in(0,1)$, \eqref{eq:const} implies that $a_0^\star>1/2$, since $H$ is increasing. With similar computations as in Section \ref{sec:example}, we compute that
\begin{itemize}
	\item[$(i)$] If $\gamma\in(0,1)$
	\begin{align*}
	p^\star (t,x) =\begin{cases} 
	\displaystyle \frac{H(b_0^\star)}{T} -N_\gamma  \left( (b_0^\star)^{\frac{1}{1-\gamma}}-x^{\frac{1}{1-\gamma}}\right), \ \text{if}\ x\in[0,b_0^\star),\\[0.8em]
	\displaystyle\tilde p^\star(t,x), \text{ if }x\in[b_0^\star,a_0^\star],\\[0.8em]
	\displaystyle  \frac{H(a_0^\star)}{T} +N_\gamma \left( \left(x-\frac12\right)^{\frac{1}{1-\gamma}}-\left(a_0^\star -\frac12\right)^{\frac{1}{1-\gamma}}\right), \ \text{if}\ x\in(a_0^\star,1],
	\end{cases}
	\end{align*}
	\item[$(ii)$] If $\gamma<0$
	\begin{align*}
	p^\star (t,x) =\begin{cases} 
	\displaystyle \frac{H(b_0^\star)}{T} -N_\gamma  \left( \left(\frac12-b_0^\star\wedge\frac12\right)^{\frac{1}{1-\gamma}}-\left(\frac12-x\wedge\frac12\right)^{\frac{1}{1-\gamma}}\right), \ \text{if}\ x\in[0,b_0^\star),\\[0.8em]
	\displaystyle\tilde p^\star(t,x), \text{ if }x\in[b_0^\star,a_0^\star],\\[0.8em]
	\displaystyle  \frac{H(a_0^\star)}{T} +N_\gamma \left(\left(1-x\right)^{\frac{1}{1-\gamma}}-\left(1-a_0^\star\right)^{\frac{1}{1-\gamma}}\right), \ \text{if}\ x\in(a_0^\star,1],
	\end{cases}
	\end{align*}
\end{itemize}
\medskip
Actually, in this case, the map $p^\star$ will be $u-$convex if and only if the following implicit assumption holds.
\begin{assumption}\label{assump:hhh}
	The solutions $(a_0^\star,b_0^\star)$ of \eqref{eq:prob} are such that
	\[b_0^\star\leq a_0^\star-\frac12.\]
\end{assumption}
Our main result in this case reads.
\begin{theorem}\label{t:exam2}
	Let Assumptions \ref{assump:power}, \ref{assump:k}, \ref{assump:ghf}, \ref{assump:cost}, \ref{assump:Hconc} and \ref{assump:hhh} hold, then the optimal tariff $p\in\widehat{\mathcal{P}}$ is given for any $(t,c)\in[0,T]\times\R_+$, when $\gamma\in(0,1)$ by
	\begin{align*}
	p(t,c)=\begin{cases}
	\displaystyle\phi(t)\frac{c^\gamma}{2\gamma}+\phi(t)L_\gamma(t)^{\gamma-1}c+ N_\gamma\left(a_0^\star-\frac12\right)^{\frac{1}{1-\gamma}}-\frac{H(a_0^\star)}{T},\ \text{if}\ L_\gamma(t)\left(a_0^\star-\frac12\right)^{\frac{1}{1-\gamma}}<c,\\[0.8em]
	\displaystyle\tilde x^\star(c)\phi(t)\frac{c^\gamma}{\gamma}-\tilde p^\star(t,\tilde x(c)), \ \text{if}\ L_\gamma(t)(b_0^\star)^{\frac{1}{1-\gamma}}<c\leq L_\gamma(t)\left(a_0^\star-\frac12\right)^{\frac{1}{1-\gamma}},\\[0.8em]
	\displaystyle\phi(t)L_\gamma(t)^{\gamma-1}c-\frac{H(b_0^\star)}{T}+N_\gamma(b_0^\star)^{\frac{1}{1-\gamma}},\ \text{if}\ 0\leq c\leq L_\gamma(t)(b_0^\star)^{\frac{1}{1-\gamma}},
	\end{cases}
	\end{align*}
	and when $\gamma<0$ by
	\begin{align*}
	p(t,c)=\begin{cases}
	\displaystyle\phi(t)L_\gamma(t)^{\gamma-1}c+N_\gamma(1-a_0^\star)^{\frac{1}{1-\gamma}}-\frac{H(a_0^\star)}{T},\ \text{if}\ 0<c\leq L_\gamma(t)(1-a_0^\star)^{\frac{1}{1-\gamma}},\\[0.8em]
	\displaystyle\tilde x^\star(c)\phi(t)\frac{c^\gamma}{\gamma}-\tilde p^\star(t,\tilde x(c)), \ \text{if}\  L_\gamma(t)(1-a_0^\star)^{\frac{1}{1-\gamma}}<c\leq L_\gamma(t)\left(\frac12-b_0^\star\right)^{\frac{1}{1-\gamma}},\\[0.8em]
	\displaystyle \phi(t)\frac{c^\gamma}{2\gamma}+\phi(t)L_\gamma(t)^{\gamma-1}c+N_\gamma\left(\frac12-b_0^\star\right)^{\frac{1}{1-\gamma}}-\frac{H(b_0^\star)}{T},\ \text{if}\ L_\gamma(t)\left(\frac12-b_0^\star\right)^{\frac{1}{1-\gamma}}< c.\end{cases}
	\end{align*}
	where $L_\gamma(t):=\left(\frac{\gamma N_\gamma}{(1-\gamma)\phi(t)}\right)^{\frac 1\gamma}$
	and where $(a_0^\star,b_0^\star)$ are maximisers of
	\begin{equation*}
	\sup_{(a_0,b_0)\in\Ac_2^\prime} C(T) R_\gamma(a_0,b_0)^\frac{n(1-\gamma)}{2-\gamma}-b_0H(b_0) + (a_0-1)H(a_0).
	\end{equation*}
	Furthermore, the Principal will only choose clients with type $x\in[0,b_0^\star]\cup[a_0^\star,1]$.
\end{theorem}

\subsection{General reservation utility}

In this section, we want to point out that the assumption of the reservation utility function $H$ being strictly concave is not mandatory in order to solve problem \eqref{Principla_obj2_sobolev}. We intend to explain in which other cases we can hope to solve the problem and what procedure can be followed to do so. 

\medskip
In order to reduce the relaxed problem $\overline U_P$ to a finite dimensional problem, we need $H$ to have at most a finite number of intersecting points with a strictly convex function. If $H$ were to satisfy this property, then we would be able to prove a result similar to Proposition \ref{prop-an-bn}, and we would conclude that the optimal set $\hat{X^\star}(p^\star)$ is a finite union of intervals contained in $[0,1]$. The next step then would be to prove that the Lagrange multipliers $\mu_t$ in Theorem \ref{thr-optimality-condition-concave-case} are equal to zero, using for instance local perturbations as we did to prove Proposition \ref{prop:mu_t=0}. This would allow to solve explicitly the optimality conditions \eqref{eq-optimality-concave-h-a0} and \eqref{eq-optimality-concave-h-an} by using the corresponding auxiliary map $A(a,b)$.

\medskip 
An interesting example of a reservation utility function satisfying the finite intersecting points property is the ``constant--linear'' case, which is presented next and leads to a 3--dimensional optimisation problem.

{\rm\begin{example} Suppose that for some $\alpha, \beta\geq0$ and $x_h\in[0,1]$, $H$ has the form
	\[
	H(x) =\begin{cases}
	\beta,\ \textrm{if } x\in[0,x_h], \\
	\alpha(x-x_h) + \beta,\ \textrm{if } x\in[x_h,1].
	\end{cases}
	\]
	Such a reservation utility accounts for the fact that all the Agents, whatever their appetence for power consumption is, should at least receive a minimal level of utility, in this case $\beta$. Though in general two convex functions can intersect at countably many points, given the specific form of $H$, it can intersect an increasing and convex function at, at most, three points, as shown in Figure \ref{fighconvex}. 
	\begin{figure}[!ht]
		\centering
		\subfigure[Case $a_1=0$.]{
			\begin{tikzpicture}
			[/pgf/declare function={ cons=2; rect=1.6*(x-4); p=.05*x*x+.1*x+4;}]
			\begin{axis}[
			axis lines=left,    
			domain=0:22, samples=100, 
			xmin=0, xmax=30, ymin=0, ymax=30, 
			no markers, 
			unit vector ratio*=1 1 1, 
			xtick=\empty, ytick=\empty, 
			cycle list={} 
			]
			\addplot [ line width=1pt, blue] {max(cons,rect,p)};
			\addplot [ line width=1pt, red] {max(cons,rect)};
			\end{axis}
			
			\foreach \a/\b [count=\i] in {0/0, 2.05/2.1, 3.6/4.5} { 
				\draw[dashed] (\a,0) -- (\a,\b); \node [below] at (\a,0) {$a_\i$};
			}
			\draw[dashed] (1,0) -- (1,.33); \node[below] at (1,0){$x_h$};
			\draw (4,5.5) node[left] {$P^\star (x)$}; \draw (2.4,2.6) node[right] {$H(x)$};
			\end{tikzpicture}}
		\hspace{2cm}
		\subfigure[Case $a_1>0$.]{
			\begin{tikzpicture}
			[/pgf/declare function={ cons=2; rect=1.7*(x-6); p=1.5*1.5*.025*x*x+1.5*.005*x+1.5;}]
			\begin{axis}[
			axis lines=left,    
			domain=0:22, samples=100, 
			xmin=0, xmax=30, ymin=0, ymax=30, 
			no markers, 
			unit vector ratio*=1 1 1, 
			xtick=\empty, ytick=\empty, 
			cycle list={} 
			]
			\addplot [ line width=1pt, blue] {max(cons,rect,p)};
			\addplot [ line width=1pt, red] {max(cons,rect)};
			\end{axis}
			\foreach \a/\b [count=\i] in {.55/.4, 2.05/1.5, 3.65/4.25} { 
				\begin{scope}
				\draw[dashed] (\a,0) -- (\a,\b); \node [below] at (\a,0) {$a_\i$};
				\end{scope}   	}
			\draw[dashed] (1.35,0) -- (1.35,.33); \node[below] at (1.35,0){$x_h$};
			\draw (4.5,5) node[right] {$H(x)$}; \draw (4,5.5) node[left] {$P^\star (x)$};
			\end{tikzpicture}}
		\caption{$X^\star(p^\star) $ for a "constant-linear" $H$.}\label{fighconvex}
	\end{figure}
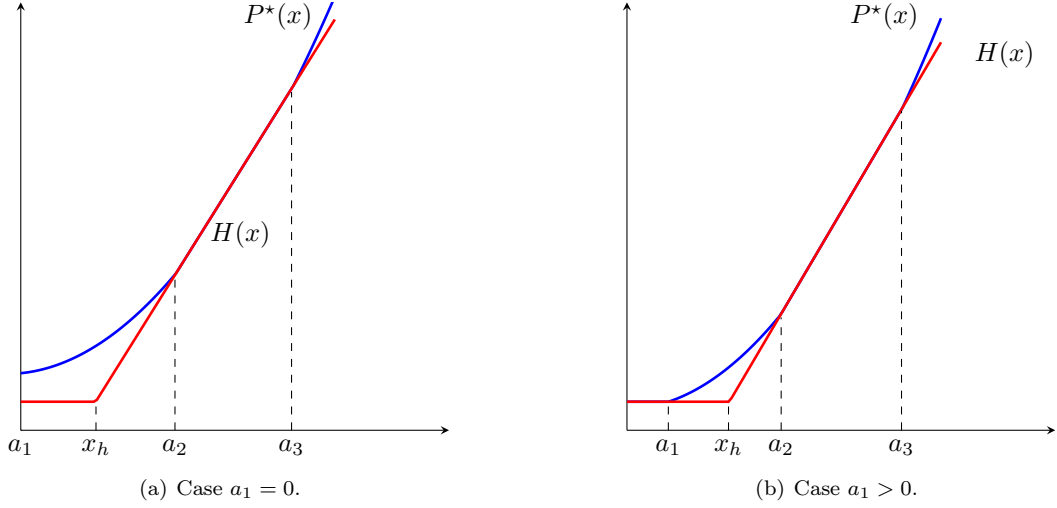
	
	\medskip
	Therefore, we deduce that $X^\star(p^\star)$ has the following form
	\[
	X^\star(p^\star)  = [a_1,a_2] \cup [a_3,1], \text{ for some }0\leq a_1\leq a_2\leq a_3\leq 1.
	\]
	We define then the set
	\[
	\Ac_3:=\left\{(a,b,c)\in[0,1]^2,\ a\leq b\leq c\right\}.
	\]
	
	After proving that the Lagrange multipliers $\mu_t$ in Theorem \ref{thr-optimality-condition-concave-case} are equal to zero, Problem $\overline U_P$ becomes, abusing notations slightly 
	\[
	\sup_{  (a_1,a_2,a_3)\in \Ac_3}\  \int_{0}^{T} \left[ \frac{\phi(t)^\frac{1}{1-\gamma}\ell(a_1,a_2,a_3)}{\gamma \left(\frac{\partial K}{\partial c}(t,A(t,a_1,a_2,a_3))\right)^\frac{\gamma}{1-\gamma}} - K\left( t,\frac{\phi(t)^\frac{1}{1-\gamma}\ell(a_1,a_2,a_3)}{\left(\frac{\partial K}{\partial c}(t,A(t,a_1,a_2,a_3)\right)^\frac{1}{1-\gamma}}\right) \right]\mathrm{d}t  + \theta(a_1,a_2,a_3), 
	\]
	where for any $(t,a_1,a_2,a_3)\in[0,T]\times\Ac_3$
	\begin{align*}
	&\ell(a_1,a_2,a_3) :=  \int_{a_1}^{a_2}\bigg(\frac{ [ g_\gamma(x)f(x)+g'_\gamma(x)F(x)]^+ }{f^\gamma(x)}\bigg)^{\frac{1}{1-\gamma}}\mathrm{d}x   + \int_{a_3}^1\bigg(\frac{ [ g_\gamma(x)f(x)+g'_\gamma(x)(F(x)-1)]^+ }{f^\gamma(x)}\bigg)^{\frac{1}{1-\gamma}}\mathrm{d}x,\\
	&\theta(a_1,a_2,a_3) :=   F(a_1)H(a_1)-F(a_2)H(a_2) + (F(a_3)-1)H(a_3),\\
	&A(t,a_1,a_2,a_3):= g_K^{(-1)}  \bigg( \phi(t)^{ \frac{1}{1-\gamma}} \int_{a_1}^{a_2}\bigg(\frac{ [ g_\gamma(x)f(x)+g'_\gamma(x)F(x)]^+ }{f^\gamma(x)}\bigg)^{\frac{1}{1-\gamma}}\mathrm{d}x   \\
	&  \hspace{10.5em} +\phi(t)^{ \frac{1}{1-\gamma}}\int_{a_3}^1\bigg(\frac{ \left[ g_\gamma(x)f(x)+g'_\gamma(x)F(x)-g'_\gamma(x) \right]^+ }{f^\gamma(x)}\bigg)^{\frac{1}{1-\gamma}}\mathrm{d}x  \bigg).
	\end{align*}
Since all the maps are continuous, the previous problem has a solution $(a_1^\star, a_2^\star, a_3^\star)$ which provides the optimal tariff. Choosing again
	$ K(t,c) = k(t)\frac{c^n}{n}$, $f(x) = 1$, $g(x)=x,$ in the case $\gamma\in(0,1)$ the functions $\ell$ and $A$ will be given by
	\begin{align*}
	\ell(a_1,a_2,a_3) & = \frac{1-\gamma}{2(2-\gamma)}\left[ (2a_2)^\frac{2-\gamma}{1-\gamma} - (2a_1)^\frac{2-\gamma}{1-\gamma} + 1 - \left([2a_3-1]^+\right)^\frac{2-\gamma}{1-\gamma} \right],\\
	A(t,a_1,a_2,a_3) & = \left( \frac{\phi(t)}{k(t)}  \right)^\frac{1}{n-\gamma} \ell(a_1,a_2,a_3)^\frac{1-\gamma}{n-\gamma}.
	\end{align*}
It can be proved that $a_3^\star\geq\frac12+a_2^\star$, $a_1^\star\leq x_h$ so the optimisation over $(a_1,a_2,a_3)\in\mathcal{A}_3$ reduces to optimising over the set $D_h^\prime = [0,x_h]\times[x_h,1]\times[\max\{\frac{1}{2},x_h\},1]$. Then, similar computations to the ones in Theorems \ref{th:main2} and \ref{t:exam2} give
\begin{align*}
	p^\star (t,x) =\begin{cases} 
	\displaystyle\tilde p^\star(t,x), \text{ if }x\in[0^\star,a_1^\star),\\[0.8em]
	\displaystyle \frac{H(a_2^\star)}{T} -N_\gamma  \left( (a_2^\star-\mu)^{\frac{1}{1-\gamma}}-(x-\mu)^{\frac{1}{1-\gamma}}\right), \ \text{if}\ x\in[a_1^\star,a_2^\star),\\[0.8em]
	\displaystyle\tilde p^\star(t,x), \text{ if }x\in[a_2^\star,a_3^\star],\\[0.8em]
	\displaystyle  \frac{H(a_3^\star)}{T} +N_\gamma \left( \left(x-\frac12\right)^{\frac{1}{1-\gamma}}-\left(a_3^\star -\frac12\right)^{\frac{1}{1-\gamma}}\right), \ \text{if}\ x\in(a_3^\star,1],
	\end{cases}
	\end{align*}
where the constant $\mu$ is such that $p^\star(t,a_1^\star)=\frac{H(a_1^\star)}{T}$. Since the map $x\longmapsto x\phi(t)c^\gamma/\gamma-p^\star(t,x)$ is concave on $[0,1]$, its maximum is attained at
\[x^\star(c):=\begin{cases}
\displaystyle1,\ \text{if}\ c> L_\gamma(t)2^{-\frac{1}{1-\gamma}},\\[0.8em]
\displaystyle\frac12+L_\gamma(t)^{\gamma-1}c^{1-\gamma},\ \text{if}\ L_\gamma(t)\left(a_3^\star-\frac12\right)^{\frac{1}{1-\gamma}}<c\leq L_\gamma(t)2^{-\frac{1}{1-\gamma}},\\[0.8em]
\displaystyle\tilde x^\star(c), \ \text{if}\ L_\gamma(t)(a_2^\star-\mu)^{\frac{1}{1-\gamma}}<c\leq L_\gamma(t)\left(a_3^\star-\frac12\right)^{\frac{1}{1-\gamma}},\\[0.8em]
\mu+L_\gamma(t)^{\gamma-1}c^{1-\gamma},\ \text{if}\ 0\leq c\leq L_\gamma(t)(a_2^\star-\mu)^{\frac{1}{1-\gamma}},
\end{cases}\]
where $\tilde x^\star(c)$ is any point in $[a_2^\star,a_3^\star]$ such that
\[\frac{\partial \tilde p^\star}{\partial x}(t,\tilde x^\star(c))=\phi(t)\frac{c^\gamma}{\gamma}.\]
We deduce that
\begin{align*}
p(t,c)=\begin{cases}
\displaystyle\phi(t)\frac{c^\gamma}{\gamma}-N_\gamma\left(2^{-\frac{1}{1-\gamma}}-\left(a_3^\star-\frac12\right)^{\frac{1}{1-\gamma}}\right)-\frac{H(a_3^\star)}{T},\ \text{if}\ c> L_\gamma(t)2^{-\frac{1}{1-\gamma}},\\[0.8em]
\displaystyle\phi(t)\frac{c^\gamma}{2\gamma}+\phi(t)L_\gamma(t)^{\gamma-1}c+ N_\gamma\left(a_3^\star-\frac12\right)^{\frac{1}{1-\gamma}}-\frac{H(a_3^\star)}{T},\ \text{if}\ L_\gamma(t)\left(a_3^\star-\frac12\right)^{\frac{1}{1-\gamma}}<c\leq \frac{L_\gamma(t)}{2^{\frac{1}{1-\gamma}}},\\[0.8em]
\displaystyle\tilde x^\star(c)\phi(t)\frac{c^\gamma}{\gamma}-\tilde p^\star(t,\tilde x(c)), \ \text{if}\ L_\gamma(t)(a_2^\star-\mu)^{\frac{1}{1-\gamma}}<c\leq L_\gamma(t)\left(a_3^\star-\frac12\right)^{\frac{1}{1-\gamma}},\\[0.8em]
\phi(t)L_\gamma(t)^{\gamma-1}c+\mu\phi(t)\frac{c^\gamma}{\gamma}-\frac{H(a_2^\star)}{T}+N_\gamma(a_2^\star-\mu)^{\frac{1}{1-\gamma}},\ \text{if}\ 0\leq c\leq L_\gamma(t)(a_2^\star-\mu)^{\frac{1}{1-\gamma}}.
\end{cases}
\end{align*}

\end{example}}

\section{Conclusion}\label{sec:conc}

In this paper, we provide an explicit formulation of optimal electricity tariffs, from the point of view of a power company offering service to consumers facing alternatives for satisfying energy requirements. It covers the adverse selection feature of Agents, which has been empirically highlighted in the literature \cite{ qiu2017risk}. From this perspective, our results differ from the existing analysis in the literature on electricity pricing, which mainly focused on the point of view of a monopoly needing to recover its costs, such as in the Ramsey--Boiteux pricing described in \cite{boyer2006partage} or in \cite{rasanen1997optimal}. Our main result is the derivation of time--dependent and consumption--based optimal tariffs. They act in several ways in order to reflect the increasing marginal cost of electricity production with the aggregate consumption, by making electricity more expensive in peak hours; increasing the marginal price of high consumption power units; excluding effectively some of the consumers, by offering a tariff which is too expensive in comparison to their alternative sources of energy.

\medskip
Our tariff is either linear or concave. Thus, its structure can compare to common electricity tariffs which are often linear functions of the consumed energy (see for example the Australia Queensland tariffs\footnote{Available at \url{https://www.dews.qld.gov.au/electricity/prices/tariffs}}). Interestingly, it can also compare to regulated tariff menus\footnote{In practice, consumers can choose for instance: between having a standing charge or not; or having a duration of engagement or not.} that are offered to consumers. For instance, residential consumers in France can choose their maximal consumption which is associated to a different fixed cost plus a common linear pricing, these family of functions could potentially be used as an approximation to our concave tariffs.\footnote{Graphically, the French menu looks like a set of parallel segments with different lengths that could be linked to concave functions by some appropriate regression. For a description see \myself.} The main distinction is the concavity of our optimal tariff with respect to the power limit, and this difference may be induced by the fact that we offer a single tariff instead of a menu of contracts. By doing so, our tariff intends to address all consumers in a single function. As a future research project, we plan to tackle the extension of our model to a menu of tariffs, allowing more flexibility for handling adverse selection. 

\medskip
This paper is but the first step in the direction we have outlined. In addition to extending the analysis to a menu of tariffs, at this stage, there is no uncertainty in our model, and Agents commit for a definite period $[ 0, T] $. At a later stage, we plan to
incorporate uncertainties in the model: there are fluctuations in the prices of primary energies and changes in uses of electricity, and domestic consumption is subject to weather conditions. Expanding the characteristics of Agents in order to model their ability to modulate their consumption is also a possible extension; see for example \cite{joskow2006retail} for an analysis of rational consumers who react imperfectly to prices or \cite{salies2013real} for a model with a portion of consumers who cannot modify their consumption. One can also consider more sophisticated contracts, allowing for instance the producer to cut supply a certain number of days during the year. One should also take into account
the existence of a spot market for electricity, which creates arbitrage opportunities for the Agents. Most tantalising are the congestion problems: peak hours are not seasonal, they result from aggregate behaviour, which is strategic in nature, and call for an analysis in terms of mean field games.

{\footnotesize
\bibliographystyle{plain}
\bibliography{bibliographyDylan}
}

\appendix

\section{Some convex analysis}\label{Appendix A}

We first recall the definition of $u-$convexity (adapted to our context, we refer the reader to the monograph by Villani \cite{villani2008optimal} on optimal transport theory for more details).

\begin{definition}
	Let $\psi$ be a map from $[0,T]\times X$ to $\R$. The $u-$transform of $\psi$, denoted by $\psi^\star:[0,T]\times \Cc\longrightarrow \R\cup\{+\infty\}$ is defined by
	\[\psi^\star(t,c):=\underset{x\in X}{\sup}\left\{u(t,x,c)-\psi(t,x)\right\}, \ \text{for any $(t,c)\in[0,T]\times\Cc$}.\]
	Similarly, if $\varphi$ is a map from $[0,T]\times\Cc$ to $\R$, its $u-$transform, still denoted by $\varphi^\star:[0,T]\times X\longrightarrow \R\cup\{+\infty\}$, is defined by
	\[\varphi^\star(t,x):=\underset{c\in \Cc}{\sup}\left\{u(t,x,c)-\varphi(t,c)\right\}, \ \text{for any $(t,x)\in[0,T]\times X$}.\]
	A map $\phi:[0,T]\times\Cc\longrightarrow \R\cup\{+\infty\}$ is then said to be $u-$convex if it is proper\footnote{That is to say not identically equal to $+\infty$.}and if there exists some $\psi:[0,T]\times X\longrightarrow \R$ such that
	\[\phi(t,c)=\psi^\star (t,c),\ \text{for any $(t,c)\in[0,T]\times\Cc$}.\]
	Similarly, a map $\Phi:[0,T]\times X\longrightarrow \R\cup\{+\infty\}$ is said to be $u-$convex if it is proper and there exists some $\Psi:[0,T]\times \Cc\longrightarrow \R$ such that
	\[\Phi(t,x)=\Psi^\star (t,x),\ \text{for any $(t,x)\in[0,T]\times X$}.\]
\end{definition}

We recall the following easy characterisation of $u-$convexity.
\begin{lemma}\label{lemma:uconvexity}
	A map $\phi:[0,T]\times\Cc\longrightarrow \R\cup\{+\infty\}$ is $u-$convex if and only if
	\[\phi(t,c)=(\phi^\star)^\star(t,c),\ \text{for any $(t,c)\in[0,T]\times\Cc$}.\]
	A similar statement holds for maps $\Phi:[0,T]\times X\longrightarrow \R\cup\{+\infty\}$.
\end{lemma}

\begin{proof}
We only prove the first statement, the other one being exactly similar. The result is an easy consequence of the fact that for any map $\phi:[0,T]\times\Cc\longrightarrow \R\cup\{+\infty\}$, we have the identity
\[\phi^\star(t,x)=((\phi^\star)^\star)^\star(t,x),\ \text{for any $(t,x)\in[0,T]\times X$}.\]
Indeed, we have by definition that for any $(t,x)\in[0,T]\times X$
\begin{align*}
((\phi^\star)^\star)^\star(t,x)&=\underset{c\in\Cc}{\sup}\left\{u(t,x,c)-\underset{x'\in X}{\sup}\left\{u(t,x',c)-\underset{c'\in\Cc}\sup\left\{u(t,x',c')-\phi(t,c'))\right\}\right\}\right\}\\
&=\underset{c\in\Cc}{\sup}\ \underset{x'\in X}{\inf}\ \underset{c'\in\Cc}\sup\left\{u(t,x,c)-u(t,x',c)+u(t,x',c')-\phi(t,c')\right\}.
\end{align*}
Choosing $x'=x$, we immediately get that $((\phi^\star)^\star)^\star(t,x)\leq \phi^\star(t,x)$, while the converse inequality is obtained by choosing $c=c'$.
\end{proof}

\medskip
Next, we can define the notion of the $u-$subdifferential  of a $u-$convex function.

\begin{definition}
	Let $\phi:[0,T]\times\Cc\longrightarrow \R\cup\{+\infty\}$ be a $u-$convex function. For any $(t,c)\in[0,T]\times\Cc$, the $u-$subdifferential of $\phi$ at the point $(t,c)$ is the set $\partial^\star\phi(t,c)\subset X$ defined by
	\[\partial^\star\phi(t,c):=\left\{x\in X,\ \phi^\star(t,x)=u(t,x,c)-\phi(t,c)\right\}.\]
	Similarly, let $\psi:[0,T]\times X\longrightarrow \R\cup\{+\infty\}$ be a $u-$convex function. For any $(t,x)\in[0,T]\times X$, the $u-$subdifferential of $\psi$ at the point $(t,x)$ is the set $\partial^\star\psi(t,x)\subset \R_+$ defined by
	\[\partial^\star\psi(t,x):=\left\{c\in \Cc,\ \psi^\star(t,c)=u(t,x,c)-\psi(t,x)\right\}.\]
\end{definition}
Since the map $u$ is continuous, a $u-$convex function is automatically lower-semicontinuous and its $u-$subdifferential is a closed set.

\section{Proofs of Section \ref{sec:model}}

{\bf Proof of Proposition \ref{prop:agent-problem}.} 
Since the space of admissible strategies for the Agent is decomposable and the integrand is normal when $p$ is admissible (see Definitions 14.59 and 14.27 in Rockafellar and Wets \cite{rockafellar2009variational} and also the particular case 14.29 of a Careth\'eodory integrand), we have from Theorem 14.60 in \cite{rockafellar2009variational} that the solution of problem \eqref{prob:Agent} is given by pointwise optimisation. Moreover, $\partial^\star p^\star(t,x)$ is non--empty for every $(t,x)\in[0,T]\times X(p)$, so we have that every optimal consumption strategy $c^\star:[0,T]\longrightarrow\R_+$ satisfies $c^\star(t)\in\partial^\star p^\star(t,x)$ for almost every $t\in[0,T]$ and
\[p^{\star}(t,x) =  u(t,x,c^\star(t)) - p(t,c^\star(t)).
\]
Since $u(t,x,0)$ does not depend on $x$, the envelop Theorem ensures that the map $x\longmapsto p^\star(t,x)$ is differentiable Lebesgue almost everywhere and that we have for almost every $(t,x)\in [0,T]\times X(p)$
\begin{equation} \label{envelope}
\frac{\partial u}{\partial x}(t,x,c^\star (t)) = \frac{\partial p^\star }{\partial x}(t,x).
\end{equation}
Indeed, if $c^\star(t)>0$, that is the classical envelop Theorem. Otherwise, when $\Cc=\R_+$, it is immediate to check, using the fact that $u(t,x,0)$ does not depend on $x$, that for any $(t,x)\in[0,T]\times X$, we have
\[
0\in \partial^\star p^\star(t,x) \Longrightarrow 0\in \partial^\star p^\star(t,x'), \text{ for all } x' \leq x,
\]
so that both terms in \eqref{envelope} are then actually equal to $0$.

\medskip
Then, since the map $c\longmapsto  \frac{\partial u}{\partial x}(t,x,c)$ is invertible, we have for almost every $(t,x)\in [0,T]\times X(p)$ that $\partial^\star p^\star(t,x)$ is a singleton, and the optimal consumption is $c^\star:[0,T]\times X(p)\longrightarrow\R_+$ defined in \eqref{optim_conso}.
\qed

\section{Proofs of Section \ref{section-constant}}\label{Appendix B}

Let us start this section with the following Lemma, which provides sufficient conditions of $u-$convexity when Assumption \ref{assump:power} holds. 

\begin{lemma}\label{lemma:sion}
	Let Assumption \ref{assump:power} hold and suppose in addition that $g_\gamma$ is concave if $\gamma\in(0,1)$ and convex if $\gamma\in(-\infty,0)$. Let $\psi:[0,T]\times X\longrightarrow \R$ be a map such that $x\longmapsto \psi(t,x)$ is non-decreasing and convex. Then $\psi$ is $u-$convex. Furthermore, if we take
	\[g_\gamma(x):=\begin{cases}
	x,\; \text{if $\gamma\in(0,1)$},\\
	1-x,\; \text{if $\gamma<0$},
	\end{cases}\]
	then any $u-$convex function $[0,T]\times X\longrightarrow \R$ is convex.
\end{lemma}

{\bf Proof of Lemma \ref{lemma:sion}.} By Lemma \ref{lemma:uconvexity}, we know that the $u-$convexity of $\psi$ is equivalent to
\begin{equation}\label{eq:sion}
\psi(t,x)=\sup_{c>0}\min_{y\in[0,1]}\left\{(g_\gamma(x)-g_\gamma(y))\phi(t)\frac{c^\gamma}{\gamma}+\psi(t,y)\right\}.
\end{equation}
First notice that since $\psi$ is convex in $y$ and $\frac{g_\gamma}{\gamma}$ is concave, then for any $(t,x,c)\in[0,T]\times[0,1]\times(0,+\infty)$, the map
\[f_{(t,x)}(y,c):=(g_\gamma(x)-g_\gamma(y))\phi(t)\frac{c^\gamma}{\gamma}+\psi(t,y),\]
is convex in $y$. Furthermore, for any $(t,x,y)\in[0,T]\times[0,1]^2$, the map $c\longmapsto f_{(t,x)}(y,c)$ is monotone on $(0,+\infty)$ and therefore quasiconcave. Since $[0,1]$ is convex and compact, we can apply Sion's minimax theorem \cite{sion1958on} to obtain that
\begin{align*}
\sup_{c>0}\min_{y\in[0,1]}\left\{(g_\gamma(x)-g_\gamma(y))\phi(t)\frac{c^\gamma}{\gamma}+\psi(t,y)\right\}&=\min_{y\in[0,1]}\sup_{c>0}\left\{(g_\gamma(x)-g_\gamma(y))\phi(t)\frac{c^\gamma}{\gamma}+\psi(t,y)\right\}\\
&=\min_{y\in[0,1]}\left\{+\infty{\bf 1}_{x>y}+\psi(t,y)\right\}=\psi(t,x),
\end{align*}
since $\psi$ is non-decreasing.

\medskip
Finally, it is easy to see that when $g_\gamma$ is defined as in the statement of the lemma, we have
\[\sup_{c>0}\min_{y\in[0,1]}\left\{(g_\gamma(x)-g_\gamma(y))\phi(t)\frac{c^\gamma}{\gamma}+\psi(t,y)\right\}=\sup_{c>0}\min_{y\in[0,1]}\left\{(x-y)c+\psi(t,y)\right\},\]
which corresponds to the classical convex conjugate, hence the desired result by Fenchel-Moreau's theorem.
\qed

\medskip
{\bf Proof of Theorem \ref{th:main}.} We optimise first with respect to $p^\star$. For fixed $x_0$, start by defining $\Psi_{x_0}: L^1([0,T]\times[0,1])\longrightarrow \R$ by
\begin{align}\label{psi}
&\nonumber\Psi_{x_0}(p^\star):= \int_{0}^{T} \bigg[ \int_{x_0}^1  \frac{\left( g_\gamma(x)f(x)+g'_\gamma(x)F(x)-g'_\gamma(x) \right)}{g'_\gamma(x)} \frac{\partial p^\star}{\partial x}(t,x)  \mathrm{d}x  \\ 
&\hspace{6em}- K\bigg( t,\int_{x_0}^1  \left( \frac{\gamma}{\phi(t)g'_\gamma(x)} \frac{\partial p^\star }{\partial x} (t,x) \right)^{\frac{1}{\gamma}}f(x)\mathrm{d}x \bigg)  \bigg] \mathrm{d}t +(F(x_0)-1)H.
\end{align}

\medskip

$\Psi_{x_0}$ is clearly continuous and Fr\'echet differentiable and it is also concave because $K$ is convex in $c$. Furthermore, for any $q\in L^1([0,T]\times[0,1])$, we have
\begin{align*}
\Psi'_{x_0}(p^\star;q)=&\ \int_0^T\int_{x_0}^1 \frac{\left( g_\gamma(x)f(x)+g'_\gamma(x)F(x)-g'_\gamma(x) \right)}{g'_\gamma(x)} \frac{\partial q}{\partial x}(t,x)  \mathrm{d}x\mathrm{d}t\\
&- \int_0^T\int_{x_0}^1  \frac{\partial q }{\partial x} (t,x) \frac1\gamma \left( \frac{\partial p^\star }{\partial x} (t,x) \right)^{\frac{1-\gamma}{\gamma}}  \left( \frac{\gamma}{\phi(t)g'_\gamma(x)} \right)^{\frac{1}{\gamma}} f(x)  \frac{\partial K}{\partial c}\left( t,A(t,x_0) \right) \mathrm{d}x\mathrm{d}t,
\end{align*}
where we defined 
\[
A(t,x_0) := \int_{x_0}^1  \left( \frac{\gamma}{\phi(t)g'_\gamma(x)} \frac{\partial p^\star }{\partial x} (t,x) \right)^{\frac{1}{\gamma}}  f(x) \mathrm{d}x.
\]
Since $\Psi_{x_0}$ is concave, the necessary and sufficient optimality condition for the problem with fixed $x_0$ is
\begin{equation}\label{eq:optim}
\Psi'_{x_0}(p^\star;q) \leq 0,~ \forall q\in T_{C^+(x_0)}(p^\star),
\end{equation}
where $T_{C^+(x_0)}(p^\star)$ denotes the tangent cone to the closed set $C^+(x_0)$ at the point $p^\star$ defined by
\[
T_{C^+(x_0)}(p^\star) := \left\{ z,\ \exists\varepsilon>0, \forall h\in[0,\varepsilon] ~\exists w(h)\in C^+(x_0),\  ||p^\star+hz-w(h) ||=\circ(h)  \right\}.
\]
Using local functions we see that the inequality \eqref{eq:optim} must be satisfied almost everywhere on $[0,T]\times[x_0,1]$, that is, for every $q\in T_{C^+(x_0)}(p^\star)$ and almost every $(t,x)\in[0,T]\times[x_0,1]$ we have
\begin{align*}
\frac{\partial q}{\partial x}(t,x) \bigg[ \frac{( g_\gamma(x)f(x)+g'_\gamma(x)(F(x)-1))}{g'_\gamma(x)}  - \frac1\gamma \left( \frac{\partial p^\star }{\partial x} (t,x) \right)^{\frac{1-\gamma}{\gamma}}  \left( \frac{\gamma}{\phi(t)g'_\gamma(x)} \right)^{\frac{1}{\gamma}} f(x)  \frac{\partial K}{\partial c}( t,A(t,x_0) ) \bigg] \leq 0.
\end{align*}
Therefore, the optimal $p^\star\in C^+(x_0)$ should verify
\begin{equation}\label{optimalpstar}
\frac{\partial p^\star}{\partial x}(t,x)=  \left(\frac{\phi(t)^{\frac1\gamma}\left[ g_\gamma(x)f(x)+g'_\gamma(x)F(x)-g'_\gamma(x)\right]^+}{f(x)\displaystyle\frac{\partial K}{\partial c}\left(t,A(t,x_0)\right)}\right)^{\frac{\gamma}{1-\gamma}} \frac{g'_\gamma(x)}{\gamma},
\end{equation}
with $a^+=\max\{a,0\}$, the positive part operator. By the above equation, we must have
\[
A(t,x_0)\left(\frac{\partial K}{\partial c}\left(t,A(t,x_0)\right)\right)^{\frac{1}{1-\gamma}}=\phi^{ \frac{1}{1-\gamma}}(t)\int_{x_0}^1\left(\frac{ \left[ g_\gamma(x)f(x)+g'_\gamma(x)F(x)-g'_\gamma(x) \right]^+ }{f^\gamma(x)}\right)^{\frac{1}{1-\gamma}}\mathrm{d}x.
\]
Now, let
\[
g_K(c):=c\left(\frac{\partial K}{\partial c}(t,c)\right)^{\frac{1}{1-\gamma}},\ c\geq 0.
\]
Since $K$ is strictly convex and increasing with respect to $c$, it can be checked directly that $g_K$ is increasing as well (on $\R_+$), so that we deduce
\[
A(t,x_0)=g_K^{(-1)}\left(\phi^{ \frac{1}{1-\gamma}}(t)\int_{x_0}^1\left(\frac{ \left[  g_\gamma(x)f(x)+g'_\gamma(x)F(x)-g'_\gamma(x) \right]^+ }{f^\gamma(x)}\right)^{\frac{1}{1-\gamma}}\mathrm{d}x\right).\]
Therefore the solution to the problem with fixed $x_0$ is given by \eqref{optimalpstar}. We have thus reduced the problem \eqref{solvableproblem} to
\begin{align*}
\widetilde U_P=&\sup_{ x_0\in[0,1]}  \int_{0}^{T}\left( \frac{\phi^{\frac{1}{1-\gamma}}(t)\ell(x_0)}{\gamma\left(\frac{\partial K}{\partial c}(t,A(t,x_0))\right)^{\frac{\gamma}{1-\gamma}}}- K\left( t,\frac{\phi^{\frac{1}{1-\gamma}}(t)\ell(x_0)}{\left(\frac{\partial K}{\partial c}(t,A(t,x_0))\right)^{\frac{1}{1-\gamma}}}\right)\right)\mathrm{d}t +(F(x_0)-1)H.
\end{align*}
Seen as a function of $x_0$, the right-hand side above is clearly a continuous function. It therefore attains its maximum over the compact set $[0,1]$ at some (possibly non unique) $x_0^\star$. We will abuse notations and denote by $x_0^\star$ a generic maximiser.

\medskip
If $p^\star$ is $u-$convex, since $p$ is also $u-$convex (by definition), then $p^\star$ is necessarily the $u-$transform of $p$ and therefore $p\in\Pc$, which means that we actually have $\widetilde U_P=U_P$. For the uniqueness result, define
\begin{align*}
\alpha(x_0) &: = \int_{0}^{T}\left( \frac{\phi^{\frac{1}{1-\gamma}}(t)\ell(x_0)}{\gamma\left(\frac{\partial K}{\partial c}(t,A(t,x_0))\right)^{\frac{\gamma}{1-\gamma}}}- K\left( t,\frac{\phi^{\frac{1}{1-\gamma}}(t)\ell(x_0)}{\left(\frac{\partial K}{\partial c}(t,A(t,x_0))\right)^{\frac{1}{1-\gamma}}}\right)\right)\mathrm{d}t  +(F(x_0)-1)H.
\end{align*}
Note that $\alpha$ does not attain its maximum over any interval outside $L$, because there $\ell$ is constant (and therefore $A(t,\cdot)$ too) and $F$ is increasing. Then, since over $L$ we have
\begin{align*}
\alpha'(x_0) & = \int_0^T \left(\frac{1-\gamma}{\gamma}\right)\phi(t)^\frac{1}{1-\gamma}\ell'(x_0)\left(\frac{\partial K}{\partial c}(t,A(t,x_0))\right)^{\frac{-\gamma}{1-\gamma}}\mathrm{d}t + f(x_0)H \\
& = - \int_0^T \left(\frac{1-\gamma}{\gamma}\right)\phi(t)^\frac{1}{1-\gamma}\beta(x_0)\left(\frac{\partial K}{\partial c}(t,A(t,x_0))\right)^{\frac{-\gamma}{1-\gamma}}\mathrm{d}t + f(x_0)H.
\end{align*}
Under the hypotheses of the theorem, $\alpha'$ is decreasing over $L$ in each one of the two cases so $\alpha$ is strictly concave.
\qed

\medskip
{\bf Proof of Theorem \ref{th:mainex}.} We divide the proof in two cases.

\medskip
\hspace{3em}{$\bullet$ \bf Case 1: $\gamma\in(0,1)$}

\medskip
In this case we have
\[
\ell(x_0)=\int_{x_0\vee\frac{1}{2}}^1\left(2x-1\right)^{\frac 1{1-\gamma}}\mathrm{d}x=\frac{1-\gamma}{2(2-\gamma)}\left(1-\left((2x_0-1)^+\right)^{\frac{2-\gamma}{1-\gamma}}\right).
\]
Hence, it is clear that $x_0\longmapsto \Phi(x_0)$ is increasing in $[0,\frac{1}{2}]$, so that it suffices to solve
\[
\sup_{ x_0\in[1/2,1]}\left\{ B_\gamma(T) \ell(x_0)^{\frac{n(1-\gamma)}{n-\gamma}}+(x_0-1)H\right\}.
\]
Let
\[
 y_0:=(2x_0-1)^{\frac{1}{1-\gamma}},~ A_\gamma(T) := B_\gamma(T)\left(\frac{1-\gamma}{2(2-\gamma)}\right)^{\frac{n(1-\gamma)}{n-\gamma}}.
\]
Defining the map $\overline \Phi:[0,1]\longrightarrow \R$ by $$\overline \Phi(y_0):=\Phi\left(\frac{y_0^{1-\gamma}+1}{2}\right),$$
we deduce
\[
\overline \Phi(y_0)=A_\gamma(T) \left(1-y_0^{2-\gamma}\right)^{\frac{n(1-\gamma)}{n-\gamma}}+\frac12\left(y_0^{1-\gamma}-1\right)H.
\]
Next, we can check directly that $\overline \Phi$ is concave on $[0,1]$, and we have for any $y_0\in[0,1]$
\[
\overline\Phi'(y_0)=\frac{1-\gamma}{2}y_0^{-\gamma}\left(H -2nA_\gamma(T)\frac{2-\gamma}{n-\gamma}y_0\left(1-y_0^{2-\gamma}\right)^{-\frac{\gamma(n-1)}{n-\gamma}}\right).
\]
Denote finally for any $y_0\in[0,1]$
\[
\chi(y_0):=H-2nA_\gamma(T)\frac{2-\gamma}{n-\gamma}(2x_0-1)^{\frac{1}{1-\gamma}}\left(1-y_0^{2-\gamma}\right)^{-\frac{\gamma(n-1)}{n-\gamma}}.
\]
We have for any $y_0\in[0,1]$
\[
\chi'(y_0)=-\frac{2n(2-\gamma)}{(n-\gamma)^2}A_\gamma(T)\left(1-y_0^{2-\gamma}\right)^{-\frac{n+\gamma (n-2)}{n-\gamma}}\left(n-\gamma+\gamma(n-1)(2-\gamma)y_0^{2-\gamma}\right)<0,
\]
since $\gamma\in(0,1)$. Thus, since in addition we have $$\chi(0)=H>0,\text{ and }\underset{y_0\uparrow 1}{\lim}\ \chi(y_0)=-\infty,$$
there is a unique $y_0^\star \in(0,1)$ (and thus a unique $x_0^\star \in(1/2,1)$) such that $\overline\Phi'(y^\star _0)=0$, at which the maximum of $\overline \Phi$ is attained. Finally, we can compute explicitly $p^\star(t,x)$ for any $(t,x)\in[0,T]\times [0,1]$ as
\[
p^\star (t,x) =   \frac HT +M(t) \left( ((2x-1)^+)^{\frac{1}{1-\gamma}}-(2x_0^\star -1)^{\frac{1}{1-\gamma}}\right),
\]
where we defined for simplicity
\[M(t):=\frac{1-\gamma}{2\gamma} \left( \frac{2(2-\gamma)}{1-\gamma}\right)^{\frac{\gamma(n-1)}{n-\gamma}}\left(\frac{\phi^n(t)}{k^\gamma(t)}\right)^{\frac{1}{n-\gamma}}\left(1-(2x_0^\star -1)^\frac{2-\gamma}{1-\gamma}\right)^{-\frac{\gamma(n-1)}{n-\gamma}}.\]
It can then be checked directly that for any $c\geq 0$, the map $x\longmapsto x\phi(t)c^\gamma/\gamma-p^\star(t,x)$ is concave on $[0,1]$ and attains its maximum at the point
\[x^\star(c):={\bf 1}_{c> \left(\frac{2\gamma M(t)}{(1-\gamma)\phi(t)}\right)^{\frac 1\gamma}}+\frac12\left(1+\left(\frac{(1-\gamma)\phi(t)}{2\gamma M(t)}\right)^{\frac{1-\gamma}{\gamma}}c^{1-\gamma}\right){\bf 1}_{c\leq \left(\frac{2\gamma M(t)}{(1-\gamma)\phi(t)}\right)^{\frac 1\gamma}}.\]
Therefore, we deduce immediately that for any $(t,c)\in[0,T]\times \R_+$
\[p(t,c)=\begin{cases}
\displaystyle\phi(t)\frac{c^\gamma}{\gamma}+M(t)\left((2x_0^\star -1)^{\frac{1}{1-\gamma}}-1\right)-\frac HT,\ \text{if $c>\left(\frac{2\gamma M(t)}{(1-\gamma)\phi(t)}\right)^{\frac 1\gamma}$},\\[0.8em]
\displaystyle\phi(t)\frac{c^\gamma}{2\gamma}+\left(\left(\frac{\phi(t)}{2}\right)^{\frac{1}{1-\gamma}}\frac{1-\gamma}{\gamma M(t)}\right)^{\frac{1-\gamma}{\gamma}}c-\frac HT+M(t)(2x_0^\star -1)^{\frac{1}{1-\gamma}},\ \text{otherwise.}
\end{cases}\]
Next, we notice that for any $(t,x)\in[0,T]\times[0,1]$, the map $c\longmapsto x\phi(t)c^\gamma/\gamma -p(t,c)$ is decreasing on $\R_+$ if $x<1/2$, and that it is concave on $\R_+$ if $x\geq 1/2$, so that it attains its maximum at the point
\[c^\star(t,x):=\left(\frac{2\gamma M(t)}{(1-\gamma)\phi(t)}\right)^{\frac 1\gamma}(2x-1)^{\frac{1}{1-\gamma}}{\bf 1}_{x\in(1/2,1]}.\]
It is also immediate that $p^\star$ is always non-decreasing and is convex, and therefore $u-$convex by Lemma \ref{lemma:sion}, so much so that we conclude that $p\in\Pc$. 

\medskip
It can easily be shown that the following suboptimal but simpler tariff will give the same results in terms of selected Agents, optimal consumption and Principal's utility. Indeed as no consumer select $c>\left(\frac{2\gamma M(t)}{(1-\gamma)\phi(t)}\right)^{\frac 1\gamma}$, to replace the relative tariff section by a higher function which means a more expensive contract does not modify the results. For any $(t,c)\in[0,T]\times \R_+$, the following tarif is also admissible
\[p(t,c)= \displaystyle\phi(t)\frac{c^\gamma}{2\gamma}+\left(\left(\frac{\phi(t)}{2}\right)^{\frac{1}{1-\gamma}}\frac{1-\gamma}{\gamma M(t)}\right)^{\frac{1-\gamma}{\gamma}}c-\frac HT+M(t)(2x_0^\star -1)^{\frac{1}{1-\gamma}}.\]

\medskip
\hspace{3em}{$\bullet$ \bf Case 2: $\gamma\in(-\infty,0)$}

\medskip
Now, we actually have
\[\ell(x_0)= 2^{\frac{1}{1-\gamma}}\int_{x_0}^1(1-x)^{\frac1{1-\gamma}}\mathrm{d}x  =2^{\frac{1}{1-\gamma}}\left(\frac{1-\gamma}{2-\gamma}\right)(1-x_0)^{\frac{2-\gamma}{1-\gamma}}.\]

The problem to solve is now
\[
\sup_{ x_0\in[0,1]}\left\{ B_\gamma(T) \ell(x_0)^{\frac{n(1-\gamma)}{n-\gamma}}+(x_0-1)H\right\}.
\]
It can be checked directly that the above map is actually strictly concave for $x_0\in[0,1]$, and therefore that it attains its maximum at 
\[\widehat x_0^\star:=\left(1-\left(\frac{H}{B_\gamma(T)}\frac{n-\gamma}{n{(1-\gamma)}}\right)^{\frac{n-\gamma}{n(1-\gamma)+\gamma}}\left(\frac{2-\gamma}{1-\gamma}\right)^{\frac{-\gamma(n-1)}{n(1-\gamma)+\gamma}} 2^\frac{-n}{n(1-\gamma)+\gamma}\right)^+.\]

Finally, we can compute explicitly $p^\star(t,x)$ for any $(t,x)\in[0,T]\times [0,1]$ as
\[
p^\star (t,x) =   \frac HT +\widehat M(t) \left( (1-\widehat x_0^\star)^{\frac{1}{1-\gamma}}-(1-x)^{\frac{1}{1-\gamma}}\right),
\]
where we defined for simplicity
\[\widehat M(t):=-\frac{1-\gamma}{\gamma} \left( \frac{2-\gamma}{1-\gamma}\right)^{\frac{\gamma(n-1)}{n-\gamma}}\left(\frac{2^\gamma\phi^n(t)}{k^\gamma(t)}\right)^{\frac{1}{n-\gamma}}(1-\widehat x_0^\star)^{-\frac{\gamma(2-\gamma)(n-1)}{(n-\gamma)(1-\gamma)}}.\]
We deduce directly that in this case the map $x\longmapsto (1-x)\phi(t)c^\gamma/\gamma-p^\star(t,x)$ is concave, so that it attains its maximum on $[0,1]$ at
\[\widehat x^\star(c):=\left(1-\left(-\frac{\phi(t)(1-\gamma)}{\gamma\widehat M(t)}\right)^{\frac{1-\gamma}{\gamma}}c^{1-\gamma}\right)^+,\]
so that
\begin{align*}
p(t,c)=\begin{cases}
\displaystyle \phi(t)\frac{c^\gamma}{\gamma}-\frac HT-\widehat M(t)(1-\widehat x_0^\star)^{\frac{1}{1-\gamma}}+\widehat M(t),\ \text{if } c>\left(-\frac{\gamma\widehat M(t)}{\phi(t)(1-\gamma)}\right)^{\frac{1}{\gamma}},\\
\displaystyle -\gamma c\left(-\frac{\phi(t)}{\gamma}\right)^{\frac1\gamma}\left(\frac{1-\gamma}{\widehat M(t)}\right)^{\frac{1-\gamma}{\gamma}}- \frac HT-\widehat M(t)(1-\widehat x_0^\star)^{\frac{1}{1-\gamma}},\ \text{otherwise.} 
\end{cases}
\end{align*}
It is also immediate in this case that $p^\star$ is always non-decreasing and is convex, and therefore $u-$convex by Lemma \ref{lemma:sion}, so much so that we conclude that $p\in\Pc$. 

\medskip
It can easily be shown that the following suboptimal but simpler tariff will give the same results in terms of selected Agents, optimal consumption and Principal's utility because no consumer selects $c>\left(-\frac{\gamma\widehat M(t)}{\phi(t)(1-\gamma)}\right)^{\frac{1}{\gamma}}$. For any $(t,c)\in[0,T]\times \R_+$, the following tariff is also admissible
\[p(t,c)= \displaystyle -\gamma c\left(-\frac{\phi(t)}{\gamma}\right)^{\frac1\gamma}\left(\frac{1-\gamma}{\widehat M(t)}\right)^{\frac{1-\gamma}{\gamma}}- \frac HT-\widehat M(t)(1-\widehat x_0^\star)^{\frac{1}{1-\gamma}}.\]

\qed

\section{Proofs of Section \ref{sec:general}}\label{Appendix C}
\subsection{Technical results}
{\bf Proof of Proposition \ref{prop:P=H}.} 
Define the following functionals on $\widehat{C}^+$
\begin{align*}
&\mathfrak K(p) := \int_0^T K\left(t, \int_{X^\star (p)}  \left( \frac{\gamma}{\phi (t)g'(x)} \frac{\partial p}{\partial x} (t,x) \right)^{\frac{1}{\gamma}} f(x)\mathrm{d}x\right) \mathrm{d}t,\\
&\mathfrak J(p):= \int_{0}^{T}  \left[ \int_{X^\star (p)}  \left( \frac{g(x)}{g'(x)}\frac{\partial p}{\partial x} (t,x)  - p(t,x)\right)   f(x) \mathrm{d}x \right] \mathrm{d}t.
\end{align*}
Proposition \ref{prop:P=H} states that if $Y^\star (p^\star )$ has positive measure, then $p^\star $ is not a maximiser of $p\longmapsto \mathfrak J(p)-\mathfrak K(p)$ over the set $\widehat{C}^+$. Indeed, in this case we can find an interval $[c,d]\subset Y^\star (p^\star )$ (remember that this is an open set with positive Lebesgue measure and that the latter is regular) and thus
\[
\int_0^T p^\star (t,x)\mathrm{d}t=H(x),\ \textrm{for every } x\in[c,d].
\]
Next, define
\[
T^+=\left\{ t\in[0,T]: p^\star(t,c) < p^\star(t,d)  \right\}.
\]
Since $H$ is strictly increasing we have that $T^+$ has positive Lebesgue measure. For every $t\in T^+$, define over $[c,d]$ a continuous and increasing function $q$ satisfying $q(t,c):=p^\star (t,c)$, $q(t,d):=p^\star (t,d)$ and $q(t,x)<p^\star (t,x)$ over $(c,d)$. Consider the following modification of $p^\star $.
\[
\hat{p}(t,x) : = 
\begin{cases}
q(t,x), \ \text{if}\ (t,x)\in T^+\times[c,d], \\
p^\star (t,x), \ \text{if}\ (t,x)\not\in T^+\times[c,d].
\end{cases}
\]
We have that $X^\star (\hat{p})=X^\star (p^\star ) \setminus (c,d)$ and therefore $\mathfrak K(\hat{p})<\mathfrak K(p^\star )$. Moreover, 
\begin{align*}
\mathfrak J(\hat{p}) & = \mathfrak J(p^\star ) - \int_{0}^{T}  \left[ \int_c^d  \left( \frac{g(x)}{g'(x)}\frac{\partial p^\star}{\partial x} (t,x)  - p^\star(t,x)\right)   f(x) \mathrm{d}x \right] \mathrm{d}t  \\
& = \mathfrak J(p^\star ) - \int_c^d  \left( \frac{g(x)}{g'(x)}H'(x)  - H(x)\right)   f(x) \mathrm{d}x   > \mathfrak J(p^\star ),
\end{align*}
where we used Assumption \ref{assump:ghf}. Since $\hat{p}$ is also non-decreasing in $x$, $\hat p\in \widehat{C}^+$, and we conclude that $p^\star $ is not optimal.
\qed

\vspace{2em}
{\bf Proof of Proposition \ref{prop:existence}.} Note that in the optimisation problem $(P_{a,b})$ we can without loss of generality restrict our attention to the feasible maps on $[0,T]\times X^\star(a,b)$. In other words, for fixed $(a,b)\in\Ac$, we define the closed and convex set $F_{a,b}$ as the set of maps $q\in W_x^{1.m}(X^\star(a,b))$ such that for every $t\in[0,T]$, $x\longmapsto q(t,x)$ is continuous and non--decreasing, $Q(x):=\int_0^T q(t,x)\mathrm{d}t\geq H(x)$ for every $x\in X^\star(a,b)$ and $Q(a_n)=H(a_n)$, $Q(b_n)=H(b_n)$ for every $n\geq 1$. 

\medskip
We show that $\Psi_{(a,b)}$, seen on the Banach space $\big(W_x^{1,m}(X^\star(a,b)),||\cdot ||_{m,X^\star(a,b)}\big)$, is coercive on $F_{a,b}$.

\medskip
\hspace{3em}{$\bullet$ \bf Case 1: $\gamma\in(0,1)$}

\medskip
Observe first that if $(q_n)_{n\in\N}\subset F_{a,b}$ is such that $||q_n||_{L^m([0,T]\times X^\star(a,b))}\underset{n\rightarrow\infty}{\longrightarrow}\infty$, then since for every $n\in\N$ the map $x\longmapsto\int_0^T q_n(t,x)\mathrm{d}t$ is bounded from below by $H$ on $X^\star(a,b)$, we have that
\begin{equation}\label{limit-1}
-\int_0^T \int_{X^\star(a,b)} q_n(t,x)f(x)\mathrm{d}x\mathrm{d}t \underset{n\rightarrow\infty}{\longrightarrow} -\infty.
\end{equation}
Next, define
\[A:=\int_{X^\star(a,b)}f(x)\left( \frac{\gamma}{\phi(t) g'(x)} \right)^\frac{1}{\gamma}\mathrm{d}x, \ B:=\int_0^T k(t)\mathrm{d}t.\]
From Jensen's inequality for the maps $\psi_A(x)=x^\frac{1}{\gamma m}$, $\psi_B(x)=x^\frac{n}{\gamma m}$, we have that (recall that $m$ is such that $\gamma m<1$ and that $n>1$)	
\begin{align*}
& \int_0^T K\left(t, \int_{X^\star(a,b)}  \left( \frac{\gamma}{\phi (t)g_\gamma'(x)} \frac{\partial q_n }{\partial x} (t,x) \right)^{\frac{1}{\gamma}} f(x)\mathrm{d}x\right)\mathrm{d}t \geq  \int_0^T k(t) \left( \int_{X^\star(a,b)}  \left( \frac{\gamma}{\phi (t)g_\gamma'(x)} \frac{\partial q_n }{\partial x} (t,x) \right)^{\frac{1}{\gamma}} f(x)\mathrm{d}x\right)^n\mathrm{d}t 
\\
&\geq  A^{n(1-\frac{1}{\gamma m})} \int_0^T k(t) \left( \int_{X^\star(a,b)}  \left( \frac{\gamma}{\phi (t)g_\gamma'(x)}  \right)^{\frac{1}{\gamma}}  \bigg| \frac{\partial q_n }{\partial x}  (t,x) \bigg|^m  f(x)\mathrm{d}x\right)^\frac{n}{\gamma m} \mathrm{d}t \\
&\geq  A^{n(1-\frac{1}{\gamma m})} B^{(1-\frac{n}{\gamma m})} \left( \int_0^T  \int_{X^\star(a,b)}  k(t) \left( \frac{\gamma}{\phi (t)g_\gamma'(x)}  \right)^{\frac{1}{\gamma}}  \bigg| \frac{\partial q_n }{\partial x}  (t,x) \bigg|^m  f(x)\mathrm{d}x\mathrm{d}t \right)^\frac{n}{\gamma m} \\
&\geq   A^{n(1-\frac{1}{\gamma m})} B^{(1-\frac{n}{\gamma m})} I  \bigg|\bigg| \frac{\partial q_n }{\partial x} \bigg|\bigg|_{L^m([0,T]\times X^\star(a,b))}^\frac{n}{\gamma}.
\end{align*}
We have therefore proved that, denoting by $m'$ the conjugate of $m$
\begin{align*}
\Psi_{(a,b)}(q_n) \leq &~ \bigg|\bigg|\frac{g_\gamma}{g_\gamma'} \bigg|\bigg|_{L^{m'}([0,T]\times X^\star(a,b))} \bigg|\bigg| \frac{\partial q_n }{\partial x} \bigg|\bigg|_{L^m([0,T]\times X^\star(a,b))}  - A^{n(1-\frac{1}{\gamma m})} B^{1-\frac{n}{\gamma m}} I  \bigg|\bigg| \frac{\partial q_n }{\partial x} \bigg|\bigg|_{L^m([0,T]\times X^\star(a,b))}^\frac{n}{\gamma}.
\end{align*}
This with \eqref{limit-1} implies clearly that $\Psi_{(a,b)}$ is indeed coercive in $F_{a,b}$.

\medskip

\hspace{3em}{$\bullet$ \bf Case 2: $\gamma<0$}

\medskip
In this case $\frac{g_\gamma}{g'_\gamma}<0$. Since $\frac{\partial q_n }{\partial x}$ is non-negative, if $|| \frac{\partial q_n }{\partial x} ||_{L^m([0,T]\times X^\star(a,b))}\longrightarrow\infty$ we have that
\begin{equation*}
\int_0^T \int_{X^\star(a,b)} \frac{g_\gamma(x)}{g'_\gamma(x)	} \frac{\partial q_n}{\partial x}(t,x) f(x)\mathrm{d}x~\mathrm{d}t \longrightarrow -\infty.
\end{equation*}
Then, from \eqref{limit-1} and the positiveness of $K$ we conclude that $\Psi_{(a,b)}$ is coercive in $F_{a,b}$. \vspace{.5em}

To conclude the proof, note that $W_x^{1,m}(X^\star(a,b))$ is a reflexive Banach space, so the coercivity of $\Psi_{(a,b)}$ implies that it possesses at least a maximiser $p^\star$ in $F_{a,b}$. Therefore any $q\in W_x^{1.m}(0,1)$, continuous and non-decreasing in $x$, which coincides with $p^\star$ in $X^\star(a,b)$ and such that $Q(x)=\int_0^T q(t,x)\mathrm{d}t$ satisfies $Q(x)<H(x)$ for $x\in[0,1]\setminus X^\star(a,b)$ is a solution to $(P_{a,b})$. 
\qed

\medskip
{\bf Proof of Proposition \ref{prop-an-bn}.} Let $p^\star$ be the solution to problem $(P_{a,b})$. For such $n_0$, $P^\star$ cannnot be strictly greater than $H$ over $(a_{n_0},b_{n_0})$. Otherwise, by Proposition \ref{optimalconvexity} it would be a convex map, which contradicts $P^\star(a_{n_0})=H(a_{n_0})$, $P^\star(b_{n_0})=H(b_{n_0})$. Then, there exists a subset of $(a_{n_0},b_{n_0})$ of positive measure over which $P^\star=H$ so we conclude by using Proposition \ref{prop:P=H}. \qed

\medskip
We state the following Lemma before proving Proposition \ref{prop:mu_t=0}.

\begin{lemma} \label{lemma:p*-constant}
	Let $p^\star$ be a solution of $(P_{a,b})$. Define
	\begin{align*}
	& T^{b_0}= \left\{ t\in[0,T] : \frac{\partial p^\star}{\partial x}(t,x)=0, \forall x\in(0,b_0)  \right\}, \ T^{a_0}= \left\{ t\in[0,T] : \frac{\partial p^\star}{\partial x}(t,x)=0, \forall x\in(a_0,1)  \right\}.
	\end{align*}
	If $T^{a_0}$ has positive Lebesgue measure, then for every $x\in(0,b_0)$, $
	g(x)f(x) + g'(x)F(x) \leq 0.$ If $T^{b_0}$ has positive Lebesgue measure, then for every $x\in(a_0,1)$, $
	g(x)f(x) + g'(x)F(x) - g'(x) \leq 0.
	$
\end{lemma}

\begin{proof}
	We consider the case in which $T^{a_0}$ has positive Lebesgue measure. Suppose there exist $[x_1,x_2]\subset[a_0,1]$ such that for every $x\in[x_1,x_2]$
	\[
	g(x)f(x) + g'(x)F(x) - g'(x) > 0.
	\]
	Then, for any $q\in W_x^{1,m}(0,1)$ satisfying $q(t,x) = 0, \forall (t,x)\not\in T^{a_0}\times [x_1,1],$ $x\longmapsto q(t,x)$ is increasing in $[x_1,x_2],$ $\forall t\in T^{a_0},$ and $q(t,x) = q(t,x_2),$ $\forall (t,x)\in T^{a_0}\times[x_2,1]$,
	the map $p^\star+\varepsilon q$ belongs to $C^+(a,b)$ for $\varepsilon\geq 0$. Therefore $
	\Psi'_{(a,b)}(p^\star;q)\leq 0,
	$, which means
	\[
	\int_{T^{a_0}}\int_{x_1}^{x_2} \frac{\partial q}{\partial x}(t,x) \frac{f(x)g(x)+g'(x)F(x)}{g'(x)} \leq 0, 
	\]
	hence a contradiction.
\end{proof}

\vspace{2em}
{\bf Proof of Proposition \ref{prop:mu_t=0}.} Let us prove the case $I\subset(a_0,1)$, the case $I\subset(0,b_0)$ being similar. From the convexity of $P^\star$ on every interval over which $P^\star$ is strictly greater than $H$ we deduce the existence of $c_0\in[a_0,1)$ such that $P^\star(x)>H(x)$ for every $x\in(c_0,1]$ and $P^\star(x)=H(x)$ for every $x\in[a_0,c_0]$. It follows from Lemma \ref{lemma:p*-constant} that  either $T^{a_0}$ is a null set or for every $t\in T^{a_0}$ the optimality conditions from Theorem \ref{thr-optimality-condition-concave-case} hold with $\mu_t=0$. Call $T_1=[0,T]\setminus\left( \Nc \cup T^{a_0}\right)$ and define for every $t\in T_1$
\[
x_1(t)=\inf\left\{ x\in  (c_0,1): \frac{\partial p^\star}{\partial x}(t,x)>0   \right\}.
\]
We have that $p^\star(t,\cdot)$ is strictly increasing in $[x_1(t),1]$ and it is given by \eqref{eq-optimality-concave-h-a0}. Define next
\[
T_1^+: =\left\{ t\in T_1,\ \mu_t > 0 \right\}, \ T_1^- :=\left\{ t\in T_1,\ \mu_t < 0 \right\}.
\]
We will prove that $T_1^-$ and $T_1^+$ have Lebesgue measure equal to zero. Consider any map $q\in W_x^{1,m}(0,1)$ satisfying
\begin{align*}
\begin{cases}
\displaystyle q(t,x) = 0, \forall (t,x)\not\in T_1^-\times(x_1(t),1] , \\
\displaystyle x\longmapsto q(t,x) \textrm{ is increasing in } [x_1(t),1], \forall t\in T_1^-.
\end{cases}
\end{align*} 
Then $p^\star + \varepsilon q \in C^+(a,b)$ for every $\varepsilon\geq 0$, so $\Psi^\prime(p^\star;q)\leq0$. Since
\[
\Psi'(p^\star;q)= \int_{T_1^-} \int_{x_1(t)}^1 -\frac{\partial q}{\partial x}(t,x) \mu_t\mathrm{d}x\mathrm{d}t,
\]
we conclude that $T_1^-$ is a null set. Next, take any $\bar{x}\in(c_0,1)$ such that $P^\star(\bar{x})>H(\bar{x})$ and for every $t\in T_1^+$ redefine if necessary the point $x_1(t)$ in order to satisfy $x_1(t)\geq\bar{x}$. Choose $\Delta>0$ and define then $q:[0,T]\times[0,1]\longrightarrow\R$ by
\[
q(t,x) := 1_{\{ t\in T_1^+,~ x\geq x_1(t) \}} \frac{\frac{\partial p^\star}{\partial x}(t,x_1(t))(x-x_1(t))+p^\star(t,x_1(t))-p^\star(t,x)}{p^\star(t,1)- p^\star(t,x_1(t)) - \frac{\partial p^\star}{\partial x}(t,x_1(t))(1-x_1(t))}\Delta.
\]
Since $p^\star(t,\cdot)$ is convex, we have that $q$ is non-increasing, $p^\star(t,\cdot)+\varepsilon q(t,\cdot)$ is non-decreasing for $\varepsilon\sim 0$ and $p(t,1)+q(t,1)=p(t,1)-\Delta$. Therefore $p^\star+\varepsilon q\in C^+(a,b)$ for $\varepsilon\sim 0$ so $\Psi^\prime(p^\star;q)\leq0$. Since
\[
\Psi'(p^\star;q)= \int_{T_1^+} \int_{x_1(t)}^1 -\frac{\partial q}{\partial x}(t,x) \mu_t~\mathrm{d}x~\mathrm{d}t,
\]
we conclude that the set $T_1^+$ has Lebesgue measure equal to zero. 
\qed

\medskip
{\bf Proof of Proposition \ref{prop:great}.} 
We show that under the conditions of the proposition, $P^\star\equiv H$ over some subset of $(0,a_0]\cup[b_0,1)$ with positive Lebesgue measure and the result follows from Proposition \ref{prop:P=H}. Suppose not, then for almost every $t\in[0,T]$ we have 
\begin{equation*}
\frac{\partial p^\star}{\partial x}(t,x)=\begin{cases}
\displaystyle\left(\frac{\phi(t)^{\frac1\gamma}\left[ g_\gamma(x)f(x)+g'_\gamma(x)F(x)\right]^+}{f(x)\frac{\partial K}{\partial c}\left(t,A(t,a_0,b_0)\right)}\right)^{\frac{\gamma}{1-\gamma}} \frac{g'_\gamma(x)}{\gamma},\; x\in (0,b_0), \\[.6cm]
\displaystyle \left(\frac{\phi(t)^{\frac1\gamma}\left[ g_\gamma(x)f(x)+g'_\gamma(x)F(x)-g'_\gamma(x)\right]^+}{f(x)\frac{\partial K}{\partial c}\left(t,A(t,a_0,b_0)\right)}\right)^{\frac{\gamma}{1-\gamma}} \frac{g'_\gamma(x)}{\gamma},\; x\in (a_0,1).
\end{cases}
\end{equation*}
Thus either $\frac{ \partial P^\star}{\partial x}(a_0)> H'(a_0)$ or $\frac{ \partial P^\star}{\partial x}(b_0)< H'(b_0)$, which contradicts that $\hat{X^\star}(p^\star)=[0,b_0) \cup (a_0,1]$.
\qed

\medskip
{\bf Proof of Theorem \ref{t:exam2}.} 
First of all, we recall that we have a degree of freedom in choosing the map $\tilde p$ to which $p$ is equal on $[b_0^\star,a_0^\star]$, since it does not play any role in criterion that $p^\star$ maximises. Of course, if we want to be able to conclude, this map has to be $u-$convex in the end. Therefore, if we can choose it so that $p^\star$ is $C^1$ and convex in $x$, we can apply Lemma \ref{lemma:sion} and conclude that $p^\star$ is indeed $u-$convex. This can be made if and only if the derivative of $p^\star$ at $a_0^\star$ is greater or equal to the derivative of $p^\star$ at $b_0^\star$, which can be shown immediately to be equivalent to, regardless of the value of $\gamma$,
\[a_0^\star-\frac12\geq b_0^\star.\]
Furthermore, if this is not satisfied, then $p^\star$ is not convex, and we can apply the second part of Lemma \ref{lemma:sion} to conclude that $p^\star$ is not $u-$convex.

\medskip
We now divide the proof in two steps.

\medskip
\hspace{2em} {$\bullet$ \bf Case} (i): $\gamma\in(0,1)$.

\medskip
Given the discussion above, in this case the only thing we have to do is to compute $p$. Denote for simplicity
\[L_\gamma(t):=\left(\frac{\gamma N_\gamma}{(1-\gamma)\phi(t)}\right)^{\frac 1\gamma}.\]

We know that the map $x\longmapsto x\phi(t)c^\gamma/\gamma-p^\star(t,x)$ is concave on $[0,1]$. Notice as well that since $a_0^\star>1/2$, we have $1/2\geq a_0^\star-1/2\geq b_0^\star$. We can then compute its maximum and obtain directly that it is attained at
\[x^\star(c):=\begin{cases}
\displaystyle1,\ \text{if}\ c> L_\gamma(t)2^{-\frac{1}{1-\gamma}},\\[0.8em]
\displaystyle\frac12+L_\gamma(t)^{\gamma-1}c^{1-\gamma},\ \text{if}\ L_\gamma(t)\left(a_0^\star-\frac12\right)^{\frac{1}{1-\gamma}}<c\leq L_\gamma(t)2^{-\frac{1}{1-\gamma}},\\[0.8em]
\displaystyle\tilde x^\star(c), \ \text{if}\ L_\gamma(t)(b_0^\star)^{\frac{1}{1-\gamma}}<c\leq L_\gamma(t)\left(a_0^\star-\frac12\right)^{\frac{1}{1-\gamma}},\\[0.8em]
L_\gamma(t)^{\gamma-1}c^{1-\gamma},\ \text{if}\ 0\leq c\leq L_\gamma(t)(b_0^\star)^{\frac{1}{1-\gamma}},
\end{cases}\]
where $\tilde x^\star(c)$ is any point in $[b_0^\star,a_0^\star]$ such that
\[\frac{\partial \tilde p^\star}{\partial x}(t,\tilde x^\star(c))=\phi(t)\frac{c^\gamma}{\gamma}.\]
We deduce that
\begin{align*}
p(t,c)=\begin{cases}
\displaystyle\phi(t)\frac{c^\gamma}{\gamma}-N_\gamma\left(2^{-\frac{1}{1-\gamma}}-\left(a_0^\star-\frac12\right)^{\frac{1}{1-\gamma}}\right)-\frac{H(a_0^\star)}{T},\ \text{if}\ c> L_\gamma(t)2^{-\frac{1}{1-\gamma}},\\[0.8em]
\displaystyle\phi(t)\frac{c^\gamma}{2\gamma}+\phi(t)L_\gamma(t)^{\gamma-1}c+ N_\gamma\left(a_0^\star-\frac12\right)^{\frac{1}{1-\gamma}}-\frac{H(a_0^\star)}{T},\ \text{if}\ L_\gamma(t)\left(a_0^\star-\frac12\right)^{\frac{1}{1-\gamma}}<c\leq \frac{L_\gamma(t)}{2^{\frac{1}{1-\gamma}}},\\[0.8em]
\displaystyle\tilde x^\star(c)\phi(t)\frac{c^\gamma}{\gamma}-\tilde p^\star(t,\tilde x(c)), \ \text{if}\ L_\gamma(t)(b_0^\star)^{\frac{1}{1-\gamma}}<c\leq L_\gamma(t)\left(a_0^\star-\frac12\right)^{\frac{1}{1-\gamma}},\\[0.8em]
\phi(t)L_\gamma(t)^{\gamma-1}c-\frac{H(b_0^\star)}{T}+N_\gamma(b_0^\star)^{\frac{1}{1-\gamma}},\ \text{if}\ 0\leq c\leq L_\gamma(t)(b_0^\star)^{\frac{1}{1-\gamma}}.
\end{cases}
\end{align*}
As in the case $H$ constant, it can easily be shown that the following simpler tariff is also admissible and produce the same results as no consumer selects $c> L_\gamma(t)2^{-\frac{1}{1-\gamma}}$:
\begin{align*}
p(t,c)=\begin{cases}
\displaystyle\phi(t)\frac{c^\gamma}{2\gamma}+\phi(t)L_\gamma(t)^{\gamma-1}c+ N_\gamma\left(a_0^\star-\frac12\right)^{\frac{1}{1-\gamma}}-\frac{H(a_0^\star)}{T},\ \text{if}\ L_\gamma(t)\left(a_0^\star-\frac12\right)^{\frac{1}{1-\gamma}}<c,\\[0.8em]
\displaystyle\tilde x^\star(c)\phi(t)\frac{c^\gamma}{\gamma}-\tilde p^\star(t,\tilde x(c)), \ \text{if}\ L_\gamma(t)(b_0^\star)^{\frac{1}{1-\gamma}}<c\leq L_\gamma(t)\left(a_0^\star-\frac12\right)^{\frac{1}{1-\gamma}},\\[0.8em]
\phi(t)L_\gamma(t)^{\gamma-1}c-\frac{H(b_0^\star)}{T}+N_\gamma(b_0^\star)^{\frac{1}{1-\gamma}},\ \text{if}\ 0\leq c\leq L_\gamma(t)(b_0^\star)^{\frac{1}{1-\gamma}}.
\end{cases}
\end{align*}

\medskip
\hspace{2em} {$\bullet$ \bf Case} (ii): $\gamma<0$.
As in the previous case, our assumptions imply that $a_0^\star\geq 1/2$ and $1/2\geq a_0^\star-1/2\geq b_0^\star$. We can then prove that the maximum of the map $x\longmapsto x\phi(t)c^\gamma/\gamma-p^\star(t,x)$ is attained at
\[x^\star(c):=\begin{cases}
\displaystyle1-L_\gamma(t)^{\gamma-1}c^{1-\gamma},\ \text{if}\ 0<c\leq L_\gamma(t)(1-a_0^\star)^{\frac{1}{1-\gamma}},\\[0.8em]
\displaystyle\tilde x^\star(c), \ \text{if}\  L_\gamma(t)(1-a_0^\star)^{\frac{1}{1-\gamma}}<c\leq L_\gamma(t)\left(\frac12-b_0^\star\right)^{\frac{1}{1-\gamma}},\\[0.8em]
\displaystyle\frac12-L_\gamma(t)^{\gamma-1}c^{1-\gamma},\ \text{if}\ L_\gamma(t)\left(\frac12-b_0^\star\right)^{\frac{1}{1-\gamma}}< c\leq L_\gamma(t)2^{-\frac{1}{1-\gamma}},\\[0.8em]
\displaystyle 0,\ \text{if}\ c>L_\gamma(t)2^{-\frac{1}{1-\gamma}},
\end{cases}\]
where $\tilde x^\star(c)$ is any point in $[b_0^\star,a_0^\star]$ such that
\[\frac{\partial \tilde p^\star}{\partial x}(t,\tilde x^\star(c))=\phi(t)\frac{c^\gamma}{\gamma}.\]
We deduce that
\begin{align*}
p(t,c)=\begin{cases}
\displaystyle\phi(t)L_\gamma(t)^{\gamma-1}c+N_\gamma(1-a_0^\star)^{\frac{1}{1-\gamma}}-\frac{H(a_0^\star)}{T},\ \text{if}\ 0<c\leq L_\gamma(t)(1-a_0^\star)^{\frac{1}{1-\gamma}},\\[0.8em]
\displaystyle\tilde x^\star(c)\phi(t)\frac{c^\gamma}{\gamma}-\tilde p^\star(t,\tilde x(c)), \ \text{if}\  L_\gamma(t)(1-a_0^\star)^{\frac{1}{1-\gamma}}<c\leq L_\gamma(t)\left(\frac12-b_0^\star\right)^{\frac{1}{1-\gamma}},\\[0.8em]
\displaystyle \phi(t)\frac{c^\gamma}{2\gamma}+\phi(t)L_\gamma(t)^{\gamma-1}c+N_\gamma\left(\frac12-b_0^\star\right)^{\frac{1}{1-\gamma}}-\frac{H(b_0^\star)}{T},\ \text{if}\ L_\gamma(t)\left(\frac12-b_0^\star\right)^{\frac{1}{1-\gamma}}< c\leq \frac{L_\gamma(t)}{2^{\frac{1}{1-\gamma}}},\\[0.8em]
\displaystyle \phi(t)\frac{c^\gamma}{\gamma}+N_\gamma\left(\left(\frac12-b_0^\star\right)^{\frac{1}{1-\gamma}}-2^{-\frac{1}{1-\gamma}}\right)-\frac{H(b_0^\star)}{T},\ \text{if}\ c>L_\gamma(t)2^{-\frac{1}{1-\gamma}}.\end{cases}
\end{align*}

As previously, it can easily be shown that the following simpler tariff is also admissible and produce the same results as no consumer selects $c>L_\gamma(t)2^{-\frac{1}{1-\gamma}}$:
\begin{align*}
p(t,c)=\begin{cases}
\displaystyle\phi(t)L_\gamma(t)^{\gamma-1}c+N_\gamma(1-a_0^\star)^{\frac{1}{1-\gamma}}-\frac{H(a_0^\star)}{T},\ \text{if}\ 0<c\leq L_\gamma(t)(1-a_0^\star)^{\frac{1}{1-\gamma}},\\[0.8em]
\displaystyle\tilde x^\star(c)\phi(t)\frac{c^\gamma}{\gamma}-\tilde p^\star(t,\tilde x(c)), \ \text{if}\  L_\gamma(t)(1-a_0^\star)^{\frac{1}{1-\gamma}}<c\leq L_\gamma(t)\left(\frac12-b_0^\star\right)^{\frac{1}{1-\gamma}},\\[0.8em]
\displaystyle \phi(t)\frac{c^\gamma}{2\gamma}+\phi(t)L_\gamma(t)^{\gamma-1}c+N_\gamma\left(\frac12-b_0^\star\right)^{\frac{1}{1-\gamma}}-\frac{H(b_0^\star)}{T},\ \text{if}\ L_\gamma(t)\left(\frac12-b_0^\star\right)^{\frac{1}{1-\gamma}}< c.\end{cases}
\end{align*}

\qed

\subsection{Proof of Theorem \ref{thr-optimality-condition-concave-case}}\label{app:general}
We give here a series of result which once combined prove Theorem \ref{thr-optimality-condition-concave-case}. To simplify the statements, we give them in a generic set $(a_n,b_n)$, the generalisation being straightforward. The first proposition shows that the existence of the interval $I$ in the theorem allows us to localise Problem $(P_{a,b})$, and replace it by a simpler one, in which the constraint $P^\star(x)\geq H(x)$ for every $x\in X^\star(a,b)$ can be ignored.
\begin{proposition}\label{prop-localization}
	Let $p^\star$ be a solution of $(P_{a,b})$ and suppose there exists $x_1\in (a_n,b_n)$ such that $P^\star(x_1)>H(x_1)$. Then, there exists $x_0\in(a_n,b_n)$, $x_0<x_1$, such that $p^\star$ is solution to the following problem
	\begin{equation}\label{problem-x0x1}
	(P_{x_0,x_1}) \qquad \sup_{q\in C(x_0,x_1)} \Psi_{(a,b)}^{x_0,x_1,p^\star}(q),
	\end{equation}
	where
	\begin{align*}
	\Psi_{(a,b)}^{x_0,x_1,p^\star}(q) := & \int_0^T \int_{x_0}^{x_1} \frac{g_\gamma(x)f(x)+g_\gamma'(x)F(x)}{g_\gamma'(x)}\frac{\partial q}{\partial x}(t,x) \mathrm{d}x\mathrm{d}t \\
	 - \int_0^T K\left(t, \int_{x_0}^{x_1} \left( \frac{\gamma}{\phi(t)g_\gamma'(x)}\frac{\partial q}{\partial x}(t,x)\right)^\frac{1}{\gamma}f(x)\mathrm{d}x + I_{(a,b)}^{x_0,x_1}(p^\star) \right)\mathrm{d}t , \\
	I_{(a,b)}^{x_0,x_1}(p^\star) := & \int_{X^\star(a,b)\setminus(x_0,x_1)} \left( \frac{\gamma}{\phi(t)g_\gamma'(x)}\frac{\partial p^\star}{\partial x}(t,x)\right)^\frac{1}{\gamma}f(x)\mathrm{d}x,
	\end{align*}
	and $C(x_0,x_1)$ denotes the set of maps $q\in W_x^{1,m}(x_0,x_1)$ such that
	\begin{itemize}
		\item $x\longmapsto q(t,x)$ is continuous and increasing for every $t\in[0,T]\setminus\Nc(q)$.
		\item $\displaystyle p^\star(t,x_0) + \int_{x_0}^{x_1} \frac{\partial q}{\partial x}(t,x)~\mathrm{d}x = p^\star(t,x_1)$ for every $t\in[0,T]\setminus\Nc(q)$.
	\end{itemize}
\end{proposition}

\begin{proof}
	Define 
	\[
	x_0:=\inf\left\{ z\in X^\star(a,b),\ P^\star(x)\geq H(x_1) \textrm{ for every } x\in[z,x_1] \right\}.
	\]
	By continuity we have that $x_0<x_1$ and $P^\star(x_0)=H(x_1)$. Notice that the restriction of  $p^\star$ to the set $[x_0,x_1]$ belongs to $C(x_0,x_1)$. Suppose the restriction is not a solution of $(P_{x_0,x_1})$, then there exists $q^\star\in C(x_0,x_1)$ such that $\Psi_{(a,b)}^{x_0,x_1,p^\star}(q^\star)>\Psi_{(a,b)}^{x_0,x_1,p^\star}(p^\star)$. Define then $\bar{p}:[0,T]\times[0,1]\longrightarrow\R$ by
	\[
	\bar{p}(t,x):=\begin{cases}
	\displaystyle p^\star(t,x),\ x\not\in[x_0,x_1], \\
	\displaystyle p^\star(t,x_0) + \int_{x_0}^x \frac{\partial q^\star}{\partial x}(t,x)\mathrm{d}x,\ x\in(x_0,x_1).
	\end{cases}
	\]
	Then, for every $x\in[x_0,x_1]$
	\[
	\int_0^T \bar{p}(t,x)\mathrm{d}t \geq \int_0^T \bar{p}(t,x_0)\mathrm{d}t \geq H(x_1) \geq H(x),
	\]
	and it is straightforward that $\bar{p}\in C^+(a,b)$. This is a contradiction with the optimality of $p^\star$ in problem $(P_{a,b})$ because
	\[
	\Psi_{(a,b)}(\bar{p}) = \Psi_{(a,b)}(p^\star) - \Psi_{(a,b)}^{x_0,x_1,p^\star}(p^\star)+\Psi_{(a,b)}^{x_0,x_1,p^\star}(q^\star) .
	\] 
\end{proof}

Now we state the optimality conditions for the problem $(P_{x_0,x_1})$.
\begin{proposition}\label{prop-optimality-condition-problemx0x1}
	Let $p^\star$ be a solution of $(P_{x_0,x_1})$ with $x_0,$ $x_1$ as in Proposition \ref{prop-localization}. Then there exists a null set $\Nc\subset[0,T]$ and a constant $\mu_t$ for every $t\in[0,T]\setminus\Nc$ such that for every $x\in (x_0,x_1)$ 
	\begin{equation}
	\frac{\partial p^\star}{\partial x}(t,x) = \displaystyle\left(\frac{\phi(t)^{\frac1\gamma}\left[ g_\gamma(x)f(x)+g'_\gamma(x)F(x)+g'_\gamma(x)\mu_t\right]^+}{f(x)\frac{\partial K}{\partial c}\left(t,A(t,a,b)\right)}\right)^{\frac{\gamma}{1-\gamma}} \frac{g'_\gamma(x)}{\gamma}.
	\end{equation}
\end{proposition}

\begin{proof}
	Notice that the set $C(x_0,x_1)$ can be written as
	\[
	C(x_0,x_1)=\left\{ q\in W_x^{1,m}(x_0,x_1),\ g(q)\in C,~ h(q)=0  \right\},
	\]
	where $g:W_x^{1,m}(x_0,x_1)\longrightarrow L^m([0,T]\times[x_0,x_1])$ is defined by $g(q)=\frac{\partial q}{\partial x}$, where $C$ is the following convex cone $C:=\left\{ q \in L^m([0,T]\times[x_0,x_1]),\ q(t,x)\geq 0, \textrm{ a.e.} \right\}$ and  $h:W_x^{1,m}(x_0,x_1)\longrightarrow L^m([0,T])$ is defined by 
	$$h(q):=\int_{x_0}^{x_1} \frac{\partial q}{\partial x}(\cdot,x) \mathrm{d}x +p^\star(\cdot,x_0)-p^\star(\cdot,x_1).$$
	It can be checked in the same way as in Remark 5 from \cite{daniele2007infinite}, that their Assumption S is satisfied in this context. Furthermore, it is a classical result that the dual of $W_x^{1,m}(x_0,x_1)$ is $W_x^{1,m/(m-1)}(x_0,x_1)$.
	
	\medskip
	Define now the Lagrangian $L:W_x^{1,m}(x_0,x_1)\times W_x^{1,m/(m-1)}(x_0,x_1)\times L^{\frac{m}{m-1}}(0,T)\longrightarrow\R$ by
	\begin{align*}
	L(q,\lambda,\mu)  :=&\   \Psi_{(a,b)}^{x_0,x_1,p^\star}(q)+\int_0^T\int_{x_0}^{x_1}\lambda(t,x)\frac{\partial q}{\partial x}(t,x)\mathrm{d}x\mathrm{d}t\\
	& +\int_0^T \mu(t) \left( \int_{x_0}^{x_1} \frac{\partial q}{\partial x}(t,x) \mathrm{d}x +p^\star(t,x_0)-p^\star(t,x_1)  \right)\mathrm{d}t.
	\end{align*}
	Then, from Corollary 2 in \cite{donato2011infinite} it follows that there exists $\lambda\in W_x^{1,m/(m-1)}(x_0,x_1)$, $\mu\in L^m(0,T)$ such that
	\[\begin{cases}
	\displaystyle 0  = 	 \frac{ g_\gamma(x)f(x)+g_\gamma'(x)F(x) }{g_\gamma'(x)} - \frac1\gamma \left( \frac{\partial p^\star }{\partial x} (t,x) \right)^{\frac{1-\gamma}{\gamma}}  \left( \frac{\gamma}{\phi(t)g_\gamma'(x)} \right)^{\frac{1}{\gamma}} f(x)  \frac{\partial K}{\partial c}\left( t,A(t,a,b) \right)	 \\ 
	\displaystyle  \hspace{1.9em}+ \mu(t) + \lambda(t,x) , \textrm{ a.e. in } [0,T]\times[x_0,x_1],\\[0.8em]
	\displaystyle \lambda(t,x) \frac{\partial p^\star}{\partial x}(t,x) = 0, ~\lambda(t,x) \geq 0, \textrm{ a.e. in }[0,T]\times[x_0,x_1].
	\end{cases}\]
	
	\medskip
	Then, when $\frac{\partial p^\star}{\partial x}(t,x)>0$ we have that $\lambda(t,x)=0$ and
	\[
	\frac{\partial p^\star}{\partial x}(t,x) = \displaystyle\left(\frac{\phi(t)^{\frac1\gamma}\left[ g_\gamma(x)f(x)+g'_\gamma(x)F(x)+g'_\gamma(x)\mu(t)\right]}{f(x)\frac{\partial K}{\partial c}\left(t,A(t,a,b)\right)}\right)^{\frac{\gamma}{1-\gamma}} \frac{g'_\gamma(x)}{\gamma}.
	\]
	In case $\frac{\partial p^\star}{\partial x}(t,x)=0$ we have that 
	\[
	\frac{ g_\gamma(x)f(x)+g_\gamma'(x)F(x) }{g_\gamma'(x)} + \mu(t) = -\lambda(t,x) \leq 0,
	\]
	which ends the proof.
\end{proof}

\medskip
We prove finally that the map $\mu$ does not depend on $x_0,x_1$ and is the same in the interval $I=(x_\ell,x_r)$. 

\begin{proposition}
	Let $I=(x_\ell,x_r)\subset(a_n,b_n)$ be as in Theorem \ref{thr-optimality-condition-concave-case}. Then for any $x_0,$ $x_1\in I$, there exist a null set $\Nc\subset[0,T]$ and a constant $\mu_t$ for every $t\in[0,T]\setminus\Nc$ such that for every $x\in (x_0,x_1)$ 
	\eqref{eq-optimality-concave-h-an} is satisfied.
\end{proposition}

\begin{proof}
	Let $y_0:=x_1$ and define by induction for $k\geq 0$
	\[
	z_k := \inf\{ z\in(a_n,b_n),\ P^\star(x) \geq H(y_k),\ \forall x\in[z,y_k] \}, ~ y_{k+1}:=\frac{z_k+y_k}{2}.
	\]
	By continuity we have that $P^\star(z_k)=H(y_k)$, so $y_{k+1}<y_k$ and the sequence $(y_k)_k$ converges necessarily to $a_n$. We conclude by applying Proposition \ref{prop-optimality-condition-problemx0x1} to every interval $(z_k,y_k)$ and noting that these intervals overlap themselves.
\end{proof}

\subsection{Other proofs}

{\bf Proof of Proposition \ref{optimalconvexity}.} 

From Assumption \ref{assump:ghf} and Theorem \ref{thr-optimality-condition-concave-case} we have that on every interval $I$ over which $P^\star>H$, there exists a null set $\Nc\subset[0,T]$ such that for every $t\in[0,T]\setminus\Nc$, $x\longmapsto\frac{\partial p^\star}{\partial x}(t,x)$ is non-decreasing on $I$. Therefore $P^\star$ is convex on $I$ since
	\[
	\frac{\partial P^{\star}}{\partial x}(x)=\int_0^T \frac{\partial p^\star}{\partial x}(t,x)\mathrm{d}t.
	\]
	\qed

\medskip
{\bf Proof of Proposition \ref{prop-an-bn}.} 
Let $p^\star$ be the solution of problem $(P_{a,b})$. We will prove that $P^\star\equiv H$ in the interval $(a_{n_0},b_{n_0})$ and the result will follow from Proposition \ref{prop:P=H}. Suppose not, then there exists $x_0\in(a_{n_0},b_{n_0})$ such that $P^\star(x_0)>H(x_0)$ and $p^\star$ is given by \eqref{eq-optimality-concave-h-a0} in a neighbourhood around $x_0$, so $P^\star$ is increasing in that neighbourhood. By Proposition \ref{optimalconvexity} we have that $P^\star(b_{n_0})>H(b_{n_0})$, because on every interval which is contained in the set $\{ x,\ P^\star(x) \geq H(x) \}$ the convex map $P^\star$ and the strictly concave map $H$ can intersect at most at one point. This contradicts the fact that $p^\star\in C^+(a,b)$.
	\qed

\end{document}